\RequirePackage{luatex85}

\documentclass[a4paper,reqno,10.5pt, oneside]{amsart}
\usepackage{amsmath,amssymb,amsthm,geometry,graphics,color,mdwlist,mathrsfs, mathtools}
\usepackage[all]{xy}

\usepackage{tikz,tikz-cd}
\usetikzlibrary{arrows}
\usetikzlibrary{calc,patterns}
\usetikzlibrary{decorations.markings}
\usetikzlibrary{patterns,decorations.pathreplacing}
\usetikzlibrary{arrows.meta,nfold}
\usepackage{here}
\usepackage{hyperref}

\DeclareMathOperator{\image}{Im}

\DeclareMathOperator{\Coker}{Coker}

\DeclareMathOperator{\Hom}{Hom}

\DeclareMathOperator{\RHom}{\mathrm{R}\!\Hom}
\DeclareMathOperator{\cHom}{\mathcal{H}\hspace{-0.3pt}om}
\DeclareMathOperator{\End}{End}

\DeclareMathOperator{\REnd}{\mathrm{R}\!\End}

\def\lotimes{\otimes^\mathrm{L}}

\def\rR{\mathrm{R}}
\def\RcHom{\rR\cHom}
\DeclareMathOperator{\rT}{\mathrm{T}}

\DeclareMathOperator{\vect}{\mathrm{vect}}

\DeclareMathOperator{\Ext}{Ext}

\DeclareMathOperator{\gldim}{gl.dim}

\DeclareMathOperator{\md}{mod}
\renewcommand{\mod}{\md}
\DeclareMathOperator{\Mod}{Mod}

\DeclareMathOperator{\fl}{fl}

\DeclareMathOperator{\QMod}{QMod}
\DeclareMathOperator{\qmod}{qmod}

\DeclareMathOperator{\coh}{coh}
\DeclareMathOperator{\Qcoh}{Qcoh}
\DeclareMathOperator{\res}{res}

\DeclareMathOperator{\thick}{thick}
\DeclareMathOperator{\per}{per}

\DeclareMathOperator{\add}{add}

\DeclareMathOperator{\rank}{rank}

\DeclareMathOperator{\Sym}{Sym}
\DeclareMathOperator{\Spec}{Spec}
\DeclareMathOperator{\Max}{Max}
\DeclareMathOperator{\Proj}{Proj}

\DeclareMathOperator{\height}{ht}

\DeclareMathOperator{\Cl}{Cl}

\def\seg{\#}

\def\rD{\mathrm{D}}

\def\Db{\mathrm{D^b}}

\newcommand\recollement[3]{\xymatrix{{#1}\ar[r]&{#2}\ar[r]\ar@/_7pt/[l]\ar@/^7pt/[l]&{#3}\ar@/_7pt/[l]\ar@/^7pt/[l] }}

\def\C{\mathcal{C}}
\def\D{\mathcal{D}}

\renewcommand{\O}{\mathcal{O}}
\def\P{\mathcal{P}}

\renewcommand\S{\mathcal{S}}
\def\T{\mathcal{T}}

\def\X{\mathcal{X}}
\def\Y{\mathcal{Y}}

\def\rL{\mathrm{L}}

\def\Ga{\Gamma}

\def\Del{\Delta}
\def\e{\varepsilon}

\def\La{\Lambda}

\def\Si{\Sigma}
\def\bSi{\mathbf{\Sigma}}

\def\a{\alpha}

\def\b{\beta}
\def\be{\beta}
\def\ga{\gamma}

\def\s{\sigma}
\def\si{\sigma}
\def\om{\omega}

\def\Z{\mathbb{Z}}
\def\ZZ{\mathbb{Z}}

\def\CC{\mathbb{C}}
\def\RR{\mathbb{R}}
\def\AA{\mathbb{A}}
\def\PP{\mathbb{P}}
\def\LL{\mathbb{L}}

\def\p{\mathfrak{p}}

\def\m{\mathfrak{m}}
\def\n{\mathfrak{n}}

\def\op{\mathrm{op}}

\def\rsimeq{\rotatebox{-90}{$\simeq$}}
\def\xsimeq{\xrightarrow{\simeq}}
\def\ysimeq{\xleftarrow{\simeq}}

\def\romanenum{\renewcommand\labelenumi{(\roman{enumi})}\renewcommand\theenumi{\roman{enumi}}}
\def\alphenum{\renewcommand\labelenumi{(\alph{enumi})}\renewcommand\theenumi{\alph{enumi}}}
\newtheorem{Thm}{Theorem}[section]
\newtheorem{Lem}[Thm]{Lemma}
\newtheorem{Prop}[Thm]{Proposition}
\newtheorem{Cor}[Thm]{Corollary}

\newtheorem{Prop-Def}[Thm]{Proposition-Definition}
\newtheorem{Thm-Def}[Thm]{Theorem-Definition}

\theoremstyle{definition}
\newtheorem{Def}[Thm]{Definition}
\newtheorem{Ex}[Thm]{Example}

\newtheorem{Assum}[Thm]{Assumption}

\newtheorem{Setup}[Thm]{Setting}

\theoremstyle{remark}
\newtheorem{Rem}[Thm]{Remark}

\newcounter{step}

\def\disoplus{\displaystyle\bigoplus}

\def\Wedge{\textstyle\bigwedge}
\makeatletter

\@addtoreset{equation}{section}
\makeatother

\renewcommand{\cHom}{\mathcal{H}\hspace{-1pt}{\it om}}

\def\rmX{\mathrm{X}}

\def\calP{\mathcal{P}}

\def\Cl{\mathrm{Cl}}

\newcommand{\bfSigma}{{\bf\Sigma}}

\newcommand{\calX}{\mathcal{X}}
\def\XX{\mathcal{X}}
\def\YY{\mathcal{Y}}
\newcommand{\Pic}{\operatorname{Pic}}
\newcommand{\DG}{\operatorname{DG}}

\newcommand{\calT}{\mathcal{T}}

\def\bSi{\mathbf{\Sigma}}
\def\rR{\mathrm{R}}
\def\RcHom{\rR\cHom}
\renewcommand{\cHom}{\mathcal{H}\hspace{-1pt}{\it om}}
\def\romanenum{\renewcommand\labelenumi{(\roman{enumi})}\renewcommand\theenumi{\roman{enumi}}}
\def\alphenum{\renewcommand\labelenumi{(\alph{enumi})}\renewcommand\theenumi{\alph{enumi}}}

\title{Non-commutative crepant resolutions via \linebreak higher rank vector bundles}
\author[N. Hanihara \and Y. Nakajima]{Norihiro Hanihara \and Yusuke Nakajima} 

\address[N. Hanihara]{Faculty of Mathematics, Kyushu University, 744 Motooka, Nishi-ku, Fukuoka, 819-0395, Japan}
\email{hanihara@math.kyushu-u.ac.jp}

\address[Y. Nakajima]{Department of Mathematics, Kyoto Sangyo University, Motoyama, Kamigamo, Kita-Ku, Kyoto, 603-8555, Japan}
\email{ynakaji@cc.kyoto-su.ac.jp}


\date{}

\begin{document}
\begin{abstract}
We give a construction of non-commutative crepant resolutions (NCCRs) using tilting bundles. Given a smooth projective variety $X$ and a vector bundle $E$ whose determinant is the canonical bundle $\om_X$, we show that the total space of $E$ is derived equivalent to an Adams graded Calabi--Yau dg algebra $\Ga$ whose Gorenstein parameter is equal to the rank of $E$. Consequently, it is a non-commutative crepant resolution of each of its idempotent subalgebras when the dg algebra is an ordinary algebra.
We next show that $\Ga$ can be regarded as the Calabi--Yau completion of a root pair, in the sense of the first author, over the derived category of the projective bundle associated with $E$. This generalizes a well-known fact for the case $E=\om_X$ that taking the total space of $\om_X$ amounts to the ordinary Calabi--Yau completion in the sense of Keller.
We also give a suitable variant for vector bundles over toric Deligne--Mumford stacks.
We apply this construction to some product of projective spaces, Hirzebruch surfaces, and their stack variants including Geigle--Lenzing projective spaces. We consequently obtain NCCRs for toric rings including some Hibi rings and toric compound Du Val singularities. 
\end{abstract}

\maketitle
\setcounter{tocdepth}{1}
\tableofcontents

\section{Introduction}
\subsection{Motivations and results}
A resolution of singularities is a powerful method for studying singularities. 
Among them, crepant resolutions are those that satisfy particularly good properties. 
A resolution of singularities $f\colon Y\to X$ of a normal Gorenstein singularity $X=\Spec R$ is said to be a \emph{crepant resolution} 
if $f^*\omega_X=\omega_Y$ holds, where $\omega_X$ and $\omega_Y$ are the canonical divisors of $X$ and $Y$, respectively. 
Crepant resolutions do not exist in general, and even when they exist, they are not necessarily unique. 
It is known that when there are some projective crepant resolutions, they are connected by an operation called a flop \cite[Theorem~1]{Kaw_flop}. 
Furthermore, motivated by observations in superstring theory and mirror symmetry, 
it was conjectured in \cite{BO_conj,Kaw_Dequiv} that crepant resolutions are derived equivalent. 
This conjecture was resolved in the three-dimensional case in \cite{Bri_flop} (see also \cite{Chen_flop}). 
As an alternative approach, Van den Bergh considered a non-commutative algebra derived equivalent to a crepant resolution, and resolved this conjecture via the derived category of this non-commutative algebra \cite{VdB_3dimflop}. 
In Van den Bergh's method, a non-commutative algebra obtained as the endomorphism ring of a tilting bundle plays a crucial role. 
By formulating the properties satisfied by such non-commutative algebras, the notion of a non-commutative crepant resolution was defined in \cite{VdB_NCCR} (see also \cite{VdB_ICM}). 

\begin{Def}[\cite{VdB_NCCR}]
Let $R$ be a Gorenstein normal domain and $M$ be a non-zero reflexive $R$-module. 
The $R$-algebra $\Gamma\coloneqq\End_R(M)$ is said to be a \emph{non-commutative crepant resolution} (NCCR) of $R$ 
if the global dimension of $\Gamma$ is finite and $\Gamma$ is a maximal Cohen--Macaulay $R$-module. 
\end{Def}

A typical example of NCCRs is the one of a quotient singularity for a finite group which is isomorphic to the skew group algebra (see \cite[Example~1.1]{VdB_NCCR}, \cite[Proposition~5.6]{IW_MM}). 
For several important singularities, the existence of NCCRs and their derived equivalence to crepant resolutions have been established (see e.g., \cite{BKR,IU_special,VdB_NCCR}). 
This makes it possible to study crepant resolutions from an algebraic perspective via derived equivalence \cite{IW_Q, Wem_MMP}. 
Moreover, even when a crepant resolution does not exist, an NCCR may exist. 
In such cases, one can approach a singularity algebraically by treating an NCCR as a substitute for a crepant resolution. 
For example, NCCRs are deeply connected to Calabi--Yau algebras \cite{Gi,Ke_CYcomp,IR_FZ}, (cluster) tilting theory, and higher Auslander--Reiten theory \cite{AIR,Iya_HAR,Iya_AC}, which are topics actively studied in representation theory and non-commutative algebraic geometry. Thus, singularities are linked to these theories via NCCRs. Furthermore, they are also important from the viewpoint of superstring theory and mirror symmetry in mathematical physics. 

\medskip

For the anti-canonical rings of some smooth projective varieties, their NCCRs can be obtained by considering the total space of the canonical bundle. 

\begin{Thm}{(cf. \cite[Section~7]{VdB_NCCR})}
\label{motivating_NCCR}
Let $X$ be a smooth projective Fano variety of dimension $d$. 
Suppose that $X$ has a tilting bundle $T$ such that the global dimension of $\End_X(T)$ is $d$. 
Let $Y$ be the total space of the canonical bundle $\om_X$ of $X$ with a canonical map $\pi\colon Y\to X$. 
Then $\End_Y(\pi^*T)$ is a $(d+1)$-Calabi--Yau algebra (see Subsection~\ref{subsec_CYalg}) and is an NCCR of the section ring $\Gamma(Y, \O_Y)=\bigoplus_{i\geq0}\Ga(X,\om_X^{-i})$. 
\end{Thm}

The crucial points of this construction are that the total space $Y$ of the canonical bundle $\om_X$ is a Calabi--Yau (CY) variety in the sense that $\om_Y\cong\O_Y$, and that $\pi^*T$ is a tilting bundle on $Y$ (cf. \cite{BF,BS_helix,BuH,Tom_2tilt,VdB_NCCR}). 
As a more general construction of a CY variety, we consider the total space of a vector bundle $E$ on a smooth projective variety (or stack) $X$. We do not require $E$ to be the canonical bundle $\om_X$, but its determinant to be $\om_X$, which ensures that the total space $Y$ of $E$ is CY. When $Y$ moreover has a tilting object $U$, the CY property of $Y$ implies that the endomorphism algebra of $U$ is a CY algebra. We refer to \cite[Section~5]{Lam_thesis} for this idea.

The construction of taking the total space $Y$ of the canonical bundle over a variety $X$ is derived equivalent to an algebraic operation of the \emph{CY completion} (or the \emph{derived preprojective algebra}) \cite{Ke_CYcomp,HIO}. 
More precisely, in the setting of Theorem~\ref{motivating_NCCR}, the variety $X$ is derived equivalent to $A:=\End_X(T)$, and the total space $Y$ of $\om_X$ is derived equivalent to the $(d+1)$-CY completion $\Pi=\Pi_{d+1}(A)$ of $A$ which is isomorphic to $\End_Y(\pi^\ast T)$.
Furthermore, the structure as the total space of a line bundle is translated to a grading on the CY algebra $\Pi$ whose Gorenstein parameter is $1$ and its degree $0$ part is a \emph{$d$-representation infinite algebra} (see \cite{HIO},\cite[Theorem~4.12]{MM}). 

The first main result of this paper is the following generalization of Theorem~\ref{motivating_NCCR}, which shows that we obtain CY algebras of higher Gorenstein parameters from higher rank vector bundles. 
For such CY algebras obtained in this way, by extending the argument for isolated singularities in \cite{AIR}, we can construct an NCCR of the section ring. 
\begin{Thm}[Theorem~\ref{CY} and Corollary~\ref{CY-tilting}]\label{intro-main1}
Let $X$ be a smooth projective variety of dimension $d$, let $E$ be a vector bundle of rank $r$ on $X$ such that $\Wedge^r E=\om_X$, and denote $\pi\colon Y=\Spec(\Sym E^\vee)\to X$. 
Suppose that $T\in\Db(\coh X)$ is a tilting object such that $\pi^\ast T\in\Db(\coh Y)$ is still tilting, and put 
\[ A=\End_X(T), \qquad \Ga=\End_Y(\pi^\ast T)=\bigoplus_{i\geq0}\Hom_X(T,T\otimes_X\Sym^iE^\vee), \]
and naturally regard $\Ga$ as a graded algebra.
\begin{enumerate}
\item There exists an isomorphism
\[ \RHom_{\Ga^e}(\Ga,\Ga^e)[d+r]\cong\Ga(r) \text{ in } \rD^\Z(\Ga^e), \]
that is, $\Ga$ is a $(d+r)$-CY algebra of Gorenstein parameter $r$. 
\item If $T$ contains a line bundle and $E^\vee$ is semiample, then $\Ga$ is a module-finite algebra and is an NCCR of the section ring $R=\Ga(Y,\O_Y)$ which is a Gorenstein normal domain. Moreover $Y$ is a crepant resolution of $\Spec R$.
\end{enumerate}
\end{Thm}
We note that the same statements for (1) hold without tilting assumptions in the setting of dg algebras. We also give an explicit bimodule resolution of the CY algebra $\Ga$ which generalizes that for the tensor algebra when $E=\om_X$, see Theorem~\ref{CY} for details.

In the setting where $E=\om_X$, we know that the endomorphism algebra of a tilting object $T$ over $X$ is $d$-representation infinite if and only if $\pi^\ast T$ is tilting over $Y$ \cite{BuH}. We give a partial generalization of this result. In the setting of Theorem~\ref{intro-main1}, we show that if $E$ consists of line bundles then $\End_X(T)$ is $d$-representation infinite (see Proposition~\ref{Prop_BuH}), while the converse does not hold in general (see Example~\ref{ex_fail_tilting}).

We also give a stack version of Theorem~\ref{intro-main1}, see Theorem~\ref{CY2} and Corollary~\ref{CY2-tilting}, which will be applied to some toric Deligne--Mumford stacks including Geigle--Lenzing projective spaces in Sections \ref{sec_toric_stack} and \ref{sec_GLspace}.

\medskip

Recall from above that when $E=\om_X$, the operation of taking the total space amounts to the CY completion. In particular in this case, the CY algebra $\Ga$ obtained in Theorem~\ref{intro-main1} is the $(d+1)$-CY completion of $A=\End_X(T)$. 
Our next result is that, even for a general vector bundle $E$, we may regard the resulting CY algebra $\Ga$ in Theorem~\ref{intro-main1} as a generalized CY completion in the sense of \cite{Ha7}. Recall that it is an operation for a {\it root pair} $(U,P)$ on a dg algebra which consists of a bimodule $U$ and a projective module $P$ satisfying some conditions, see Definition~\ref{root}, which produces a CY algebra as the CY completion (Theorem~\ref{root-CY}). 
A prototypical example of a root pair is obtained from Beilinson's tilting bundle on $\PP^n$ (see Example \ref{root-ex}).
We show that the CY algebra $\Ga$ in Theorem~\ref{intro-main1} is the CY completion of a root pair obtained from the projective bundle over $X$ associated with $E$.

\begin{Thm}
Let $X$ be a smooth projective variety of dimension $d$, $E$ an anti-ample vector bundle of rank $r$ on $X$ such that $\Wedge^r E=\om_X$. Let $\pi\colon Y=\Spec(\Sym E^\vee)\to X$ be the associated total space and $p\colon Z=\Proj(\Sym E^\vee)\to X$ its projectivization. 
Suppose that $T\in\Db(\coh X)$ is a tilting bundle and that $\pi^\ast T\in\Db(\coh Y)$ is also tilting. Put $\Ga=\End_Y(\pi^\ast T)=\bigoplus_{i\geq0}\Hom_X(T,T\otimes_X\Sym^iE^\vee)$. 
\begin{enumerate}
\item The object $U=\bigoplus_{i=0}^{r-1}p^\ast T(i)\in\Db(\coh Z)$ is tilting, $B=\End_Z(U)$ is $(d+r-1)$-representation infinite, and has a strict $r$-th root pair (see Definitions \ref{root}, \ref{strict}).
\item $\Ga$ is the Calabi--Yau completion of the root pair on $B$ in (1).
\end{enumerate}
\end{Thm}
In fact, we show (Theorem \ref{proj}) that the non-commutative projective scheme $\qmod^\Z\!\Ga$ associated with $\Ga$ in the sense of Artin--Zhang \cite{AZ} is equivalent to the category $\coh Z$ of coherent sheaves over $Z$, which is an essential observation for the proof.

\subsection{Applications}

We also explicitly provide a vector bundle $E$ and a tilting bundle $T$ satisfying the conditions in Theorem~\ref{intro-main1} for several varieties  (stacks), 
thereby constructing an NCCR. Note that vector bundles $T$ and $E$ must be chosen carefully in this process. 
For instance, as shown in Example~\ref{ex_fail_tilting}, even if the pull-back of $T$ to the total space of the canonical bundle is tilting, it may fail to be tilting on the total space of a different vector bundle. 

As an application of Theorem~\ref{intro-main1}, we are particularly interested in NCCRs of toric singularities. 
It is conjectured \cite[Conjecture~3.7]{VdB_ICM} that every toric singularity has an NCCR although it is open, in general. 
For example, we know the existence of NCCRs for the following classes of toric rings (singularities): 
\begin{itemize}
\setlength{\parskip}{0pt} 
\setlength{\itemsep}{3pt}
\item Toric rings whose class groups are isomorphic to $\ZZ$, see \cite[Theorem~8.9]{VdB_NCCR}. 
\item Three-dimensional Gorenstein toric rings, see e.g., \cite{Bro_dimer, IU_special,SpVdB_toricNCCR2}. 
\item Toric rings defined by quasi-symmetric representations, see \cite[Theorem~1.19]{SpVdB_NCCR}. 
\item Segre products of polynomial rings with suitable gradings \cite{Han_NCCR,HigNak_conic}. 
\end{itemize}
Besides these results, there are several results concerning the existence of NCCRs of toric rings, 
e.g., \cite{HigMat,MS_NCCRalmost,MS_NCCR,Mat_conic,Nak_HibiNCCR,SpVdB_toricNCCR1,Tom_rankone}. 
In many cases, a module that gives an NCCR of a toric ring is a direct sum of rank one reflexive modules, 
in which case an NCCR is called a \emph{toric} (or \emph{splitting}) NCCR (cf. \cite{Boc_toricNCCR,IN_steady}). 
However, there exists a toric ring that has an NCCR, but does not have a toric one, see \cite[Example~5.4]{SpVdB_toricNCCR1}. 

\medskip

In this paper, we apply Theorem~\ref{intro-main1} to various smooth projective varieties and toric Deligne--Mumford stacks. It yields NCCRs for some toric singularities which are in fact toric NCCRs.

Our first result arises from a product of projective spaces. 
We define the ring $R$ as follows. Let $r,s\geq1$ and $n_1,\ldots,n_s\geq1$ be integers, and consider the polynomial ring
\[ S=k[x_{1,0},\ldots,x_{1,n_1},x_{2,0},\ldots,x_{2,n_2},\ldots,x_{s,0},\ldots,x_{s,n_s},y_1,\ldots,y_r]. \] 
We give a $\Z^s$-grading on $S$ by
\begin{align*}
\deg x_{i,j}&\coloneqq\e_i=(0, \dots, 0, \overset{i}{\check{1}}, 0,\dots, 0)\in\Z^s, \quad 0\leq j\leq n_i\\
\deg y_l&\coloneqq c_l=(c_{1,l}, \dots, c_{s,l})\in\Z^s,
\end{align*}
and set $R=S_{(0,\ldots,0)}$, the degree $0\in\Z^s$ part of $S$.

For each $m\geq0$ we denote by $T^m:=\bigoplus_{i=0}^m\O(i)\in\coh\PP^m$ Beilinson's tilting bundle.
\begin{Thm}[Theorem~\ref{NCCR_toric_PxP} and Corollary~\ref{CCR_toric_PxP}]
\label{NCCR_toric_PxP_intro}
Let $n_1, \dots, n_s$ be positive integers. 
Let $X= \PP^{n_1}\times\cdots\times\PP^{n_s}$ and define the tilting bundle $T$ and a vector bundle $E$ over $X$ by
\[ T=T^{n_1}\boxtimes T^{n_2}\boxtimes\cdots\boxtimes T^{n_s}, \qquad E=\bigoplus_{l=1}^r\O_X(c_l), \]
where $c_l\coloneqq (c_{1,l}, \dots, c_{s,l})\in\ZZ^s$ with $c_{1,l}, \dots, c_{s,l}\le 0 \ (l=1, \dots, r)$ and $\sum_{l=1}^rc_l=(-n_1-1, \dots, -n_s-1)$. 
We consider the total space $Y$ of $E$ and the map $\pi\colon Y\rightarrow X$, thus the section ring $\Ga(Y,\O_Y)$ is isomorphic to $R$ defined above.
\begin{enumerate}
\item $\pi^\ast T$ is a tilting bundle on $Y$ and $\End_Y(\pi^\ast T)$ is an NCCR of $R$. 
\item The canonical map $Y\to\Spec R$ is a crepant resolution.
\end{enumerate}
\end{Thm}

Note that Theorem~\ref{intro-main1} gives NCCRs of some families of toric rings. 
For example, we can obtain NCCRs of some Hibi rings as shown in Subsection~\ref{subsec_Hibi}. 
NCCRs of some classes of Hibi rings have already been constructed in \cite{HigNak_conic,Nak_HibiNCCR}. 
Our theorem not only recovers the results of these papers, but also gives NCCRs for a broader class of Hibi rings. 
Also, the construction via Theorem~\ref{intro-main1} uses an approach different from those employed in \cite{HigNak_conic,Nak_HibiNCCR}. 

\medskip

By applying the stack version of Theorem~\ref{intro-main1} (Theorem~\ref{CY2} and Corollary~\ref{CY2-tilting}) to toric Deligne--Mumford (DM) stacks admitting tilting bundles, we can also construct NCCRs of some toric rings. 
For example, using a suitable tilting bundle $T$ and a vector bundle $E$ on a weighted projective stack, 
we obtain an NCCR of a toric ring whose class group is $\ZZ$, see Theorem~\ref{NCCR_toricZ_revisited}. 
Such an NCCR was already constructed in \cite{VdB_NCCR} (see also \cite{Tom_rankone}). Our construction uses a different approach, which also gives an interpretation of such NCCRs as being derived equivalent to the total space of a vector bundle on a weighted projective stack.
In addition, we construct several NCCRs by applying our method to toric DM stacks of Picard rank two that admit tilting bundles. 
For example, we apply our theorem to toric DM stacks whose coarse moduli spaces are products of weighted projective spaces (Subsection~\ref{Subsec_ex_toricDM1}) 
or the Hirzebruch surface of degree one (Subsection~\ref{Subsec_ex_toricDM2}), and obtain NCCRs of some toric rings.

\medskip

We can also construct some NCCRs by applying the stack version of Theorem~\ref{intro-main1} to a certain class of Geigle--Lenzing projective spaces introduced in \cite{HIMO}. We focus on the case where the Geigle--Lenzing complete intersection $(R,\LL)$ is isomorphic to a polynomial ring.

Let $R=k[X_1,\ldots,X_n]$ be the polynomial ring graded by the abelian group 
\[\LL=\bigoplus_{i=1}^n\Z\vec{x}_i/(p_1\vec{x}_1=\cdots=p_n\vec{x}_n\eqqcolon \vec{c})\]
by setting $\deg X_i=\vec{x}_i$, in which the scheme $\Spec R$ has an action of a reductive algebraic group $G\coloneqq \Spec k[\LL]$.
Then, we consider the Geigle--Lenzing projective space, which is the quotient stack 
\[ \XX=[(\Spec R\setminus \{R_+\})/G] \]
for the $G$-invariant maximal ideal $R_+=(X_1,\ldots,X_n)$. 
In this case, the Geigle--Lenzing projective space $\X$ has a tilting bundle $T=\bigoplus_{0\leq\vec{x}\leq (n-1)\vec{c}}\O_\XX(\vec{x})$, whose endomorphism algebra has global dimension at most $n-1=\dim\X$ \cite[Theorems~1.13 and 1.14]{HIMO}. 
Moreover, we observe that such a Geigle--Lenzing projective space can be described as a toric Deligne--Mumford stack. Applying the stack version of Theorem~\ref{intro-main1} leads to the following result.

\begin{Thm}[Theorem~\ref{NCCRfromGL}]\label{NCCR_GL_intro}
Let $\X$ be a Geigle--Lenzing projective space as above.
Let $E=\bigoplus_{j=1}^r\O_\XX(-\vec{y}_j)$ be a vector bundle of rank $r$ on $\XX$ such that $\vec{y}_1, \dots, \vec{y}_r\in\LL_+$ and 
$\vec{x}_1+\cdots+\vec{x}_n=\vec{y}_1+\cdots+\vec{y}_r$. 
We consider the total space $\YY$ of $E$ and the map $\pi\colon \YY\to \XX$. 
Let $T$ be the tilting bundle on $\X$ as above. Then $\pi^\ast T$ is a tilting bundle on $\YY$, and 
$\End_\YY(\pi^\ast T)$ is an NCCR of the section ring $\Ga(\YY,\O_\YY)$.
\end{Thm}

When each $y_j$ is non-zero and $r\ge 2$, then the rank of the class group of the section ring $\Ga(\YY,\O_\YY)$ is one, so an NCCR of $\Ga(\YY,\O_\YY)$ can also be constructed using the methods in \cite{Tom_rankone,MS_NCCRalmost}. 

In the case where $\X$ is a classical Geigle--Lenzing weighted projective line, taking the vector bundle $E$ 
in Theorem~\ref{NCCR_GL_intro} to be the canonical bundle $\om_\X$, we obtain an NCCR of a Kleinian singularity of type A, 
which also follows from \cite[Corollary~1.20]{HIMO}. 
By choosing a different bundle $E$, we can also obtain an NCCR of a toric compound Du Val singularity of type A (see Example~\ref{Ex_WPL}). 
Although various constructions of NCCRs of toric compound Du Val singularities are known (see e.g., \cite{IW_cAn,Nag,Nak_cAn,Tom_rankone,Wem_MMP}), 
our Theorem~\ref{NCCR_GL_intro} gives a geometric interpretation of such constructions.


\section{Preliminaries on Calabi--Yau algebras and tilting theory}
\subsection{Tilting objects in triangulated categories}
We recall some basic facts on derived equivalences of varieties and algebras, given by the following fundamental notion.
\begin{Def}
Let $\T$ be a triangulated category. 
\begin{enumerate}
\item A {\it generator} of $\T$ is an object $T\in\T$ such that the thick subcategory $\thick T$ generated by $T$ is the whole $\T$.
\item A {\it tilting object} of $\T$ is a generator $T$ of $\T$ such that $\Hom_\T(T,T[i])=0$ for all $i\neq0$.
\end{enumerate}
If $X$ is a variety and $T$ is a vector bundle on $X$ which gives a tilting object in $\per X$, then we call such $T$ a {\it tilting bundle}.
\end{Def}

The following result explains the significance of tilting objects.
\begin{Thm}[\cite{Ke}]
Let $\T$ be an algebraic triangulated category. 
If $T\in\T$ is a tilting object, then there is a triangle equivalence $\T\simeq\per\End_\T(T)$ up to direct summands.
\end{Thm}

In this paper, we focus particularly on the following setting: 

\begin{Setup}\label{general_settingXY}
Let $X$ be a smooth projective variety of dimension $d$, let $E$ be a vector bundle over $X$ of rank $r$, and let $\pi\colon Y=\Spec(\Sym E^\vee)\to X$ be the natural projection. 
\begin{enumerate}
\item\label{general} $\Wedge^rE=\om_X$. 
\item\label{semiample} In addition to (1), $E^\vee$ is semiample, that is, $\Sym^mE^\vee$ is globally generated for some $m>0$.
\end{enumerate}

\end{Setup}
Under the setting (\ref{general}), the total space $Y$ is a Calabi--Yau variety since $\om_Y=\pi^*(\om_X\otimes \Wedge^rE^\vee)=\O_Y$. 
Also, under (\ref{semiample}), the section ring $\Ga(Y,\O_Y)$ is finitely generated by \cite[Example 2.1.30]{Lazarsfeld}.
We shall use semiampleness to ensure certain finiteness, for example, the finite generation of section rings. 

In this setting, we will consider whether the pull-back $\pi^*T$ of a tilting bundle $T$ on $X$ is a tilting bundle on $Y$. 
We first note that each generator pulls back to a generator, 
but, as we will see in Example~\ref{ex_fail_tilting}, this is not necessarily the case for a tilting bundle. 

\begin{Lem}
\label{pullback_generator}
Let $X$ be a smooth projective variety, $E$ a vector bundle on $X$, and $\pi\colon Y=\Spec (\Sym E^\vee)\to X$. 
Suppose that $X$ has a generator $T$. Then $\pi^\ast T$ is a generator of $\Db(\coh Y)$. 
\end{Lem}

\begin{proof}
The argument is similar to the one in \cite[Proposition~4.1]{Bri_tstructure}. 
By \cite{Ne92}, it is enough to show that if $V\in\Db(\coh Y)$ satisfies $\RHom_Y(\pi^\ast T,V)=0$ then $V=0$. By the adjunction $\mathrm{L}\pi^\ast=\pi^\ast\dashv \pi_\ast=\mathrm{R}\pi_\ast$, we have 
\[
0=\RHom_Y(\pi^\ast T,V)=\RHom_X(T, \pi_\ast V). 
\]
Thus,  we have $\pi_\ast V\cong 0$ since $T$ is a generator of $\per X$, and hence a compact generator of $\rD(\Qcoh X)$. 
Since $\pi$ is an affine morphism, this implies $V=0$. 
\end{proof}

\subsection{Calabi--Yau algebras and $d$-representation infinite algebras}\label{subsec_CYalg}
Recall that a differential graded (dg) algebra is a (cohomologically) graded algebra $A=\bigoplus_{i\in\Z}A^i$ endowed with a differential of degree $1$ satisfying the Leibniz rule $d(ab)=da\cdot b+(-1)^{\deg a}a\cdot db$. For example, an ordinary algebra is a dg algebra concentrated in degree $0$.
As in the case of ordinary algebras, one has the {\it derived category} $\rD(A)$ of a dg algebra $A$ which is defined as the localization of the homotopy category of dg $A$-modules at quasi-isomorphisms.

For a dg algebra $A$, we denote by $A^e:=A^\op\otimes_kA$ the enveloping algebra. A dg $k$-algebra $\Ga$ is called {\it smooth} if it is perfect as a bimodule: $\Ga\in\per\Ga^e$. Following \cite{Gi,Ke_CYcomp}, a dg algebra $\Ga$ is {\it $n$-Calabi--Yau ($n$-CY)} if it is smooth and satisfies
\[ \RHom_{\Ga^e}(\Ga,\Ga^e)[n]\cong\Ga \]
in $\rD(\Ga^e)$.
A systematic construction of such (dg) algebras is given by the concept of {\it Calabi--Yau (CY) completions} \cite{Ke_CYcomp}. Given a smooth dg algebra $A$, its {\it $n$-Calabi--Yau ($n$-CY) completion} is the derived tensor algebra
\[ \Pi=\Pi_{n}(A):=\rT^\rL_A(\RHom_{A^e}(A,A^e)[n-1]), \]
that is, the tensor algebra $\rT_A\!\theta$ of a cofibrant resolution $\theta\to\RHom_{A^e}(A,A^e)[n-1]$. Here, the bimodule $\RHom_{A^e}(A,A^e)$ is called the {\it inverse dualizing bimodule}. By \cite{Ke_CYcomp}, the dg algebra $\Pi_n(A)$ is indeed $n$-CY. 

Notice that, in addition to the ordinary cohomological grading on dg algebras, the CY completion $\Pi$ has an additional grading, called the {\it Adams grading}, given by the tensor grading. The proof of the CY property of $\Pi$ shows that, with respect to this Adams grading, $\Pi$ has Gorenstein parameter $1$ in the sense that there is an isomorphism in $\rD^\Z(\Pi^e)$: 
\[ \RHom_{\Pi^e}(\Pi,\Pi^e)[n]\cong\Pi(1). \]

As explained in \cite{Ke_CYcomp,IQ} the CY completion should be regarded as the non-commutative analogue of the total space of the canonical bundle, that is, if a smooth variety $X$ of dimension $n-1$ is derived equivalent to a dg algebra $A$, then the total space $Y$ of $\om_X$ is derived equivalent to the $n$-CY completion $\Pi_n(A)$.
We shall observe in Section \ref{projbundle} that an analogue of the CY completion can be obtained by considering the total space of a vector bundle which is not necessarily the canonical bundle. 

As we will see in Proposition~\ref{Prop_BuH} below, there is a relationship between $\pi^\ast T$ being a tilting bundle on $Y$ and the algebra $\End_X(T)$ being a $d$-representation infinite algebra in the following sense. 

\begin{Def}[\cite{HIO}]\label{def_RI}
Let $A$ be a finite dimensional algebra over $k$ and $n\geq0$. We say $A$ is {\it $n$-representation infinite} if $\gldim A\leq n$ and 
\[ \nu_n^{-i}A\in\mod A \]
for every $i\geq0$, where $\nu_n^{-1}=-\lotimes_A\RHom_A(DA,A)[n]$.
\end{Def}

We have the following, in which the case $E=\om_X$ is treated in \cite{BF,BuH}.
\begin{Prop}[{cf.\! \cite{BF,BuH}}]\label{Prop_BuH}
Let $X$, $E$, and $Y$ be the same as in Setting~\ref{general_settingXY}(\ref{general}).
Suppose that $X$ has a tilting bundle $T$ and consider the following conditions.
\begin{enumerate}
\renewcommand\labelenumi{(\roman{enumi})}
\renewcommand\theenumi{\roman{enumi}} 
\item\label{RI} $A=\End_X(T)$ is $d$-representation infinite.
\item\label{omExt} $H^i\RHom_X(T,T\lotimes_X\om_X^{-l})=0$ for all $l\geq0$ and $i\neq0$.
\item\label{Y_tilting} $\pi^\ast T$ is a tilting bundle on $Y$.
\end{enumerate}
Then we have the following.
\begin{enumerate}
\item The conditions (\ref{RI}) and (\ref{omExt}) are equivalent.
\item If the vector bundle $E$ consists of line bundles, then the implication (\ref{Y_tilting})$\Rightarrow$(\ref{RI}) holds. 
\item If $E=\om_X$, then all conditions are equivalent.
\end{enumerate}
\end{Prop}
\begin{proof}
	(1)  Since $T\in\Db(\coh X)$ is a tilting bundle we have an equivalence $\Db(\coh X)\simeq\Db(\mod A)$. Comparing the Serre functors, we obtain a commutative diagram
	\[ \xymatrix@R=5mm{
		\Db(\coh X)\ar[r]^-\simeq\ar[d]_-{-\lotimes_X\om_X^{-1}}&\Db(\mod A)\ar[d]^-{\nu_d^{-1}}\\
		\Db(\coh X)\ar[r]^-\simeq&\Db(\mod A). } \]
	By definition, $A$ is $d$-representation infinite if and only if $H^i\RHom_A(A,\nu_d^{-l}A)=0$ for every $l\geq0$ and $i\neq0$. By the above diagram, this is equivalent to $H^i\RHom_X(T,T\lotimes_X\om_X^{-l})=0$. This proves the assertion.

	(2)  By (1) it is enough to prove the conclusion of (\ref{omExt}). Write $E$ as the sum $\bigoplus_{i=1}^rL_i$ of line bundles. 
	Since $\Wedge^rE=\om_X$, we have $L_1\cdots L_r=\om_X$. We have 
	\[\RHom_Y(\pi^\ast T,\pi^\ast T)=\RHom_X(T,\pi_\ast\pi^\ast T)=\bigoplus_{i\geq0}\RHom_X(T,T\otimes_X\Sym^iE^\vee)\] by the projection formula. 
	If $\pi^\ast T$ is a tilting bundle, then the above complex concentrates in degree $0$. Note that since $E=\bigoplus_{i=1}^rL_i$, the $r$-th symmetric product $\Sym^rE^\vee$ contains $L_1^{-1}\cdots L_r^{-1}=\om_X^{-1}$ as a direct summand. In particular, every power $\om_X^{-l}$ appears as a direct summand of some symmetric power $\Sym^iE^\vee$. Therefore the vanishing of extensions implies $H^i\RHom_X(T,T\lotimes_X\om_X^{-l})=0$ for all $i\neq0$ and $l\geq0$.

	(3)  It is enough to show (\ref{omExt})$\Rightarrow$(\ref{Y_tilting}). As in the above computation, we have 
	\[\RHom_Y(\pi^\ast T,\pi^\ast T)=\RHom_X(T,\pi_\ast\pi^\ast T)=\bigoplus_{i\geq0}\RHom_X(T,T\otimes_X\Sym^iE^\vee)=\bigoplus_{i\geq0}\RHom_X(T,T\otimes_X\om_X^{-i}).\] 
	Therefore, vanishing of the non-zero cohomologies in the last term implies $\pi^\ast T$ is a tilting bundle.
\end{proof}

We note the following example where (\ref{RI})$\Rightarrow$(\ref{Y_tilting}) is not true in general.

\begin{Ex}\label{ex_fail_tilting}
Let $X=\PP^1\times\PP^1$, so that $\om_X=\O_X(-2,-2)$. 
We consider two vector bundles $E_1=\om_X$ and $E_2=\O_X(-1,-1)^{\oplus2}$ on $X$ with the total spaces $\pi_i\colon Y_i=\Spec(\Sym E_i^\vee)\to X$ $(i=1,2)$. 
Note that $\Wedge^2 E_2=E_1$. 
Take a tilting bundle $T=\O_X\oplus\O_X(1,0)\oplus\O_X(2,1)\oplus\O_X(3,1)$ on $X$.
It is easy to see that $\End_X(T)$ has global dimension $2=\dim X$, and thus it is $2$-representation infinite and $\pi_1^\ast T$ is a tilting bundle on $Y_1$ by \cite[Theorem~3]{BuH}. 
The algebra $\End_{Y_1}(\pi_1^\ast T)$ is $3$-CY which can be obtained as the CY completion (cf. \cite{Ke_CYcomp}) of 
the $2$-representation infinite algebra $\End_X(T)$, and it is also obtained from a consistent dimer model (see e.g., \cite[5.22]{Nak_reflexive}). 
Figure~\ref{counterEX_repinf} shows the quivers of $\End_X(T)$ and $\End_{Y_1}(\pi_1^\ast T)$. 

On the other hand, let us compute the extensions of the pull-back $\pi_2^\ast T$. 
We have $\RHom_{Y_2}(\pi_2^\ast T,\pi_2^\ast T)=\bigoplus_{i\geq0}\RHom_X(T,T\otimes_X\Sym^iE_2^\vee)$. 
It contains $\RHom_X(\O_X(3,1),\O_X\otimes_X\O_X(1,1))=\RHom_X(\O_X,\O_X(-2,0))$ as a direct summand, which has non-zero cohomology in degree $1$. 
Thus, $\pi_2^\ast T$ is not a tilting bundle on $Y_2$. 

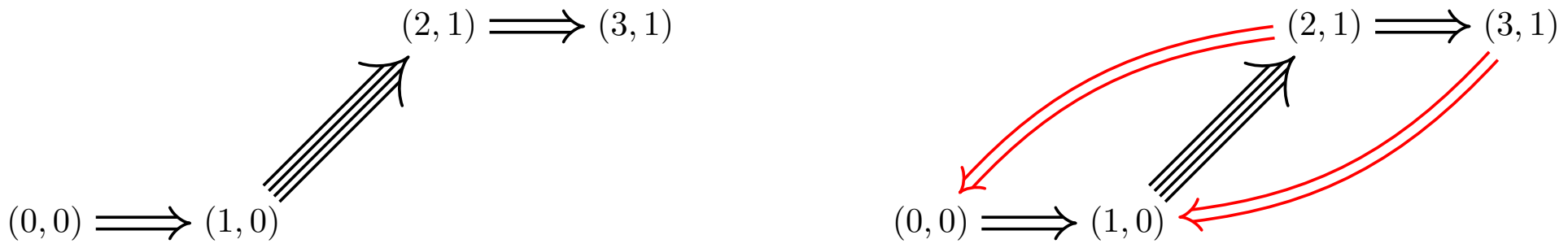
\begin{figure}[H]
\begin{center}
{\scalebox{0.8}{
\begin{tikzpicture}
\newcommand{\arrowwidth}{0.03cm} 

\node at (0,0) {
\begin{tikzpicture}
\foreach \a/\b in {0/0,1/0,2/1,3/1} {\node (V\a\b) at (\a*2,\b*2) {$(\a,\b)$} ;}; 
\foreach \s/\t in {00/10,21/31} {\draw[double distance=3pt, scaling nfold=2, -Implies, line width=\arrowwidth] (V\s)--(V\t);};
\foreach \s/\t in {10/21} {\draw[double distance=1.8pt, scaling nfold=4, -Implies, line width=\arrowwidth]  (V\s)--(V\t);};
\end{tikzpicture}
};

\node at (9,0) {
\begin{tikzpicture}
\foreach \a/\b in {0/0,1/0,2/1,3/1} {\node (V\a\b) at (\a*2,\b*2) {$(\a,\b)$} ;}; 
\foreach \s/\t in {00/10,21/31} {\draw[double distance=3pt, scaling nfold=2, -Implies, line width=\arrowwidth] (V\s)--(V\t);};

\foreach \s/\t in {10/21} {\draw[double distance=1.8pt, scaling nfold=4, -Implies, line width=\arrowwidth]  (V\s)--(V\t);};

\draw[double distance=3pt, scaling nfold=2, -Implies, line width=\arrowwidth, red] (V21)  [bend right=20] to (V00); 
\draw[double distance=3pt, scaling nfold=2, -Implies, line width=\arrowwidth, red] (V31) [bend left=20] to (V10); 
\end{tikzpicture}
};
\end{tikzpicture}
} }
\caption{The quivers of the $2$-representation infinite algebra $\End_X(T)$ (left) and the $3$-CY algebra $\End_{Y_1}(\pi_1^\ast T)$ (right)}
\label{counterEX_repinf}
\end{center}
\end{figure}
\end{Ex}

\section{NCCRs via Calabi--Yau varieties and tilting}
We give a construction of CY algebras which in turn give NCCRs of certain singularities, see Theorem \ref{CY} and Corollary \ref{CY-tilting}. We do so by tilting a CY variety obtained as the total space of a vector bundle $E$ whose determinant is the canonical bundle $\om_X$. This generalizes a well-known construction for $E=\om_X$.

\subsection{NCCRs from Calabi--Yau algebras}
In this subsection, we show (see Proposition \ref{eBe}) that a module-finite CY algebra yields an NCCR of its idempotent subalgebra by extending the approach of \cite{AIR} to the non-isolated singularities.

Let us prepare some terminology on (symmetric) orders. 
Let $R$ be a commutative Gorenstein ring. A module-finite $R$-algebra $\La$ is called an {\it $R$-order} if $\La$ is (maximal) Cohen--Macaulay as an $R$-module. An $R$-order is {\it symmetric} if $\Hom_R(\La,R)\cong\La$ in $\Mod\La^\op\otimes_R\La$. We say that an $R$-order is {\it locally symmetric} if $\La_\m$ is a symmetric $R_\m$-order for every $\m\in\Max R$.

Also, if $\Ga$ is a module-finite algebra over a commutative Noetherian ring $R$, we say that a $\Ga$-module $M$ is a {\it height $0$ progenerator} if $M_\p\in\mod\Ga_\p$ is a progenerator for all $\p\in\Spec R$ with $\height\p=0$.
If such a $\Ga$ is a prime ring, then $e\Ga$ is a height $0$ progenerator for every non-zero idempotent $e\in\Ga$.

\begin{Lem}\label{codim}
Let $\Ga$ be a module-finite algebra over a finitely generated $k$-algebra. Suppose that $\Ga$ is $d$-CY as a $k$-algebra.
\begin{enumerate}
\item There exists a central subalgebra $R\subset\Ga$ which is isomorphic to $k[x_1,\ldots,x_d]$ such that $\Ga$ is a locally symmetric $R$-order.
\item For each idempotent $e\in \Ga$ such that $e\Ga$ is a height $0$ progenerator, we have $\dim\Ga/(e)\leq d-2$.
\end{enumerate}
\end{Lem}
\begin{proof}
	(1)  By \cite[Theorem 3.1(6)]{IR_FZ}, the dimension of $\Ga$ is equal to the CY dimension $d$. Let $R$ be a Noether normalization of the center of $\Ga$, so we have $R\cong k[x_1,\ldots,x_d]$. Then by \cite[Theorem 3.3(2)]{IR_FZ}, $\Ga_\m$ is a symmetric $R_\m$-order of global dimension $d$ for each $\m\in\Max R$. 
	
	(2)  Let $R$ be a central subalgebra as in (1).
    It follows from (1) that $\Ga_\p$ is a symmetric $R_\p$-order of global dimension $1$ for every $\p\in\Spec R$ of height $1$. 
	Then the completion $\widehat{\Ga_\p}$ is a hereditary order over a complete discrete valuation ring $\widehat{R_\p}$ which is symmetric, so it has to be a product of matrix algebras over maximal $\widehat{R_\p}$-orders \cite{Reiner}, see also \cite[Proposition 3.12]{IR_FZ}. We claim that the idempotent $e\in\widehat{\Ga_\p}$ has non-zero components in every matrix factors. This will imply $\widehat{\Ga_\p}/(e)=0$, and therefore the assertion.

    For to prove the claim it is enough to show that the canonical map $\widehat{\Ga_\p}\to\End_{e\widehat{\Ga_\p}e}(\widehat{\Ga_\p}e)$ is injective. 
    But this follows from the assumption that $e\Ga$ is a height $0$ progenerator as this yield that the canonical map $\Ga\to\End_{e\Ga e}(\Ga e)$ is an isomorphism after tensoring the quotient field of $R$. 
\end{proof}

We also use the following non-isolated singularity version of \cite[Theorem 2.3]{AIR} which produces an NCCR from a CY algebra.
\begin{Prop}\label{eBe}
Let $\Ga$ be a module-finite algebra over a finitely generated $k$-algebra, and let $e\in\Ga$ be an idempotent. Suppose $\Ga$ is $d$-CY as a $k$-algebra. 
\begin{enumerate}
\item If $e\Ga$ is a height $0$ progenerator, the canonical map $\Ga\to\End_{e\Ga e}(\Ga e)$ is an isomorphism.
\item We have $e\Ga\xsimeq\RHom_{e\Ga e}(\Ga e,e\Ga e)$ and $\Ga e\xsimeq\RHom_{e\Ga^\op e}(e\Ga,e\Ga e)$.
\item The algebra $e\Ga e$ is $d$-Iwanaga-Gorenstein.
\end{enumerate}
\end{Prop}
\begin{proof}
	By Lemma \ref{codim}(1), we may assume that $\Ga$ is a locally symmetric order over $R=k[x_1,\ldots,x_d]$. 
	Put $\La=e\Ga e$ and consider the triangle
	\[ \xymatrix{ \Ga e\lotimes_\La e\Ga \ar[r]&\Ga\ar[r]& X } \]
	in $\rD(\Ga^\op\otimes_R\Ga)$ obtained by extending the multiplication map. Then $X$ (which is in fact the dg quotient of $\Ga$ by $e$) is annihilated by $e$ from each side.
	Applying $\RHom_\Ga(-,\Ga)$ to the triangle yields
	\[ \xymatrix{ \RHom_\Ga(X,\Ga)\ar[r]& \Ga\ar[r]&\RHom_\La(\Ga e,\Ga e) }. \]
   
    To prove (1), we show that $H^0\RHom_\Ga(X,\Ga)=0$ and $H^1\RHom_\Ga(X,\Ga)=0$, which gives the assertion by the above triangle.
	Since each cohomology of $X$ is a module over $\Ga/(e)$, it has codimension $\geq2$ by Lemma \ref{codim}(2). It follows that $H^0\RHom_\Ga(X,\Ga)=\Hom_\Ga(H^0X,\Ga)\subset\Hom_R(H^0X,\Ga)=0$. Also, we have an exact sequence $\Ext^1_\Ga(H^0X,\Ga)\to H^1\RHom_\Ga(X,\Ga)\to\Hom_\Ga(H^{-1}X,\Ga)$ in which both end terms vanish by the codimension reason. Thus we obtain the assertion.

    We next show the first isomorphism of (2). By the above triangle, it is enough to show that $\RHom_\Ga(X,\Ga)$ is annihilated by $e$ from the {\it left}. We may localize at each prime ideal of $R$ and thus may assume that $\Ga$ is a symmetric $R$-order. Then $\RHom_\Ga(X,\Ga)=\RHom_R(X,R)$, which is then annihilated from the left by $Xe=0$.

    Once we know (2), the assertion (3) follows as in \cite{AIR}. We include the proof for the convenience of the reader. We prove $\RHom_\La(M,\La)$ is concentrated in degrees $\leq d$ for every $M\in\mod\La$. Let $N=M\otimes_\La e\Ga$ so that $Ne=M$. Then $\RHom_\La(M,\La)=\RHom_\La(N\lotimes_\Ga\Ga e,\La)$, which is isomorphic to $\RHom_\Ga(N,\RHom_\La(\Ga e,\La))$ by adjunction, and thus to $\RHom_\Ga(N,e\Ga)$ by (2). By $\gldim\Ga\leq d$, this is concentrated in degree $\leq d$, as desired.
  \end{proof} 

We next note the following information about the idempotent subring which will arise from our setup.
\begin{Lem}\label{secring}
Let $X$ be an integral normal variety, $E$ a vector bundle, and $Y$ its total space with $\pi\colon Y\to X$ the projection.
Suppose that $T\in\coh X$ contains a line bundle and let $e\in\Ga=\End_Y(\pi^\ast T)$ be the corresponding idempotent. Then the subalgebra $e\Ga e$ is isomorphic to the section ring $\bigoplus_{i\geq0}\Ga(X,\Sym^iE^\vee)=\bigoplus_{i\geq0}\Hom_X(\O_X,\O_X\otimes\Sym^iE^\vee)$, which is a normal domain.
\end{Lem}

\begin{proof}
	By shifting $T$, we may assume that $T$ contains $\O_X$ as a direct summand. Then by $\Ga=\End_Y(\pi^\ast T)=\bigoplus_{i\geq0}\Hom_X(T,T\otimes_X\Sym^iE^\vee)$, the idempotent subring $e\Ga e$ is precisely the section ring. 
    Since this is the ring of global sections of $Y$, which is an integral normal variety, it is a normal domain. 
\end{proof}


\subsection{Calabi--Yau algebras from vector bundles}
The CY property is well-known and the detailed proofs can be found, for example, in \cite{IQ,HK}. The novelty of the following result is an explicit bimodule resolution of the CY algebra (relative to the base). It also yields that the Gorenstein parameter of the CY algebra is the rank of the vector bundle, which is an important aspect in our discussion.

\begin{Thm}\label{CY}
Let $X$ be a smooth projective variety of dimension $d$, let $E$ be a vector bundle of rank $r$ on $X$ such that $\Wedge^r E=\om_X$, and denote $\pi\colon Y=\Spec(\Sym E^\vee)\to X$.
Let $T\in\Db(\coh X)$ be a generator and put 
\[ A=\REnd_X(T), \qquad \Ga=\REnd_Y(\pi^\ast T)=\bigoplus_{i\geq0}\RHom_X(T,T\lotimes_X\Sym^iE^\vee), \]
and naturally regard $\Ga$ as an Adams graded dg algebra.
\begin{enumerate}
\item The total of the double complex 
\begin{align*}
\Ga\lotimes_A\RHom_X(T,T\otimes_X\Wedge^rE^\vee)\lotimes_A\Ga(-r)&\longrightarrow\cdots \\
&\cdots\longrightarrow\Ga\lotimes_A\RHom_X(T,T\otimes_XE^\vee)\lotimes_A\Ga(-1)\longrightarrow\Ga\lotimes_A\Ga
\end{align*}
is isomorphic to $\Ga$ in $\rD^\Z(\Ga^e)$.
\item 
Denote by $(-)^\ast=\RHom_{A^e}(-,A^e)$ on $\rD(A^e)$. The above complex is isomorphic to its $[d]$-shifted bimodule dual over $\Ga$ given by the following. 
\begin{align*}
\Ga\lotimes_AA^\ast[d]\lotimes_A\Ga(-r)\longrightarrow \Ga\lotimes_A\RHom_X(T,T\otimes_X & E^\vee)^\ast[d] \lotimes_A\Ga(-r+1)\longrightarrow \cdots \\
&\cdots\longrightarrow \Ga\lotimes_A\RHom_X(T,T\otimes_X\Wedge^rE^\vee)^\ast[d]\lotimes_A\Ga 
\end{align*}
\item There exists an isomorphism
\[ \RHom_{\Ga^e}(\Ga,\Ga^e)[d+r]\cong\Ga(r) \text{ in } \rD^\Z(\Ga^e) \]
and $\Ga$ is a $(d+r)$-CY dg algebra of Gorenstein parameter $r$.
\end{enumerate}
\end{Thm}

Consequently, we obtain a CY algebra concentrated in degree $0$, and therefore an NCCR when $\pi^\ast T$ is tilting.
\begin{Cor}\label{CY-tilting}
Under the setting of Theorem \ref{CY}, suppose that $\pi^\ast T\in\Db(\coh Y)$ is tilting.
\begin{enumerate}
\item The algebra $\Ga=\End_Y(\pi^\ast T)=\bigoplus_{i\geq0}\Hom_X(T,T\otimes_X\Sym^iE^\vee)$ 
is $(d+r)$-CY of Gorenstein parameter $r$.
\item If $T$ contains a line bundle and $E^\vee$ is semiample, then $\Ga$ is a module-finite algebra and is an NCCR of $R:=\Ga(Y,\O_Y)$ which is a Gorenstein normal domain.
\item The canonical morphism $Y\to \Spec R$ is birational. Therefore, $Y$ is a crepant resolution of $\Spec R$.
\end{enumerate}
\end{Cor}
\begin{proof}
(1)  This is immediate from Theorem \ref{CY}.

(2)  By semiampleness, the section ring $R=\Ga(X,\Sym E^\vee)$ is a finitely generated $k$-algebra, over which $\Ga$ is module-finite. The ring $R$ is Gorenstein by Proposition \ref{eBe}(3) and is a normal domain by Lemma \ref{secring}. 

We claim that every non-zero idempotent $e\in\Ga=\End_Y(\pi^\ast T)$ gives a height $0$ progenerator $e\Ga$.
Let $K$ be the quotient field of $R$, and consider the pull-back below.
\[ \xymatrix{
Y_K\ar[r]^-j\ar[d]_-g&Y\ar[d]^-f\\
\Spec K\ar[r]^-i&\Spec R } \]
Since $U:=\pi^\ast T$ is a tilting object, so is $U_K:=j^\ast U$ over $Y_K$ since $j$ is flat. Also, since $f$ is crepant, so is $g$, and hence $Y_K$ is derived equivalent to a semisimple $K$-algebra. Since we have that $\Ga(Y_K,\O_{Y_K})=\Ga(Y,\O_Y)\otimes_RK=K$ is connected, we deduce that $Y_K$ is (derived equivalent to) a simple $K$-algebra. This shows that every non-zero idempotent $e\in\End_Y(U)$ is a height $0$ progenerator.

Now the assertion follows from (1) and Proposition \ref{eBe}(1).

(3) 
It follows from the above proof of (2) that $Y\to\Spec R$ is birational. Since $Y\to\Spec R$ is also projective, the remaining assertion follows from \cite[Corollary 4.15]{IW_Q}.
\end{proof}

We start the proof from the geometry side. Let $F$ be a locally free sheaf of rank $n$ on a variety $X$, and let $Y=\Spec(\Sym F^\vee)\to X$ be the associated vector bundle. Then we have the Koszul bimodule complex
\begin{equation}\label{K}
\xymatrix{ 0\ar[r]&\O_Y\otimes_X\Wedge^nF^\vee\otimes_X\O_Y\ar[r]&\cdots\ar[r]&\O_Y\otimes_X F^\vee\otimes_X\O_Y\ar[r]&\O_Y\otimes_X\O_Y\ar[r]&0 },
\end{equation}
which we view as a resolution of $\O_Y$ over $Y\times_X Y$. The first observation is the following self-duality of the Koszul complex in the relative setting. We denote by $p\colon Y\times_XY\to X$ the projection, and by $\Wedge^nF^\vee=:L$ the determinant of $F^\vee$.
\begin{Lem}\label{Koszul}
The Koszul complex \eqref{K} and its dual by $(-)^\dagger:=\cHom_{Y\times_XY}(-,p^\ast L)$ are isomorphic, that is, we have a commutative diagram 
\[ \xymatrix@C=5mm{
	0\ar[r]&(\O_Y\otimes_X\O_Y)^\dagger\ar[r]\ar[d]&\cdots\ar[r]&(\O_Y\otimes_X \Wedge^{n-1}F^\vee\otimes_X\O_Y)^\dagger\ar[r]\ar[d]&(\O_Y\otimes_X\Wedge^nF^\vee\otimes\O_Y)^\dagger\ar[r]\ar[d]&0\\
	0\ar[r]&\O_Y\otimes_X\Wedge^nF^\vee\otimes_X\O_Y\ar[r]&\cdots\ar[r]&\O_Y\otimes_X F^\vee\otimes_X\O_Y\ar[r]&\O_Y\otimes_X\O_Y\ar[r]&0. } \]
\end{Lem}
\begin{proof}
	For each $0\leq l\leq n$, we have an isomorphism
	\[
	\begin{aligned}
	\cHom_{Y\times_XY}(\O_Y\otimes_X\Wedge^lF^\vee\otimes_X\O_Y,p^\ast L)
	&\cong\cHom_{X}(\Wedge^lF^\vee,p_\ast p^\ast L)\\
	&\cong\cHom_{X}(\Wedge^lF^\vee,\O_Y\otimes_XL\otimes_X\O_Y)\\
	&\cong\O_Y\otimes_X\cHom_X(\Wedge^lF^\vee,\Wedge^nF^\vee)\otimes_X\O_Y\\
	&\cong\O_Y\otimes_X\Wedge^{n-l}F^\vee\otimes_X\O_Y.
	\end{aligned}
	\]
	The compatibility with differentials can be checked locally, and thus follows from the case of polynomial rings.
\end{proof}

We next translate the above resolution to the algebra side. We need the following compatibility of dualities. We denote by $i\colon Y\times_XY\to Y\times_kY$ the canonical map which is a closed immersion. In what follows, we will simply write $Y\times Y$ for $Y\times_kY$.
\begin{Lem}\label{duals}
The following diagram is commutative.
\[ \xymatrix@C=30mm{
	\rD(Y\times_XY)\ar[r]^-{\rR i_\ast} \ar[d]_-{(-)^\dagger=\RcHom_{Y\times_XY}(-,p^\ast\om_X^{-1})}&\rD(Y\times Y)\ar[r]^-{\RHom_{Y\times Y}(\pi^\ast T\boxtimes \pi^\ast T^\vee,-)}\ar[d]_-{\RcHom_{Y\times Y}(-,\O_{Y\times Y})[d]}&\rD(\Ga^e)\ar[d]^-{\RHom_{\Ga^e}(-,\Ga^e)[d]}\\
	\rD(Y\times_XY)\ar[r]^-{\rR i_\ast}&\rD(Y\times Y)\ar[r]^-{\RHom_{Y\times Y}(\pi^\ast T^\vee\boxtimes \pi^\ast T,-)}&\rD(\Ga^e) } \]
\end{Lem}
\begin{proof}
	Note that the equivalence $\RHom_Y(\pi^\ast T,-)\colon\rD(Y)\to\rD(\Ga)$ induces equivalences 
	\begin{align*}
	&\RHom_{Y\times Y}(\pi^\ast T\boxtimes\pi^\ast T^\vee,-)\colon\rD(Y\times Y)\to\rD(\Ga^e),\\
	&\RHom_{Y\times Y}(\pi^\ast T^\vee\boxtimes \pi^\ast T,-)\colon\rD(Y\times Y)\to\rD(\Ga^e). 
	\end{align*}
	This implies the commutativity of the right square.
	
	We turn to the left square. 
	By Grothendieck duality, we have an isomorphism 
	\[\RcHom_{Y\times Y}(\rR i_\ast(-),\O_{Y\times Y})=\rR i_\ast\RcHom_{Y\times_XY}(-,i^!\O_{Y\times Y}),\] so it remains to show $i^!\O_{Y\times Y}[d]=p^\ast\om_X^{-1}$ in $\rD(Y\times_XY)$.
	Since the square
	\[ \xymatrix{
		Y\times_XY\ar[r]^-i\ar[d]_-p&Y\times Y\ar[d]^-{\pi\times\pi}\\
		X\ar[r]^-\Del&X\times X } \]
	is a pull-back, and the horizontal maps are closed immersions of codimension $d$, we have $i^!\O_{Y\times Y}[d]=\om_{Y\times_XY/Y\times Y}=p^\ast\om_{X/X\times X}=p^\ast(\om_X\otimes_X\Del^\ast\om_{X\times X}^{-1})$. Noting that $\Del^\ast\om_{X\times X}=\om_X^2$, the assertion follows. 
\end{proof}

We next prepare some calculus of bimodules which relates the geometric and the algebraic side.
\begin{Lem}\label{comm}
	Let $X$ be a smooth projective variety of dimension $d$, $T\in\Db(\coh X)$ a generator, and put $\REnd_X(T)=A$. For each  $F\in\rD(\Qcoh X)$, we set
	\[ M(F):=\RHom_X(T,T\lotimes_XF) \in \rD(A^e). \]
	\begin{enumerate}
		\item\label{F} For every $F\in\rD(\Qcoh X)$ we have a commutative diagram
		\[ \xymatrix@C=18mm{
			\rD(\Qcoh X)\ar[r]^-{\RHom_X(T,-)}\ar[d]_-{-\lotimes_XF}&\rD(A)\ar[d]^-{-\lotimes_A M(F)}\\
			\rD(\Qcoh X)\ar[r]^-{\RHom_X(T,-)}&\rD(A) } \]
		\item\label{comp} For $F_1,F_2\in\rD(\Qcoh X)$, we have a canonical isomorphism $M(F_1\lotimes_XF_2)\cong M(F_1)\lotimes_A M(F_2)$ in $\rD(A^e)$.
		\item\label{can} We have an isomorphism $M(\om_X^{-1})\cong\RHom_A(DA,A)[d]$ in $\rD(A^e)$.
	\end{enumerate}
\end{Lem}
\begin{proof}
	(\ref{F})  Since $\RHom_X(T,-)\colon\rD(\Qcoh X)\to\rD(A)$ and $-\lotimes_A T\colon\rD(A)\to\rD(\Qcoh X)$ are mutually inverse equivalences, the functor on $\rD(A)$ corresponding to $-\lotimes_XF$ on $\rD(\Qcoh X)$ is \[\RHom_X(T,-\lotimes_A T\lotimes_XF)=-\lotimes_A\RHom_X(T,T\lotimes_XF).\]

	(\ref{comp})  We have 
	\begin{align*}M(F_1)\lotimes_AM(F_2)&=\RHom_X(T,T\lotimes_XF_1)\lotimes_A\RHom_X(T,T\lotimes_XF_2)\\
	&=\RHom_X(T,\RHom_X(T,T\lotimes_XF_1)\lotimes_AT\lotimes_XF_2)\\
	&=\RHom_X(T,T\lotimes_XF_1\lotimes_XF_2)=M(F_1\lotimes_XF_2).
	\end{align*}

	(\ref{can})  
	We have an isomorphism $\RHom_X(T,T\otimes_X\om_X[d])=D\RHom_X(T,T)=DA$ in $\rD(A^e)$ by Serre duality. 
	It follows that \[\RHom_A(DA,A)[d]\ysimeq\RHom_X(T\otimes_X\om_X,T)=\RHom_X(T,T\otimes_X\om_X^{-1})=M(\om_X^{-1}),\] as desired.
\end{proof}

We are now ready to prove the main statement of this section.
\begin{proof}[Proof of Theorem \ref{CY}]
	(1)  We start from the Koszul bimodule complex
	\[ \xymatrix{ 0\ar[r]&\O_Y\otimes_X\Wedge^rE^\vee\otimes_X\O_Y(-r)\ar[r]&\cdots\ar[r]&\O_Y\otimes_X E^\vee\otimes_X\O_Y(-1)\ar[r]&\O_Y\otimes_X\O_Y\ar[r]&0 }, \]
	where we regard $\O_Y=\Sym_XE^\vee$ as a graded sheaf which induces the grading on $\Ga$.
	Applying $M(-)=\RHom_X(T,T\lotimes_X-)$, each term $\O_Y\otimes_X\Wedge^lE^\vee\otimes_X\O_Y$ becomes 
	\[M(\O_Y)\lotimes_AM(\Wedge^lE^\vee)\lotimes_AM(\O_Y)=\Ga\lotimes_AM(\Wedge^lE^\vee)\lotimes_A\Ga\] by Lemma \ref{comm}(\ref{comp}). We thus obtain a sequence
	\begin{equation}\label{resol}
	\xymatrix{ \Ga\lotimes_AM(\Wedge^rE^\vee)\lotimes_A\Ga(-r)\ar[r]&\cdots\ar[r]&\Ga\lotimes_AM(E^\vee)\lotimes_A\Ga(-1)\ar[r]&\Ga\lotimes_A\Ga }
	\end{equation}
	of Adams graded $\Ga^e$-modules whose total is quasi-isomorphic to $M(\O_Y)=\Ga$.
	
	(2)(3)  By Lemma \ref{Koszul}, the Koszul bimodule complex is isomorphic to its dual by $(-)^\dagger=\RcHom_{Y\times_XY}(-,p^\ast\om_X^{-1})$:
	\[\small \xymatrix@C=5mm{
		0\ar[r]&(\O_Y\otimes_X\O_Y)^\dagger(-r)\ar[r]\ar[d]^-\rsimeq&\cdots\ar[r]&(\O_Y\otimes_X \Wedge^{r-1}E^\vee\otimes_X\O_Y)^\dagger(-1)\ar[r]\ar[d]^-\rsimeq&(\O_Y\otimes_X\Wedge^rE^\vee\otimes\O_Y)^\dagger\ar[r]\ar[d]^-\rsimeq&0\\
		0\ar[r]&\O_Y\otimes_X\Wedge^rE^\vee\otimes_X\O_Y(-r)\ar[r]&\cdots\ar[r]&\O_Y\otimes_X E^\vee\otimes_X\O_Y(-1)\ar[r]&\O_Y\otimes_X\O_Y\ar[r]&0. } \]
    
    Note first that the commutative square
    \[ \xymatrix{
		Y\times_XY\ar[r]^-i\ar[d]_-p&Y\times Y\ar[d]^-{\pi\times\pi}\\
		X\ar[r]^-\Del&X\times X } \]
    shows that the horizontal composite in Lemma \ref{duals} is 
    \[
    \begin{aligned}
    \RHom_{Y\times Y}(\pi^\ast T\boxtimes \pi^\ast T^\vee,\rR i_\ast(-))&=\RHom_{X\times X}(T\boxtimes T^\vee,(\pi\times\pi)_\ast\rR i_\ast(-))\\
    &=\RHom_{X\times X}(T\boxtimes T^\vee,\Del_\ast p_\ast(-))\\
    &=\RHom_X(T\lotimes_X T^\vee,p_\ast(-))\\
    &=\RHom_X(T,T\lotimes_X p_\ast(-)).
    \end{aligned}
    \]
    We now apply
    $\RHom_{Y\times Y}(\pi^\ast T^\vee\boxtimes\pi^\ast T,\rR i_\ast(-))$ to the above isomorphic sequences. Then the above computation and Lemma \ref{comm} show that we get the sequence \eqref{resol} from the second row. Similarly, by Lemma \ref{duals} the first row gives the $d$-shifted bimodule dual over $\Ga$ of \eqref{resol}. We conclude that we have a commutative diagram of double complexes
	\[\small \xymatrix{
		\Ga\lotimes_AA^\ast[d]\lotimes_A\Ga(-r)\ar[r]\ar[d]^-\rsimeq&\cdots\ar[r]&\Ga\lotimes_AM(\Wedge^{r-1}E^\vee)^\ast[d]\lotimes_A\Ga(-1)\ar[r]\ar[d]^-\rsimeq&\Ga\lotimes_AM(\Wedge^rE^\vee)^\ast[d]\lotimes_A\Ga\ar[d]^-\rsimeq\\
		\Ga\lotimes_AM(\Wedge^rE^\vee)\lotimes_A\Ga(-r)\ar[r]&\cdots\ar[r]&\Ga\lotimes_AM(E^\vee)\lotimes_A\Ga(-1)\ar[r]&\Ga\lotimes_A\Ga. } \]
	Looking at their totals, we deduce $\RHom_{\Ga^e}(\Ga,\Ga^e)[d+r]\cong\Ga(r)$ in $\rD^\Z(\Ga^e)$.
\end{proof}

We give some explicit examples of such resolutions. The first one is the fundamental case of $E=\om_X$.
\begin{Ex}
Let $X$ be a smooth projective variety, $E=\om_X$, and $\pi\colon Y=\Spec(\Sym E^\vee)\to X$. Then for each generator $T\in\Db(\coh X)$, the resolution of the CY dg algebra $\Ga=\REnd_Y(\pi^\ast T)$ in Theorem \ref{CY} is
\[	\xymatrix{ \Ga\lotimes_A\RHom_X(T,T\otimes_X\om_X^{-1})\lotimes_A\Ga(-1)\ar[r]&\Ga\lotimes_A\Ga\ar[r]&\Ga}. \]

On the other hand, we know by Theorem \ref{CY} that $\Ga$ is a $(d+1)$-CY dg algebra of Gorenstein parameter $1$. Then by \cite[Theorem 3.9]{Ha7}, such a dg algebra is the $(d+1)$-CY completion of its (Adams) degree $0$ part which is $A=\REnd_X(T)$ in this case; $\Ga\cong\rT^\rL_A(\RHom_A(DA,A)[d])$. For such tensor algebras, we have the standard bimodule resolution
\[	\xymatrix{ \Ga\lotimes_A\RHom_A(DA,A)[d]\lotimes_A\Ga(-1)\ar[r]&\Ga\lotimes_A\Ga\ar[r]&\Ga }. \]
Noting that $\RHom_X(T,T\otimes_X\om_X^{-1})\cong\RHom_A(DA,A)[d]$ by Lemma \ref{comm}(\ref{can}), we see that the resolution from Theorem \ref{CY} is nothing but the standard one for the tensor algebra.
\end{Ex}

The previous example suggests that our bimodule resolution is a higher rank analogue of the one for tensor algebras. Let us describe some of them.
\begin{Ex}
Let $X=\PP^1$ and consider the vector bundle $E=\O(-1)^{\oplus2}$. Then the resolution from Theorem \ref{CY} becomes
\[ \xymatrix{ \Ga\lotimes_A\RHom_{\PP^1}(T,T(2))\lotimes_A\Ga(-2)\ar[r]&\Ga\lotimes_A\RHom_{\PP^1}(T,T(1))^{\oplus2}\lotimes_A\Ga(-1)\ar[r]&\Ga\lotimes_A\Ga\ar[r]&\Ga }. \]
If we take the tilting bundle $T=\O\oplus\O(1)$ as a generator of $\Db(\coh \PP^1)$, the algebra $A$ is the path algebra of the Kronecker quiver
{\raisebox{-0.15cm}{\scalebox{1}{
\begin{tikzpicture}
\newcommand{\arrowwidth}{0.02cm} 
\foreach \a/\b/\n in {0/0/1,2/0/2} {\node (V\a\b) at (\a,\b)  {$\n$}; }; 
\draw[double distance=2pt, scaling nfold=2, -Implies, line width=\arrowwidth] (V00)--(V20);
\end{tikzpicture}}}}.
Now note that the bimodule projective resolutions of $A$, $\RHom_{\PP^1}(T,T(1))=e_2A\oplus\tau^{-1}(e_1A)$, and $\RHom_{\PP^1}(T,T(2))=\tau^{-1}A=\Ext^1_A(DA,A)$ are given by
\[\small \xymatrix@R=3mm{
	0\ar[r]& (Ae_2\otimes e_1A)^{\oplus2}\ar[r]&(Ae_1\otimes e_1A)\oplus(Ae_2\otimes e_2A)\ar[r]&A\ar[r]& 0, \\
	0\ar[r]& Ae_2\otimes e_1A\ar[r]& Ae_1\otimes e_2A\ar[r]&\RHom_{\PP^1}(T,T(1))\ar[r]& 0, \\
	0\ar[r]& (Ae_1\otimes e_1A)\oplus(Ae_2\otimes e_2A)\ar[r]&(Ae_1\otimes e_2A)^{\oplus2}\ar[r]&\RHom_{\PP^1}(T,T(2))\ar[r]& 0. } \]
We then see that the bimodule resolution of $\Ga$ is the total of the following double complex.
\[ \xymatrix{
	(\Ga e_1\otimes e_1\Ga)(-2)\oplus(\Ga e_2\otimes e_2\Ga)(-2)\ar[r]\ar[d]&(\Ga e_2\otimes e_1\Ga)^{\oplus2}(-1)\ar[r]\ar[d]&(\Ga e_2\otimes e_1\Ga)^{\oplus2}\ar[d]\\
	(\Ga e_1\otimes e_2\Ga)^{\oplus2}(-2)\ar[r]&(\Ga e_1\otimes e_2\Ga)^{\oplus2}(-1)\ar[r]&(\Ga e_1\otimes e_1\Ga)\oplus(\Ga e_2\otimes e_2\Ga) } \]
We see that the total is indeed isomorphic to its dual.
\end{Ex}

\section{Calabi--Yau completions of higher rank and projective bundles}
\label{projbundle}
We show that our construction (Theorem \ref{CY}) of CY algebras of higher Gorenstein parameter can be interpreted algebraically by CY completions of higher rank \cite{Ha7} in terms of the associated projective bundle. This generalizes a well-established correspondence for ordinary CY completions \cite{Ke_CYcomp} and canonical bundles.

\begin{Thm}
Let $X$ be a smooth projective variety of dimension $d$, let $E$ be a vector bundle of rank $r$ on $X$ such that $\Wedge^r E=\om_X$ and $E^\vee$ is ample.
Suppose that $T\in\Db(\coh X)$ is a tilting bundle and that $\pi^\ast T\in\Db(\coh Y)$ is also tilting. Put $\Ga=\End_Y(\pi^\ast T)=\bigoplus_{i\geq0}\Hom_X(T,T\otimes_X\Sym^iE^\vee)$.
Let $p\colon Z=\Proj(\Sym E^\vee)\to X$ be the associated projective bundle.
\begin{enumerate}
\item The object $U=\bigoplus_{i=0}^{r-1}p^\ast T(i)\in\Db(\coh Z)$ is tilting, $B=\End_Z(U)$ is $(d+r-1)$-representation infinite, and has a strict $r$-th root pair (see Definition \ref{root}, \ref{strict}).
\item $\Ga$ is the CY completion of the root pair on $B$ in (1).
\end{enumerate}
\end{Thm}
In fact, we show (Theorem \ref{proj}) that the non-commutative projective scheme $\qmod^\Z\!\Ga$ in the sense of Artin--Zhang \cite{AZ} is equivalent to the category $\coh Z$ of coherent sheaves over $Z$, which is an essential observation for the proof.

\subsection{Calabi--Yau completion of higher rank}
Let us now explain the ``higher rank'' analogue which yields Adams graded CY dg algebras of higher Gorenstein parameter.
\begin{Def}[\cite{Ha7}]\label{root}
Let $A$ be a smooth dg algebra and $r$ a positive integer, and $d$ an integer.
\begin{enumerate}
\item A bimodule $U\in\rD(A^e)$ is an {\it $r$-th root} of the $d$-shifted inverse dualizing bimodule $\RHom_{A^e}(A,A^e)[d]$ if $U^{\lotimes_Ar}\cong\RHom_{A^e}(A,A^e)[d]$ in $\rD(A^e)$.
\item A pair $(U,P)$ consisting of $U\in\rD(A^e)$ and $P\in\add A$ is called an {\it $r$-th root pair} if $U$ is an $r$-th root of $\RHom_{A^e}(A,A^e)[d]$ and there is a semi-orthogonal decomposition
\[ \per A=\left\langle \thick P ,\thick(P\lotimes_AU),\ldots,\thick(P\lotimes_AU^{\lotimes_A(r-1)})\right\rangle . \]
\item We denote by $e$ the idempotent of $A$ corresponding to $P$. The {\it Calabi--Yau (CY) completion} of an $r$-th root pair is the dg algebra 
\[ \Pi_{d+1}^{(1/r)}(A;U,P)=\Pi_{d+1}^{(1/r)}(A):=e(\rT^\rL_AU)e, \]
the idempotent truncation of the derived tensor algebra $\rT^\rL_AU$.
\end{enumerate}
\end{Def}
For example, if $r=1$, the conditions force $U=\RHom_{A^e}(A,A^e)[d]$ and $\add P=\add A$, and thus the CY completion of this root pair is the ordinary CY completion.
While any (smooth) dg algebra can be seen as a non-commutative (smooth) variety, we shall see later that each dg algebra admitting a root pair can be regarded as a non-commutative variety which has a structure of a projective bundle.

A prototypical example of a root pair is given by the projective space. 
\begin{Ex}\label{root-ex}
Let $A$ be the $n$-Beilinson algebra; thus it is the endomorphism algebra of a tilting object $T:=\bigoplus_{i=0}^n\O(i)$ over $\PP^n$. Since the derived category $\Db(\mod A)\simeq\Db(\coh\PP^n)$ has a Serre functor $\nu=-\lotimes_ADA\leftrightarrow-\otimes_{\PP^n}\O(-n-1)[n]$, the bimodule $U:=\RHom_{\PP^n}(T,T(1))$ corresponding to $\O(1)$ gives an $(n+1)$-st root of the inverse of $\nu_n=\nu\circ[-n]$.
If we denote by $P$ the projective $A$-module corresponding to the summand $\O$, the semi-orthogonal decomposition
\[ \Db(\coh\PP^n)=\left\langle \thick\O, \thick\O(1),\cdots,\thick\O(n) \right\rangle \]
translates into one
\[ \per A=\left\langle \thick P ,\thick(P\lotimes_AU),\ldots,\thick(P\lotimes_AU^{\lotimes_An})\right\rangle \]
over $A$. This shows that $(U,P)$ is an $(n+1)$-st root pair on $A$. It is easy to see that the CY completion is the polynomial ring $k[x_0,x_1,\ldots,x_n]$.
\end{Ex}

As a consequence of this construction, one indeed gets a CY algebra under a certain invariance property. We refer to \cite[Definition~3.2]{Ha7} for the definition of cyclic invariance.
\begin{Thm}[\cite{Ha7}]\label{root-CY}
Let $A$ be a smooth dg algebra, and let $(U,P)$ be a cyclically invariant $r$-th root pair for $\RHom_{A^e}(A,A^e)[d]$. Then the CY completion $\Pi_{d+1}^{(1/r)}(A)$ is $(d+1)$-CY of Gorenstein parameter $r$.
\end{Thm}

Finally, we recall the following special root pairs on finite dimensional algebras.
\begin{Def}\label{strict}
Let $A$ be a finite dimensional algebra, $r\in\Z_{>0}$ and $d\in\Z$. We say that an $r$-th root pair $(U,P)$ is {\it strict} if it satisfies
\[ \add A=\add\left(\bigoplus_{i=0}^{r-1}P\lotimes_AU^{\lotimes_Ai}\right). \]
\end{Def}
If a $d$-representation infinite algebra has a strict root pair, then its CY completion $\Pi_{d+1}^{(1/r)}(A)=\bigoplus_{i\geq0}\RHom_A(P,P\lotimes_AU^{\lotimes_Ai})$ is concentrated in degree $0$.

\subsection{Non-commutative projective schemes}
In the previous section we discussed a vector bundle $\pi\colon Y=\Spec(\Sym E^\vee)\to X$ over a variety $X$, together with a tilting bundle $T\in\Db(\coh X)$ and its pull-back $\pi^\ast T$ to $Y$.
In this subsection, we consider the projectivization $Z:=\PP(E)=\Proj(\Sym E^\vee)$ of this vector bundle. Let
\[ \xymatrix{ p\colon Z\ar[r]& X } \]
denote the projection. We show (Theorem \ref{proj}) that the category $\coh Z$ can be described as a non-commutative projective scheme of $\Ga=\End_Y(\pi^\ast T)$.

\subsubsection{Artin--Zhang theory}
Let us first recall the notion of non-commutative projective schemes in the sense of Artin--Zhang \cite{AZ}.
Recall that a graded ring $S=\bigoplus_{i\geq0}S_i$ is {\it right graded coherent} if the category $\mod^\Z\!S$ of finitely presented right $S$-modules is abelian. For a right graded coherent $k$-algebra $S=\bigoplus_{i\geq0}S_i$ such that $S_0\in\mod^\Z\!S$, we denote by $\mod_0^\Z\!S$ the full subcategory of $\mod^\Z\!S$ consisting of graded modules $M$ such that $MS_{\geq n}=0$ for sufficiently large $n$. It is a Serre subcategory of $\mod^\Z\!S$, and when $S_0$ is finite dimensional over a field $k$, this is nothing but the category of finite dimensional modules over $S$.
\begin{Def}
For a right graded coherent $k$-algebra $S=\bigoplus_{i\geq0}S_i$ such that $S_0\in\mod^\Z\!S$, we call the Serre quotient
\[ \qmod^\Z\!S:=\mod^\Z\!S/\mod_0^\Z\!S \]
the {\it non-commutative projective scheme} of $S$.
\end{Def}
The non-commutative projective scheme $\qmod^\Z\!S$ is endowed with the {``structure sheaf''} $S\in\qmod^\Z\!S$ and the {``twist functor''} $(1)$ which is an autoequivalence of $\qmod^\Z\!S$.
The following non-commutative analogue of ampleness of line bundles is important.
\begin{Def}\label{AZample}
Let $\C$ be a $k$-linear, $\Hom$-finite abelian category. We say that a pair $(\O,(1))$ consisting of an object $\O\in\C$ and an automorphism $(1)\colon\C\to\C$ is {\it ample} if it satisfies the following conditions.
\begin{enumerate}
\romanenum
\item For every $M\in\C$ there exist $l_1,\ldots,l_p>0$ and an epimorphism $\bigoplus_{i=1}^p\O(-l_i)\to M$ in $\C$.
\item For every epimorphism $M\to N$ in $\C$, the map $\Hom_\C(\O,M(l))\to\Hom_\C(\O,N(l))$ is surjective for every sufficiently large $l$.
\end{enumerate}
\end{Def}

This notion is motivated by the usual notion of ampleness of line bundles on a variety. In fact, the following lemma shows that the non-commutative analogue is much more flexible.
\begin{Lem}\label{Vample}
Let $X$ be a projective variety and $L$ a line bundle on $X$ which is ample. Then for every non-zero vector bundle $V$ on $X$, the pair $(V,-\otimes_XL)$ is ample in $\coh X$ in the sense of Definition \ref{AZample}.
\end{Lem}
\begin{proof}
	Certainly, the usual ampleness implies that the pair $(\O_X,-\otimes_XL)$ is ample in the sense of Definition \ref{AZample}. If $0\neq V\in\vect X$, then for every $M\in\coh X$ there is a surjection $\O_X(-l)^{\oplus}\to \cHom_X(V,M)$ for $l\gg0$. Applying $-\otimes_XV$, we get a surjection $V(-l)^{\oplus}\twoheadrightarrow \cHom_X(V,M)\otimes_XV\twoheadrightarrow M$ when $V\neq0$. The second condition is proved in a similar way as follows. Since $L$ is ample in the usual sense, the epimorphism $V^\vee\otimes_XM\to V^\vee\otimes_XN$ in $\coh X$ yields a surjection $\Hom_X(\O_X,V^\vee\otimes_XM\otimes_XL^l)\to\Hom_X(\O_X,V^\vee\otimes_XN\otimes_XL^l)$ for every sufficiently large $l$. We conclude that $\Hom_X(V,M\otimes_XL^l)\to\Hom_X(V,N\otimes_XL^l)$ is surjective, as desired.
\end{proof}

Now we state Artin--Zhang's Morita-type theorem for non-commutative projective schemes which characterizes the categories $\qmod^\Z\!S$ among abelian categories. 
\begin{Thm}[\cite{AZ}]\label{Morita}
Let $\C$ be a $k$-linear, $\Hom$-finite abelian category. Suppose that there is a pair $(\O,(1))$ which is ample and put
\[ \Ga=\bigoplus_{l\geq0}\Hom_\C(\O,\O(l)), \]
which we naturally regard as a graded algebra.
\begin{enumerate}
\item The algebra $\Ga$ is graded coherent and $\Ga_0\in\mod^\Z\!\Ga$.
\item The functor $\bigoplus_{l\geq0}\Hom_\C(\O,-(l))\colon\C\to\qmod^\Z\!\Ga$ is an equivalence, and is compatible with the degree shift functors.
\end{enumerate}
\end{Thm}

\subsubsection{Our results}
Let us return to our setting from the previous subsection.
Let $X$ be a smooth projective variety of dimension $d$, and let $E$ be a vector bundle of rank $r$ such that $\Wedge^rE=\om_X$ and $E^\vee$ is ample. Let $Y=\Spec(\Sym E^\vee)$ be the total space and $\pi\colon Y\to X$ the projection. 
Also, let $Z=\Proj(\Sym E^\vee)$ and $p\colon Z\to X$ the projectivization, and denote by $\O_Z(1)$ the associated line bundle. 

The following results show that the projective bundle $Z$ gives a geometric realization of the non-commutative projective scheme of $\Ga$. Moreover, it allows us to regard $\Ga$ as the CY completion of a root pair.
\begin{Thm}\label{proj}
Under the above setting, we assume that $X$ has a tilting bundle $T\in\vect X$ such that $\pi^\ast T\in\vect Y$ is still a tilting bundle. We put $\Ga=\End_Y(\pi^\ast T)$, on which we give a grading as in Theorem \ref{CY}.
\begin{enumerate}
	\item\label{ab} There is an equivalence $\coh Z\simeq\qmod^\Z\!\Ga$ taking $p^\ast T$ to $\Ga$ such that the diagram below is commutative.
	\[ \xymatrix{ \coh Z\ar[r]^-\simeq\ar@(dl,dr)[]_-{\O_Z(1)}&\qmod^\Z\!\Ga\ar@(dl,dr)[]_-{(1)} } \]
	In particular, there is a triangle equivalence
	\[ \xymatrix{ \Db(\coh Z)\ar[r]^-\simeq&\Db(\qmod^\Z\!\Ga) }. \]
	\item\label{tri} The variety $Z$ has a tilting bundle $V:=\bigoplus_{i=0}^{r-1}p^\ast T(i)$ whose endomorphism ring $B$ is $(d+r-1)$-representation infinite.
    \item\label{pair} Let $e$ be the idempotent of $B$ corresponding to the summand $p^\ast T$, and put $U=\Hom_Z(V,V(1))\in\mod B^e$. Then $(U,eB)$ is a strict $r$-th root pair on $B$ and $\Ga$ is its CY completion. 
    \end{enumerate}
\end{Thm}

We note the computation of section rings.
\begin{Lem}\label{Ga}
The canonical map
\[ \xymatrix{ \Ga=\disoplus_{l\geq0}\Hom_X(T,T\otimes_X\Sym^lE^\vee)\ar[r]&\disoplus_{l\geq0}\Hom_Z(p^\ast T,p^\ast T(l)) } \]
is an isomorphism of graded algebras.
\end{Lem}
\begin{proof}
	By adjunction, we have $\RHom_Z(p^\ast T,p^\ast T(l))=\RHom_X(T,\rR p_\ast (p^\ast T(l)))$, and $\rR p_\ast(p^\ast T(l))\cong (\rR p_\ast\O_Z(l))\otimes_XT$ by the projection formula. Now, the assertion follows from the isomorphism $\rR p_\ast \O_Z(l)\ysimeq\Sym^lE^\vee$ for every $l\geq0$.
\end{proof}

We now prove the main result of this section.
\begin{proof}[Proof of Theorem \ref{proj}]
	(\ref{ab})  Since $X$ is projective and by the assumption that $E^\vee$ is ample, the line bundle $\O_Z(1)$ is ample. Also, since $T$ is a vector bundle over $X$, so is $p^\ast T$ over $Z$. It follows from Lemma \ref{Vample} that the pair $(p^\ast T, (1))$ on $\coh Z$ is ample. Then by Theorem \ref{Morita} we get an equivalence
	\[ \xymatrix{ \coh Z\ar[r]^-\simeq&\qmod^\Z\!S } \]
	for the graded algebra $S=\bigoplus_{l\geq0}\Hom_Z(p^\ast T,p^\ast T(l))$. Now by Lemma \ref{Ga}, this algebra $S$ is isomorphic to $\Ga=\End_Y(\pi^\ast T)=\bigoplus_{l\geq0}\Hom_X(T,T\otimes_X\Sym^lE^\vee)$ as graded algebras.
	
	(\ref{tri})  Since $\Ga$ is $(d+r)$-CY of Gorenstein parameter $r$ by Theorem \ref{CY}, we know by \cite[Theorem 4.12]{MM} that $\Db(\qmod^\Z\!\Ga)$ has a tilting object $W:=\bigoplus_{i=0}^{r-1}\Ga(i)$ whose endomorphism algebra $B$, described as the matrix algebra
	\[ \begin{pmatrix} \Ga_0&0&\cdots&0 \\ \Ga_1&\Ga_0&\cdots&0\\ \vdots&\vdots&\ddots&\vdots\\ \Ga_{r-1}&\Ga_{r-2}&\cdots&\Ga_0 \end{pmatrix}, \]
	is $(d+r-1)$-representation infinite.
	Along the equivalence $\qmod^\Z\!\Ga\simeq\coh Z$, we deduce that $\bigoplus_{i=0}^{r-1}p^\ast T(i)\in\Db(\coh Z)$ is a tilting object whose endomorphism algebra is $B$.

    (\ref{pair})  By \cite[Corollary 3.10]{Ha7}, the pair $(U,P)$ with $U=\Hom_{\qmod\Ga}^\Z(W,W(1))$ and $P=\Hom_{\qmod\Ga}^\Z(W,\Ga)$ is an $r$-th root pair on $B$ whose CY completion is $\Ga$. We get the assertion via the equivalence $\qmod^\Z\!\Ga\simeq\coh Z$.
\end{proof}

We illustrate how the CY algebra of higher Gorenstein parameter obtained in Theorem \ref{CY} can be regarded as the CY completion of a root pair on the projective bundle.	
\begin{Ex}
Let $X=\PP^1$ and let $E=\O(-1)^{\oplus 2}$. Let $\pi\colon Y=\Spec(\Sym E^\vee)\to X$ be the associated total space as in Setting \ref{general_settingXY}. 

Consider a tilting bundle $T=\O\oplus\O(1)\in\Db(\coh X)$, which pulls back to a tilting bundle $\pi^\ast T\in\Db(\coh Y)$. By Theorem \ref{CY} this gives a $3$-CY algebra of Gorenstein parameter $2$ presented by the following graded quiver with commutativity relations:
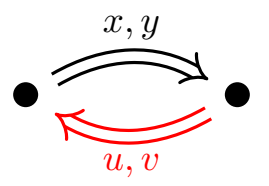
\[ \raisebox{-5mm}{{\scalebox{0.7}{
\begin{tikzpicture}
\newcommand{\arrowwidth}{0.03cm} 

\foreach \a/\b in {0/0,1/0} {\node (V\a\b) at (\a*2,\b*2)  {\huge$\bullet$}; }; 
\draw[double distance=3pt, scaling nfold=2, -Implies, line width=\arrowwidth] (V00) [bend left=30] to (V10);
\draw[double distance=3pt, scaling nfold=2, -Implies, line width=\arrowwidth, red] (V10) [bend left=30] to (V00);
\node at (1,0.65) {$x,y$}; \node[red] at (1,-0.65) {$u,v$}; 
\end{tikzpicture}
}}} \hspace{15pt} \deg x=\deg y=0,\, \deg u=\deg v=1. \]

By Theorem \ref{proj}(\ref{pair}), this can be described as the CY completion of a strict square root pair on an algebra defined as follows.
Let $p\colon Z=\Proj(\Sym E^\vee)\to X$ be the projective bundle associated with $E$. It has a tilting bundle $p^\ast T\oplus p^\ast T(1)\in\Db(\coh Z)$, whose endomorphism algebra $B$ is $2$-representation infinite.
Via an isomorphism $Z\simeq\PP^1\times\PP^1$, we have $p^\ast T\oplus p^\ast T(1)\simeq\O\oplus\O(1,0)\oplus\O(1,1)\oplus\O(2,1)$, which shows that $B$ is presented by the following quiver with commutativity relations.
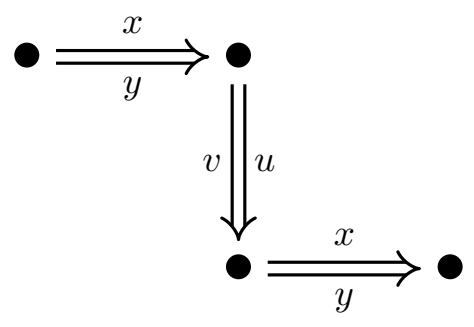
\begin{center}
{\scalebox{0.7}{
\begin{tikzpicture}
\newcommand{\arrowwidth}{0.03cm} 

\foreach \a/\b in {0/0,1/0, 1/-1, 2/-1} {\node (V\a\b) at (\a*2,\b*2)  {\huge$\bullet$}; }; 
\foreach \s/\t in {00/10, 1-1/2-1} {\draw[double distance=3pt, scaling nfold=2, -Implies, line width=\arrowwidth] (V\s)--(V\t) node[midway,xshift=0cm,yshift=-0.3cm] {$y$} node[midway,xshift=0cm,yshift=0.3cm] {$x$};}; 
\foreach \s/\t in {10/1-1} {\draw[double distance=3pt, scaling nfold=2, -Implies, line width=\arrowwidth] (V\s)--(V\t) node[midway,xshift=0.25cm,yshift=0cm] {$u$} node[midway,xshift=-0.25cm,yshift=0cm] {$v$};}; 
\end{tikzpicture}
}}
\end{center}
The square root pair $(U,P)$ on $B$ is given as follows. First, the first two vertices give the projective module $P$, and the square root $U$ of $\RHom_B(DB,B)[2]$ on $B$ is given by $\O_Z(1)$, which corresponds to $\O(1,1)$ on $\Db(\coh\PP^1\times\PP^1)$. We see that this is indeed a strict root pair whose CY completion gives the above CY algebra. 
\end{Ex}

\subsection{Calabi--Yau completions of root pairs}
Finally, we present a geometric counterpart of the construction of CY completion of root pairs.
We start with the following general isomorphism which should be standard. 
\begin{Lem}\label{iso}
Let $X$ be a variety and $E$ a vector bundle with the total space $Y=\Spec(\Sym E^\vee)$. Let $Z=\Proj(\Sym E^\vee)$ be its projectivization. Then the following varieties are canonically isomorphic over $X$.
\begin{enumerate}
\romanenum
\item The blow-up of $Y$ along $X$ (i.e.\,the zero-section).
\item The total space of the line bundle $\O_Z(-1)$ over $Z$.
\end{enumerate}
\end{Lem}
\begin{proof}
	We may assume that $X$ is affine and $E$ is a free module; let $X=\Spec R$ and let $\{x_1,\ldots,x_r\}$ be a basis of $E^\vee$. Then we have $Y=\Spec R[x_1,\ldots,x_r]$, and the blow-up of $Y$ along $X$ is a gluing of $\Spec R[x_1/x_i,\ldots,x_i,\ldots,x_r/x_i]$ for $1\leq i\leq r$.
	On the other hand, we have that $Z$ is a gluing of $\Spec R[x_1/x_i,\ldots,x_r/x_i]$ for $1\leq i\leq r$, and the total space is a gluing of $\Spec  R[x_1/x_i,\ldots,x_r/x_i,x_iy]$ for $1\leq i\leq r$.
	Then we obtain the desired isomorphism from the $R$-algebra isomorphisms
	\[ R\left[ \frac{x_1}{x_i},\ldots,x_i,\ldots,\frac{x_r}{x_i}\right]\cong R\left[\frac{x_1}{x_i},\ldots,\frac{x_r}{x_i},x_iy\right],\qquad \frac{x_j}{x_i}\leftrightarrow \frac{x_j}{x_i}, \ x_i\leftrightarrow x_iy,\]
	which can easily be seen to glue together.
\end{proof}

\begin{Lem}\label{line}
	Let $X$ be a variety, let $T\in\per X$ be a generator, and put $A=\REnd_X(T)$. Let $L$ be a line bundle on $X$ and let $Y=\Spec(\Sym L^\vee)$ be its total space. Then we have a triangle equivalence $\rD(\Qcoh Y)\simeq\rD(\Si)$ for $\Si=\rT^\rL_A(\RHom_X(T,T\otimes_XL^\vee))$.
\end{Lem}
\begin{proof}
	We denote by $\pi\colon Y\to X$ the projection. By Lemma \ref{pullback_generator} we know that $\pi^\ast T$ is a generator of $\per Y$. Then, the adjunction and the projection formula yield 
	\[\REnd_Y(\pi^\ast T)=\RHom_X(T,\pi_\ast\pi^\ast T)=\RHom_X(T,T\otimes_X\pi_\ast\O_Y)=\bigoplus_{l\geq0}\RHom_X(T,T\otimes_X(L^\vee)^{\otimes l}).\] 
	In view of Lemma \ref{comm}(\ref{comp}), the last complex is precisely the derived tensor algebra of $\RHom_X(T,T\otimes_XL^\vee)$ over $\REnd_X(T)$.
\end{proof}

Now we return to the following setting of Theorem \ref{proj}: $X$ is a smooth projective variety and $E$ an anti-ample vector bundle of rank $r$ such that $\Wedge^rE=\om_X$. Consider the total space $Y=\Spec(\Sym E^\vee)\xrightarrow{\pi}X$ and its projectivization $Z=\Proj(\Sym E^\vee)\xrightarrow{p}X$. Moreover suppose that there is a tilting bundle $T$ on $X$ such that $\pi^\ast T$ is a tilting bundle on $Y$, in which case $V=\bigoplus_{i=0}^{r-1}p^\ast T(i)$ is a tilting bundle on $Z$.

Let $B$ be the $(d+r-1)$-representation infinite algebra and $(U,P)$ the $r$-th root pair with $P=eB$ given in Theorem \ref{proj}. Then our CY algebra $\Ga$ is isomorphic to the idempotent truncation $e(\rT^\rL_B\!U)e$ of the (derived) tensor algebra $\rT^\rL_B\!U=\rT_B\!U$.
The following result gives a geometric interpretation of the CY completion of root pairs in our setting.

\begin{Thm}
Under the above setting, we denote by $W$ the varieties in Lemma \ref{iso}, and $\Si=\rT_B\!U$ the tensor algebra. 
\begin{enumerate}
\item There is a triangle equivalence
\[ \xymatrix{ \Db(\coh W)\ar[r]^-\simeq& \Db(\mod \Si) }. \]
\item The following commutative diagrams consisting of derived pull-backs and derived tensor products, respectively, are isomorphic.
\[ \xymatrix{
	\Db(\coh Y)\ar[d]&\Db(\coh X)\ar[l]\ar[d]\\
	\Db(\coh W)&\Db(\coh Z)\ar[l] }
\qquad
\xymatrix{
	\Db(\mod\Ga)\ar[d]&\Db(\mod A)\ar[d]\ar[l]\\
	\Db(\mod\Si)&\Db(\mod B)\ar[l] } \]
\item The vertical functors are fully faithful, and their images are semi-orthogonal components.
\end{enumerate}
\end{Thm}
\begin{proof}
	(1)  By the description (ii) of $W$ in Lemma \ref{iso}, the map $W\to Z$ is the total space of the line bundle $\O_Z(-1)$. We apply Lemma~\ref{line} to this line bundle. We have a tilting bundle $V=\bigoplus_{i=0}^{r-1}p^\ast T(i)$ on $Z$ with $\End_Z(V)=B$, and $\RHom_Z(V,V(1))=U$ as bimodules over $B$. Then the assertion follows from Lemma \ref{line}.
	
	(2)  We have seen that these derived categories are indeed equivalent. We verify the compatibility of the equivalences.
	
	Since $\pi\colon Y\to X$ and $W\to Z$ are vector bundles, it is easy to see that the diagrams
	\[ \xymatrix@R=5mm{
		\Db(\coh X)\ar@{-}[r]^-\simeq\ar[d]&\Db(\mod A)\ar[d]\\
		\Db(\coh Y)\ar@{-}[r]^-\simeq&\Db(\mod\Ga) }
	\qquad
	\xymatrix@R=5mm{
		\Db(\coh Z)\ar@{-}[r]^-\simeq\ar[d]&\Db(\mod B)\ar[d]\\
		\Db(\coh W)\ar@{-}[r]^-\simeq&\Db(\mod\Si) } \]
	are commutative.
	Also, since $p\colon Z\to X$ is a projective bundle and $W\to Y$ is a blow-up, the commutativity of the squares
	\[ \xymatrix@R=5mm{
		\Db(\coh X)\ar@{-}[r]^-\simeq\ar[d]&\Db(\mod A)\ar[d]\\
		\Db(\coh Z)\ar@{-}[r]^-\simeq&\Db(\mod B) }
	\qquad
	\xymatrix@R=5mm{
		\Db(\coh Y)\ar@{-}[r]^-\simeq\ar[d]&\Db(\mod \Ga)\ar[d]\\
		\Db(\coh W)\ar@{-}[r]^-\simeq&\Db(\mod \Si) }
	 \]
	follows from \cite{orlov1993projective}.
	
	(3)  This can be seen either on the geometric or the algebraic side. On the geometric side, the assertion is a consequence of \cite{orlov1993projective}, while on the algebraic side, it is clear since $A$ (resp. $\Ga$) is an idempotent subring of $B$ (resp. of $\Si$). 
\end{proof}

\begin{Ex}
Let $X=\PP^2$ and consider the vector bundle $E=\O(-1)\oplus\O(-2)$.
Consider first the tilting bundle $T=\O\oplus\O(1)\oplus\O(2)$ on $X$ with $A=\End_X(T)$ the $2$-Beilinson algebra presented by the quiver below with commutativity relations.\\
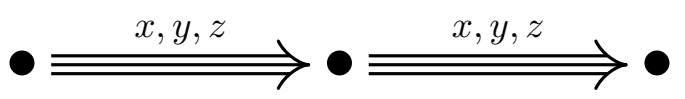
\begin{center}
{\scalebox{0.7}{
\begin{tikzpicture}
\newcommand{\arrowwidth}{0.03cm} 

\foreach \a/\b in {0/0,1/0,2/0} {\node (V\a\b) at (\a*3,\b*3)  {\huge$\bullet$}; }; 
\foreach \s/\t in {00/10,10/20} {\draw[double distance=2pt, scaling nfold=3, -Implies, line width=\arrowwidth] (V\s)--(V\t);}; 
\node at (1.5,0.3) {$x,y,z$}; \node at (4.5,0.3) {$x,y,z$}; 
\end{tikzpicture}
}}
\end{center}

The total space $Y=\Spec(\Sym E^\vee)$ of $E$ over $X$ is derived equivalent to the algebra 
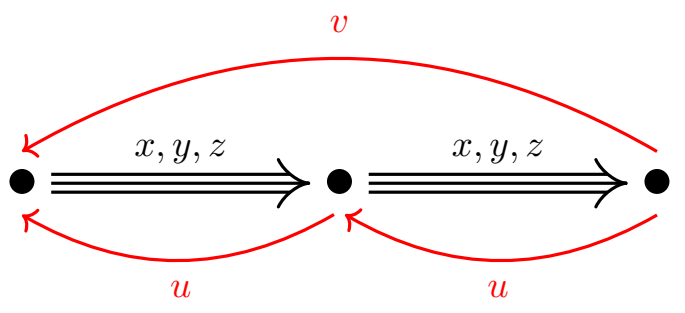
\begin{center}
{\scalebox{0.7}{
\begin{tikzpicture}
\newcommand{\arrowwidth}{0.03cm} 

\foreach \a/\b in {0/0,1/0,2/0} {\node (V\a\b) at (\a*3,\b)  {\huge$\bullet$}; }; 
\foreach \s/\t in {00/10,10/20} {\draw[double distance=2pt, scaling nfold=3, -Implies, line width=\arrowwidth] (V\s)--(V\t);}; 
\node at (1.5,0.3) {$x,y,z$}; \node at (4.5,0.3) {$x,y,z$}; 

\foreach \a/\b in {0/0,2/0} {\path (V\a\b) ++(90:0.3cm) coordinate (V\a\b_n);} 
\draw[->, line width=\arrowwidth, red] (V20_n) [bend right=30] to (V00_n); 
\foreach \a/\b in {0/0,2/0} {\path (V\a\b) ++(270:0.3cm) coordinate (V\a\b_s);} 
\foreach \a in {10} {\path (V\a) ++(280:0.3cm) coordinate (V\a_s1); \path (V\a) ++(260:0.3cm) coordinate (V\a_s2);  
\draw[->, line width=\arrowwidth, red] (V20_s) [bend left=30] to (V\a_s1) ; 
\draw[->, line width=\arrowwidth, red] (V\a_s2) [bend left=30]  to (V00_s) ; };
\node[red] at (1.5,-1) {$u$}; \node[red] at (4.5,-1) {$u$}; \node[red] at (3,1.5) {$v$}; 
\end{tikzpicture}
}}
\end{center}
which is $4$-CY of Gorenstein parameter $2$ by Theorem \ref{CY}.

Next consider the projectivization $p\colon Z=\PP(E)\to X$ of $Y$. By Theorem \ref{proj}(\ref{tri}), the object $V=p^\ast T\oplus p^\ast T(1)$ is a tilting bundle over $Z$. Its endomorphism algebra $B=\End_Z(V)$ is $3$-representation infinite, and is presented by the following quiver with commutativity relations.\\
\begin{center}
{\scalebox{0.7}{
\begin{tikzpicture}
\newcommand{\arrowwidth}{0.03cm} 

\foreach \a/\b in {0/0,1/0,2/0, -1/1,0/1,1/1} {\node (V\a\b) at (\a*3,\b*3)  {\huge$\bullet$}; }; 
\foreach \s/\t in {00/10,10/20, -11/01, 01/11} {\draw[double distance=2pt, scaling nfold=3, -Implies, line width=\arrowwidth] (V\s)--(V\t);}; 
\foreach \s/\t in {1.5/0.3,4.5/0.3, -1.5/3.3, 1.5/3.3} {\node at (\s,\t) {$x,y,z$};}; 

\foreach \s/\t in {01/00,11/10,11/00} {\draw[->, line width=\arrowwidth] (V\s)--(V\t);}; 
\node at (-0.3,1.5) {$u$}; \node at (3.3,1.5) {$u$}; \node at (1.8,1.5) {$v$}; 
\end{tikzpicture}
}}
\end{center}

By Theorem \ref{proj}(\ref{pair}), the pair $(U,P)$ with $P=\Hom_Z(V,p^\ast T)$ and $U=\Hom_Z(V,V(1))$ is a square root pair on $B$.
The tensor algebra $\Si=\rT_B\!U$ is presented by the following quiver. \\
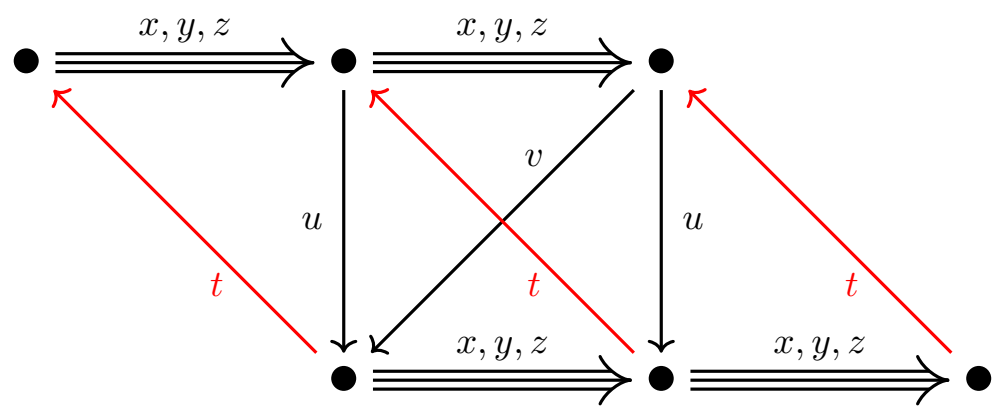
\begin{center}
{\scalebox{0.7}{
\begin{tikzpicture}
\newcommand{\arrowwidth}{0.03cm} 

\foreach \a/\b in {0/0,1/0,2/0, -1/1,0/1,1/1} {\node (V\a\b) at (\a*3,\b*3)  {\huge$\bullet$}; }; 
\foreach \s/\t in {00/10,10/20, -11/01, 01/11} {\draw[double distance=2pt, scaling nfold=3, -Implies, line width=\arrowwidth] (V\s)--(V\t);}; 
\foreach \s/\t in {1.5/0.3,4.5/0.3, -1.5/3.3, 1.5/3.3} {\node at (\s,\t) {$x,y,z$};}; 
\foreach \s/\t in {01/00,11/10,11/00} {\draw[->, line width=\arrowwidth] (V\s)--(V\t);}; 
\node at (-0.3,1.5) {$u$}; \node at (3.3,1.5) {$u$}; \node at (1.8,2.1) {$v$}; 

\foreach \s/\t in {00/-11,10/01,20/11} {\draw[->, line width=\arrowwidth, red] (V\s)--(V\t) node[midway,xshift=0.3cm,yshift=-0.6cm] {$t$};}; 
\end{tikzpicture}
}}
\end{center}

This is derived equivalent to the total space of $\O_Z(-1)$ as well as the blow-up of $Y$ along $X$.

Finally, letting $e$ be the idempotent of $\Si$ corresponding to the three vertices in the first row, we get the CY completion $e\Si e$ of the root pair $(U,P)$ of $B$. We see that this is exactly the CY algebra derived equivalent to $Y$, which confirms Theorem \ref{proj}(\ref{pair}).
\end{Ex}

\section{NCCRs arising from products of projective spaces}

\subsection{Constructions of toric NCCRs from products of projective spaces}

\begin{Setup}
\label{setting_toric_PxP}
Let $k$ be a field, $r,s\geq1$ and $n_1,\ldots,n_s\geq1$ be integers, and consider the polynomial ring
\[ S=k[x_{1,0},\ldots,x_{1,n_1},x_{2,0},\ldots,x_{2,n_2},\ldots,x_{s,0},\ldots,x_{s,n_s},y_1,\ldots,y_r]. \] 
Let the algebraic torus $G\coloneqq(k^\times)^s$ of dimension $s$ act on $S$ by weights
\begin{align*}
\e_i&\coloneqq(0, \dots, 0, \overset{i}{\check{1}}, 0,\dots, 0)\in\Z^s \text{ on } x_{i,j},\, 0\leq j\leq n_i\\
c_l&\coloneqq(c_{1,l}, \dots, c_{s,l})\in\Z^s \text{ on } y_l,
\end{align*}
where we mean that $g=(g_1,\ldots,g_s)\in G$ acts on (the one-dimensional vector space spanned by) $v$ with weight $w=(w_1,\ldots,w_s)\in\Z^s$ when $g\cdot v=g_1^{w_1}\cdots g_s^{w_s}\cdot v$. 
We assume that these weights satisfy $c_{1,l}, \dots, c_{s,l}\le 0$ and
\begin{equation}
	\label{CYcond_weight}
	\sum_{i=1}^s(n_i+1)\e_i+\sum_{l=1}^rc_l=(0,\dots,0), 
\end{equation}
that is, $\sum_{l=1}^rc_l=(-n_1-1, \dots, -n_s-1)$. 
Now let
\[ R\coloneqq S^G \]
be the invariant subring, which is a toric ring of dimension $\sum_{i=1}^sn_i+r$, and the condition \eqref{CYcond_weight} asserts that $R$ is Gorenstein (see e.g., \cite[Subsec.~6.4]{BH_book}). 
\end{Setup}

The purpose of this section is to prove that the toric ring $R$ defined in Setting~\ref{setting_toric_PxP} admits a toric NCCR. 

\begin{Rem}
If $c_l=(-1, \dots, -1)$ for all $l=1, \dots, r$, then we have $n_1=\cdots=n_s=r-1$ by the condition \eqref{CYcond_weight}, 
in which case $R$ is isomorphic to the Segre product of $(s+1)$ polynomial rings of $r$ variables. 
A toric NCCR of such a toric ring was already given in \cite[Corollary~3.11]{Han_NCCR}, \cite[Theorem~3.6]{HigNak_conic}. 
\end{Rem}

To construct an NCCR of $R$, we use the framework given in Corollary~\ref{CY-tilting}. 
First, we consider the product $X\coloneqq \PP^{n_1}\times\cdots\times\PP^{n_s}$ of projective spaces. 

\begin{Lem}
\label{tilting_PxP}
The vector bundle $T\coloneqq\displaystyle\bigoplus_{\substack{0\le a_i\le n_i \\ 1\le i\le s}}\O_X(a_1, \dots, a_{s})$ is a tilting bundle on $X$. 
\end{Lem}

\begin{proof}
By \cite{beilinson_tilting}, we see that $\bigoplus_{a=0}^n\O_{\PP^n}(a)$ is a tilting bundle of $\PP^n$. 
Thus, the assertion follows from \cite[Proposition~3.4]{CostaMiroRoig}. 
\end{proof}

For the vector bundle $E\coloneqq \bigoplus_{l=1}^r\O_X(c_l)$ on $X$, 
we let $Y$ be the total space of $E$ and $\pi\colon Y\rightarrow X$ be the bundle map. 
The next step is to show that $\pi^*T$ is a tilting bundle of $Y$. To show this, we observe cohomology of sheaves on (products of) projective spaces. 

\begin{Lem}[{cf. \cite[Chap.~3, Sect.~5]{hartshorne_AG}}] 
\label{cohomology_Pn}
We have the following: 
\begin{enumerate}
\item $H^p(\PP^n, \O_{\PP^n}(a))=0$ for any integer $a\in\ZZ$ if $p\neq 0, n$; 
\item $H^n(\PP^n, \O_{\PP^n}(a))=0$ if $a\ge -n$. 
\end{enumerate}
\end{Lem}

\begin{Lem}
\label{cohomology_PxP}
Let $X= \PP^{n_1}\times\cdots\times\PP^{n_s}$. We suppose that $a_i\ge -n_i$ for any $i=1, \dots, s$. 
For any integer $p>0$, we have $H^p(X, \O_X(a_1, \dots, a_s))=0$. 
\end{Lem}

\begin{proof}
We show the assertion for the case $s=2$, and the general case follows from the inductive arguments. 
Let $X_i\coloneqq\PP^{n_i}$ $(i=1,2)$. Then, we have 
\begin{align*}
\bigoplus_{p>0} H^p(X_1\times X_2, \O_{X_1\times X_2}(a_1,a_2))&=\bigoplus_{q_1+q_2=p>0}\big(H^{q_1}(X_1, \O_{X_1}(a_1))\otimes H^{q_2}(X_2, \O_{X_2}(a_2))\big)\\
&=\bigoplus_{\substack{(q_1,q_2)=\\ (n_1, 0), (0,n_2) \text{ or } (n_1,n_2)}}\big(H^{q_1}(X_1, \O_{X_1}(a_1))\otimes H^{q_2}(X_2, \O_{X_2}(a_2))\big)=0. 
\end{align*}
Here, the first equation follows from the K\"{u}nneth formula (see e.g. \cite[Proposition 9.2.4]{kempf}), the second one follows from 
Lemma~\ref{cohomology_Pn}(1), and the third one follows from Lemma~\ref{cohomology_Pn}(2). 
\end{proof}

\begin{Thm}
\label{NCCR_toric_PxP}
Let $X= \PP^{n_1}\times\cdots\times\PP^{n_s}$ and $T$ be the tilting bundle as in Lemma~\ref{tilting_PxP}. 
For the vector bundle $E=\bigoplus_{l=1}^r\O_X(c_l)$ on $X$, we consider the total space $Y$ of $E$ and the bundle map $\pi\colon Y\rightarrow X$, 
where $c_l$ is the weight defined in Setting~\ref{setting_toric_PxP}. 
Then, $\pi^\ast T$ is a tilting bundle on $Y$ and $\End_Y(\pi^\ast T)$ is an NCCR of the toric ring $R$ given in Setting~\ref{setting_toric_PxP}. 
\end{Thm}

\begin{proof}
By Lemma~\ref{pullback_generator}, $\pi^\ast T$ is a generator of $\Db(\coh Y)$; thus we only have to show that $\Ext^p_Y(\pi^\ast T, \pi^\ast T)=0$ for any $p>0$. 
Since we have 
\[
\Ext^p_Y(\pi^\ast T, \pi^\ast T)\cong\Ext^p_X(T, \pi_\ast\pi^\ast T)\cong\Ext^p_X(T, T\otimes \Sym E^\vee) 
\]
by the adjunction and the projection formula, we will show $\Ext^p_X(T, T\otimes_X\Sym E^\vee)=0$ for any $p>0$. 

We let $Q\coloneqq\{(a_1, \dots, a_s)\in\RR^s \mid 0\le a_i\le n_i, \, i=1, \dots, s\}\subset\RR^s$, 
and then we can write $T=\bigoplus_{(a_1, \dots, a_s)\in Q\cap\ZZ^s}\O_X(a_1, \dots, a_s)$. 
Also, we consider the Minkowski sum $\widetilde{Q}\coloneqq Q+(-Q)$. 
Then, we see that $\Ext^p_X(T, T\otimes \Sym E^\vee)\cong H^p(X, T^\vee\otimes T\otimes \Sym E^\vee)$, 
where $T^\vee\otimes T\otimes \Sym E^\vee$ is the direct sum of line bundles with the form: 
\begin{equation}\label{eq_linebdl_summand}
\O_X\Big(\widetilde{a}_1-\sum_{l=1}^r m_lc_{1,l}, \, \dots, \, \widetilde{a}_s-\sum_{l=1}^r m_l c_{s,l}\Big)
\end{equation}
for some $m_l\in\ZZ_{\ge 0}$ $(l=1, \dots, r)$ and $(\widetilde{a}_1, \dots, \widetilde{a}_s)\in \widetilde{Q}\cap\ZZ^s$. 
Since $\widetilde{a}_i\ge -n_i$ and $c_{i,l}\le 0$ for any $i=1, \dots, s$ and $l=1, \dots, r$, 
we have $\widetilde{a}_i-\sum_{l=1}^r m_l c_{i,l}\ge -n_i$. 
Thus, we see that $H^p(X, L)=0$ for any line bundle $L$ with the form \eqref{eq_linebdl_summand} by Lemma~\ref{cohomology_PxP}, 
and hence $\pi^\ast T$ is a tilting bundle on $Y$. 

Therefore, $\End_Y(\pi^\ast T)$ is an NCCR of $\Ga(Y,\O_Y)$ by Corollary \ref{CY-tilting}. It is easy to see that the algebra $\Ga(Y,\O_Y)$ is nothing but $R$ defined in Setting \ref{setting_toric_PxP}.
\end{proof}

\begin{Cor}\label{CCR_toric_PxP}
Under the setting of Theorem~\ref{NCCR_toric_PxP}, the canonical map $Y\to\Spec R$ is a crepant resolution.
\end{Cor}
\begin{proof}
This follows from Theorem~\ref{NCCR_toric_PxP} and Corollary~\ref{CY-tilting}.
\end{proof}

\begin{Ex}
Let $s=2$, $n_1=2$, and $n_2=1$. Following Setting~\ref{setting_toric_PxP}, we consider the weights 
\[
w_{1,0}=w_{1,1}=w_{1,2}=(1,0), \quad w_{2,0}=w_{2,1}=(0,1). 
\]
We then set additional weights 
\[c_1=(-1,-2), \quad c_2=(-2,0)\]
which satisfy the condition \eqref{CYcond_weight}. 
Using these weights, we consider the action of $g=(g_1, g_2)\in G=(k^\times)^2$ on $S=k[x_{1,0}, x_{1,1}, x_{1,2}, x_{2,0}, x_{2,1}, y_1, y_2]$ 
defined as 
\begin{align*}
&g\cdot x_{1,j_1}=g_1x_{1,j_1}\,\, (j_1=0,1,2), \hspace{-3cm}&&g\cdot x_{2,j_2}=g_2x_{2,j_2}\,\, (j_2=0,1), \\
&g\cdot y_1=g_1^{-1}g_2^{-2}y_1, &&g\cdot y_2=g_1^{-2}y_2. 
\end{align*}
The invariant subring $R=S^G$ is a $5$-dimensional Gorenstein toric ring, and an NCCR of $R$ can be constructed from 
the tilting bundle $T=\bigoplus_{\substack{a=0, 1, 2 \\ b=0, 1}}\O_X(a,b)$ on $X=\PP^2\times \PP^1$ as shown in Theorem~\ref{NCCR_toric_PxP}. 
Precisely, for the vector bundle $E=\O_X(-1,-2)\oplus\O_X(-2,0)$ and the bundle map $\pi\colon Y=\Spec (\Sym E^\vee)\to X$, 
$\End_Y(\pi^\ast T)$ is an NCCR of $R$. 
The left of Figure~\ref{fig_ex_NCCR_P2P1} is the quiver of $\End_X(T)$ and the right is the quiver of $\End_Y(\pi^\ast T)$. 
In both quivers, the horizontal (resp. vertical) arrows correspond to the multiplication by $x_{1,0}$, $x_{1,1}$, and $x_{1,2}$ (resp. $x_{2,0}$ and $x_{2,1}$). 
In the right quiver, the bent red arrows correspond to multiplication by $y_2$, and the arrows directed to the southwest correspond to 
multiplication by $x_{2,0}y_1$ and $x_{2,1}y_1$. 

\begin{figure}[H]
\begin{center}
{\scalebox{0.8}{
\begin{tikzpicture}
\newcommand{\arrowwidth}{0.03cm} 

\node at (0,0) {
\begin{tikzpicture}
\foreach \a/\b in {0/0,1/0,2/0,0/1,1/1,2/1} {\node (V\a\b) at (\a*3,\b*3)  {$(\a,\b)$};}; 
\foreach \s/\t in {00/10,10/20,01/11,11/21} {\draw[double distance=2pt, scaling nfold=3, -Implies, line width=\arrowwidth] (V\s)--(V\t);}; 
\foreach \s/\t in {00/01,10/11,20/21} {\draw[double distance=3pt, scaling nfold=2, -Implies, line width=\arrowwidth] (V\s)--(V\t);};  
\end{tikzpicture}}; 

\node at (9,0) {
\begin{tikzpicture}
\foreach \a/\b in {0/0,1/0,2/0,0/1,1/1,2/1} {\node (V\a\b) at (\a*3,\b*3)  {$(\a,\b)$};}; 
\foreach \s/\t in {00/10,10/20,01/11,11/21} {\draw[double distance=2pt, scaling nfold=3, -Implies, line width=\arrowwidth] (V\s)--(V\t);}; 
\foreach \s/\t in {00/01,10/11,20/21} {\draw[double distance=3pt, scaling nfold=2, -Implies, line width=\arrowwidth] (V\s)--(V\t);};  

\foreach \s/\t in {11/00,21/10} {\draw[double distance=3pt, scaling nfold=2, -Implies, line width=\arrowwidth, red] (V\s)--(V\t);};

\foreach \a/\b in {0/0,2/0} {\path (V\a\b) ++(270:0.3cm) coordinate (V\a\b_s);} 
\draw[->, line width=\arrowwidth, red] (V20_s) [bend left=15] to (V00_s); 
\foreach \a/\b in {0/1,2/1} {\path (V\a\b) ++(90:0.3cm) coordinate (V\a\b_n);} 
\draw[->, line width=\arrowwidth, red]  (V21_n) [bend right=15] to (V01_n); 
\end{tikzpicture}}; 
\end{tikzpicture}
} }
\caption{}
\label{fig_ex_NCCR_P2P1}
\end{center}
\end{figure}
\end{Ex}

\subsection{NCCRs for some Hibi rings}
\label{subsec_Hibi}

Let $P=\{p_1,\dots,p_{d-1}\}$ be a partially ordered set (poset) equipped with a partial order $\preceq$. 
For $p_i, p_j \in P$, we say that $p_i$ {\em covers} $p_j$ if $p_j \prec p_i$ and there does not exist $p' \in P$ such that $p_j \prec p' \prec p_i$. 
We consider the unique minimal (resp. maximal) element $\widehat{0}$ (resp. $\widehat{1}$) not contained in $P$ 
and let $\widehat{P}\coloneqq P \cup \{\, \widehat{0}, \widehat{1}\, \}$. 
We let $p_0=\widehat{0}$ and $p_d=\widehat{1}$. 
We denote by $H(\widehat{P})$ the \emph{Hasse diagram} of $\widehat{P}$, which is a graph whose set of vertices is $\{p_0,p_1,\dots,p_d\}$ and 
the vertices $p_i$ and $p_j$ are connected by an edge if $p_i$ covers $p_j$ or $p_j$ covers $p_i$. 
We say that $e=\{p_i,p_j\}$, where $0 \leq i \neq j \leq d$, is an \emph{edge} of $\widehat{P}$ if $e$ is an edge of $H(\widehat{P})$. We say that a sequence $(p_{k_1},\cdots,p_{k_m})$ is a \emph{cycle} in $\widehat{P}$ if it forms a cycle in $H(\widehat{P})$. 

For each edge $e=\{p_i,p_j\}$ of $\widehat{P}$ with $p_i \prec p_j$, 
let $\rho_e\in \RR^d$ be defined by 
\begin{align*}
\rho_e\coloneqq
\begin{cases}
\e_i-\e_j, \;&\text{ if }j \not= d, \\
\e_i, \; &\text{ if }j=d,
\end{cases}
\end{align*}
where $\e_0, \dots, \e_{d-1}$ are the unit vectors in $\RR^d$.
Let $\si_P=\mathrm{Cone}(\rho_e \mid e \text{ is an edge of }\widehat{P}) \subset \RR^d$. 
Then, the toric ring
\[ k[P]\coloneqq k[\si_P^\vee \cap \ZZ^d] \]
is called the \emph{Hibi ring} associated with a poset $P$.
We refer to \cite[Section~2]{Hibi} for another description using distributive lattices.
Using Theorem~\ref{NCCR_toric_PxP}, we can construct an NCCR for some Hibi rings.

Following the general theory of toric rings, let us first describe the divisor class group $\Cl (k[P])$.
First, let $\{e_1, \dots, e_n\}$ be the set of edges of $\widehat{P}$. Then $\Cl (k[P])$ is generated by the divisors $D_{e_1},\ldots,D_{e_n}$ corresponding to the rays of $\sigma_P$, which satisfy the relations 
\begin{align}
\label{relation_divisor}
\sum_{q \in U(p)} D_{\{q,p\}}=\sum_{q' \in L(p)} D_{\{p,q'\}} \text{ for }p \in \widehat{P} \setminus \{\,\widehat{0}, \widehat{1}\,\}, 
\quad\text{and}\quad \sum_{q \in U(p_0)}D_{\{q,p_0\}}=0, 
\end{align}
where $U(p)$ (resp. $L(p)$) is the set of all elements in $\widehat{P}$ which cover $p$ (resp. are covered by $p$) 
for each $p \in P \cup\{\,\widehat{0}\,\}$ (resp. $p \in P \cup\{\,\widehat{1}\,\}$). 
It was shown in \cite{HHN_divisor} that the class group $\Cl(k[P])$ is generated by prime divisors corresponding to edges not contained in a spanning tree of $\widehat{P}$. 
In particular, $\Cl(k[P])$ is isomorphic to $\ZZ^{n-d}$, where $d=|P|+1$ and $n$ is the number of edges of $\widehat{P}$. 
Here, a \emph{spanning tree} of $\widehat{P}$ is a set $\{e_{i_1},\cdots,e_{i_d}\}$ of $d$ edges of $\widehat{P}$ 
such that any element in $\widehat{P}$ is an endpoint of some edge in $\{e_{i_1},\cdots,e_{i_d}\}$ and these edges do not form cycles. 
We note that a spanning tree of $\widehat{P}$ is not unique. 

We next give a description of $k[P]$ as a torus invariant. We consider the algebraic torus $G_P\coloneqq \Hom(\Cl(k[P]), k^\times)\cong(k^\times)^{n-d}$ and the character group $\rmX(G_P)$. 
For $i=1, \dots, n$, let $\chi_{e_i}\in\rmX(G_P)$ be the character identified with $D_{e_i}\in\Cl(k[P])$ via an isomorphism $\rmX(G_P)\cong \Cl(k[P])$, 
and hence these characters satisfy the relations \eqref{relation_divisor}. 
We consider the action of $G_P$ on $S\coloneqq k[x_1, \dots, x_n]$ defined by 
\[
g\cdot x_i =\chi_{e_i}(g)x_i \quad \text{for any $g\in G_P$ \,\, ($i=1, \dots, n$)}. 
\]
Then, the invariant subring $S^{G_P}$ is isomorphic to the Hibi ring $k[P]$, see e.g., \cite[Section~3]{HigNak_conic}. 
Note that $S^{G_P}$ is a Cohen-Macaulay domain with $\dim S^{G_P}=d$ $(=|P|+1)$, and $S^{G_P}$ is Gorenstein if and only if $P$ is pure, see \cite[Section~3]{Hibi}. 
Here, $P$ is said to be \emph{pure} if all maximal chains $p_{i_1} \prec \cdots \prec p_{i_\ell}$ of $P$ have the same length. 
For a pure poset $P$ and each $p\in\widehat{P}$, let $d(p)$ be the length of a saturated chain from $\widehat{0}$ to $p$, which is well-defined since $P$ is pure.

\begin{Setup}
\label{setting_poset}
Let $P$ be a pure poset with $d(\,\widehat{1}\,)=r$. Let $V_\ell=\{p\in P \mid d(p)=\ell\}$ for $\ell=1, \dots, r-1$. 
We assume that $P$ satisfies the following conditions: 
\begin{itemize}
\setlength{\parskip}{0pt} 
\setlength{\itemsep}{3pt}
\item[(a)] For each $\ell$, the set $V_\ell$ contains at most one element whose valency in $H(\widehat{P})$ is greater than two. 
\item[(b)] There exists a maximal chain $\gamma\colon \widehat{0}\prec q_1\prec \cdots\prec q_{r-1}\prec \widehat{1}$ in $\widehat{P}$ with $q_\ell\in V_\ell$ such that  the valency of each vertex in $V_\ell\setminus\{q_\ell\}$ is at most two.
\end{itemize}
Note that if $V_\ell$ contains an element whose valency is greater than two, then it is necessarily contained in the chain $\gamma$, and it is denoted by $q_\ell$. 
\end{Setup}

In the following, we consider a pure poset $P$ and a maximal chain $\ga$ defined in Setting~\ref{setting_poset}. 
We denote the set of edges appearing in $\ga$ by the same letter $\ga$. 
Then, we take a spanning tree $\calT(P)$ of $\widehat{P}$ determined by the following rule. 
For $p\in P\cup\{\,\widehat{0}\,\}$, we set the subset $E(p)\coloneqq\{\{p,p^\prime\} \mid p^\prime\in U(p)\}$ of edges. 

\begin{itemize}
\setlength{\parskip}{0pt} 
\setlength{\itemsep}{5pt}
\item If $\gamma\cap E(p)\neq \varnothing$, in which case $|\gamma\cap E(p)|=1$ by definition, 
then we include the edge $e\in \gamma\cap E(p)$ in $\calT(P)$ and exclude the remaining edges in $E(p){\setminus}\{e\}$ from $\calT(P)$. 
\item If $\gamma\cap E(p)= \varnothing$, in which case $|E(p)|=1$ by definition, then we include the edge $e\in E(p)$ in $\calT(P)$. 
\end{itemize}

Let us consider the Hibi ring $k[P]$ associated to the pure poset $P$ given in Setting~\ref{setting_poset}. 
As we mentioned above, the class group $\Cl(k[P])$ (the character group $\rmX(G_{P})$) is generated 
by the divisors (characters) corresponding to edges of $\widehat{P}$ not contained in a spanning tree. 
Thus, we now use the spanning tree $\calT(P)$ to consider these groups. 
Let $n$ be the number of edges of $\widehat{P}$, and let $d+1=|\widehat{P}|$. 
Then, we have the edges $e_1, \dots, e_{n-d}$ of $\widehat{P}$ not contained in $\calT(P)$. 
In the following, we identify a character in $\rmX(G_{P})$ with the corresponding weight via an isomorphism 
$\rmX(G_{P})\cong \ZZ^{n-d}$. 
Since the characters $\chi_{e_1}, \dots, \chi_{e_{n-d}}$ are generators of $\rmX(G_{P})$, 
we let $\chi_{e_i}=(0,\dots,0,\overset{i}{\check{1}},0,\dots,0)\in \ZZ^{n-d}$ for $i=1, \dots, n-d$. 
Then, for any edge $e$ in $\widehat{P}$, the character $\chi_e$ is determined by the relations \eqref{relation_divisor} among prime divisors. 

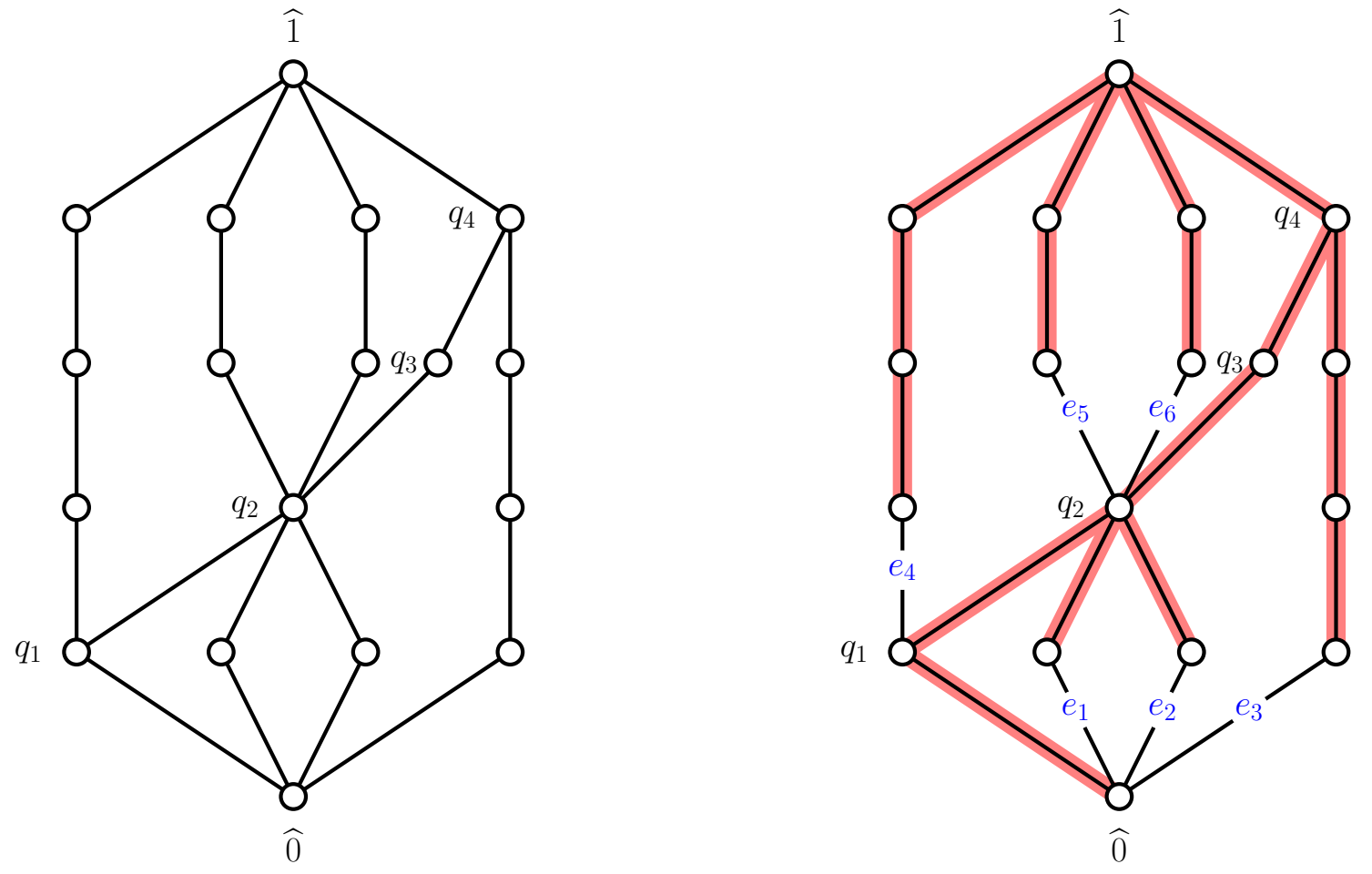
\begin{figure}[H]
\begin{center}
\newcommand{\edgewidth}{0.08cm} 
\newcommand{\nodewidth}{0.08cm} 
\newcommand{\noderad}{0.26} 
\newcommand{\lettersize}{\Huge} 

\begin{tikzpicture}

\node at (0,0) 
{\scalebox{0.35}{
\begin{tikzpicture}
\foreach \a/\b in {1/1,1/2,1/3,1/4,2/1,2/3,2/4,3/1,3/3,3/4,4/1,4/2,4/3,4/4} {\coordinate (V\a\b) at (\a*3,\b*3); };
\coordinate (Min) at (7.5,0); \coordinate (Max) at (7.5,15); \coordinate (Min-) at (7.5,-1); \coordinate (Max+) at (7.5,16);
\coordinate (Q2) at (2.5*3,2*3); \coordinate (Q3) at (3.5*3,3*3); 

\foreach \s/\t in {11/12,12/13,13/14,23/24,33/34,41/42,42/43,43/44} {\draw[line width=\edgewidth]  (V\s)--(V\t);};
\foreach \t in {1,2,3,4} {\draw[line width=\edgewidth]  (Min)--(V\t1);};
\foreach \s in {1,2,3,4} {\draw[line width=\edgewidth]  (V\s4)--(Max);}; 
\foreach \s/\t in {1/1,2/1,3/1,2/3,3/3} {\draw[line width=\edgewidth] (Q2)--(V\s\t); };
\draw[line width=\edgewidth] (Q2)--(Q3); \draw[line width=\edgewidth] (Q3)--(V44); 

\foreach \a/\b in {1/1,1/2,1/3,1/4,2/1,2/3,2/4,3/1,3/3,3/4,4/1,4/2,4/3,4/4} {
\draw [line width=\nodewidth, fill=white] (V\a\b) circle [radius=\noderad]; };
\draw [line width=\nodewidth, fill=white] (Min) circle [radius=\noderad]; \draw [line width=\nodewidth, fill=white] (Max) circle [radius=\noderad]; 
\draw [line width=\nodewidth, fill=white] (Q2) circle [radius=\noderad]; \draw [line width=\nodewidth, fill=white] (Q3) circle [radius=\noderad]; 

\coordinate (Q1-) at (1*3-1,1*3); \coordinate (Q2-) at (2.5*3-1,2*3); 
\coordinate (Q3-) at (3.5*3-0.7,3*3); \coordinate (Q4-) at (4*3-1,4*3); 
\foreach \a in {1,2,3,4} {\node at (Q\a-) {\lettersize$q_{\a}$}; };
\node at (Min-) {\lettersize$\widehat{0}$}; \node at (Max+) {\lettersize$\widehat{1}$};
\end{tikzpicture}
} } ;

\node at (6,0) 
{\scalebox{0.35}{
\begin{tikzpicture}
\newcommand{\treewidth}{0.4cm} 
\newcommand{\treecolor}{red!50} 

\foreach \s/\t in {12/13,13/14,23/24,33/34,41/42,42/43,43/44} {\draw[line width=\treewidth, \treecolor]  (V\s)--(V\t);};
\foreach \s in {11,21,31} {\draw[line width=\treewidth, \treecolor]  (V\s)--(Q2);};
\foreach \s in {14,24,34,44} {\draw[line width=\treewidth, \treecolor]  (V\s)--(Max);}; 
\draw[line width=\treewidth, \treecolor]  (Q3)--(V44); \draw[line width=\treewidth, \treecolor]  (Q2)--(Q3); \draw[line width=\treewidth, \treecolor]  (V11)--(Min);
\foreach \s/\t in {11/12,12/13,13/14,23/24,33/34,41/42,42/43,43/44} {\draw[line width=\edgewidth]  (V\s)--(V\t);};
\foreach \t in {1,2,3,4} {\draw[line width=\edgewidth]  (Min)--(V\t1);};
\foreach \s in {1,2,3,4} {\draw[line width=\edgewidth]  (V\s4)--(Max);}; 
\foreach \s/\t in {1/1,2/1,3/1,2/3,3/3} {\draw[line width=\edgewidth] (Q2)--(V\s\t); };
\draw[line width=\edgewidth] (Q2)--(Q3); \draw[line width=\edgewidth] (Q3)--(V44); 

\foreach \a/\b in {1/1,1/2,1/3,1/4,2/1,2/3,2/4,3/1,3/3,3/4,4/1,4/2,4/3,4/4} {
\draw [line width=\nodewidth, fill=white] (V\a\b) circle [radius=\noderad]; };
\draw [line width=\nodewidth, fill=white] (Min) circle [radius=\noderad]; \draw [line width=\nodewidth, fill=white] (Max) circle [radius=\noderad]; 
\draw [line width=\nodewidth, fill=white] (Q2) circle [radius=\noderad]; \draw [line width=\nodewidth, fill=white] (Q3) circle [radius=\noderad]; 

\coordinate (Q1-) at (1*3-1,1*3); \coordinate (Q2-) at (2.5*3-1,2*3); 
\coordinate (Q3-) at (3.5*3-0.7,3*3); \coordinate (Q4-) at (4*3-1,4*3); 
\foreach \a in {1,2,3,4} {\node at (Q\a-) {\lettersize$q_{\a}$}; };
\node at (Min-) {\lettersize$\widehat{0}$}; \node at (Max+) {\lettersize$\widehat{1}$};
\foreach \a/\b/\n in {6.6/1.8/e_1, 8.4/1.8/e_2, 10.2/1.8/e_3,3/4.7/e_4, 6.6/8/e_5, 8.4/8/e_6}{
\filldraw[white] (\a,\b) circle [radius=0.4]; \node[blue] at (\a,\b) {\lettersize$\n$}; };
\end{tikzpicture}
} } ;
\end{tikzpicture}
\end{center}
\caption{An example of $P$ satisfying Setting~\ref{setting_poset} and its spanning tree $\calT(P)$}
\label{Fig_Segre_54}
\end{figure}

\begin{Ex}
\label{ex_sptree}
We consider the poset $P$ whose Hasse diagram takes the form as described in the left of Figure~\ref{Fig_Segre_54}, 
which has a maximal chain $\gamma\colon \widehat{0}\prec q_1\prec q_2\prec q_3\prec q_4\prec\widehat{1}$ satisfying the conditions in Setting~\ref{setting_poset}. 
In this case, the spanning tree $\calT(P)$ consists of the colored edges in the right of Figure~\ref{Fig_Segre_54}. 
Thus, the edges $e_1, \dots, e_6$ correspond to generators $\chi_{e_1}, \dots, \chi_{e_6}$ of $\rmX(G_P)$. 
In particular, using the relations \eqref{relation_divisor}, we have 
\begin{align}\label{ex_char}
&\chi_{\{\widehat{0},q_1\}}=(-1,-1,-1,0,0,0), &&\chi_{\{q_1,q_2\}}=(-1,-1,-1,-1,0,0), \\ \nonumber
&\chi_{\{q_2,q_3\}}=\chi_{\{q_3,q_4\}}=(0,0,-1,-1,-1,-1), &&\chi_{\{q_4, \widehat{1}\}}=(0,0,0,-1,-1,-1).
\end{align}
\end{Ex}

\begin{Lem}
\label{lem_weight_Hibi}
Let the notation be the same as above. For the character $\chi_e$ corresponding to an edge $e$ in $\widehat{P}$, 
we have that 
\begin{itemize}
\setlength{\parskip}{0pt} 
\setlength{\itemsep}{5pt} 
\item if $e$ is not contained in $\gamma$, then $\chi_e=(0,\dots,0,\overset{i}{\check{1}},0,\dots,0)\in\ZZ^{n-d}$ for some $i=1, \dots, n-d$, 
\item if $e$ is contained in $\gamma$, then $\chi_e\in\ZZ^{n-d}_{\le 0}$. 
\end{itemize}
\end{Lem}

\begin{proof}
We set \[E_j\coloneqq \bigcup_{\substack{p\in\widehat{P}{\setminus}\{\,\widehat{1}\,\} \\ d(p)=j-1}}E(p)\]
for $j=1, \dots, r$. We note that 
\begin{equation}
\label{relation_character_j}
\sum_{e\in E_j}\chi_e=(0, \dots, 0)\in\ZZ^{n-d}
\end{equation}
holds for any $j=1, \dots, r$ by \eqref{relation_divisor}. Then, we show the assertion for $e\in E_j$ by induction on $j$. 

For the case $j=1$, we see that $E_1=E(\,\widehat{0}\,)$ and $|\gamma\cap E_1|=1$. Let $\{e^1_\gamma\}\coloneqq \gamma\cap E_1\subset \calT(P)$. 
Since all edges in $E_1{\setminus}\{e^1_\gamma\}$ are not contained in $\calT(P)$ by definition, 
they are a part of the edges $\{e_1, \dots, e_{n-d}\}$, and hence they correspond to characters with the form $(0,\dots,0,\overset{i}{\check{1}},0,\dots,0)$ for some $i$. 
By \eqref{relation_character_j}, we have 
\[
\chi_{e^1_\gamma}+\sum_{e\in E_1{\setminus}\{e^1_\gamma\}}\chi_{e}=(0, \dots, 0), 
\]
thus $\chi_{e^1_\gamma}\in\ZZ_{\le 0}^{n-d}$. 

For $j\geq1$, we assume that the assertion holds for any edge in $E_j$, and we consider $E_{j+1}$. 
Since $|\gamma\cap E_{j+1}|=1$, we let $\{e^{j+1}_\gamma\}\coloneqq \gamma\cap E_{j+1}$, in which $e^{j+1}_\gamma=\{q_j, q_{j+1}\}$. 
If an edge $e\in E_{j+1}{\setminus}\{e^{j+1}_\gamma\}$ does not share the endpoint $q_j$ with $e^{j+1}_\gamma$, 
then $\chi_e=\chi_{e^\prime}$ for some edge $e^\prime\in E_j$ not contained in $\gamma$. 
By assumption, $\chi_{e^\prime}$ takes the form of a unit vector, and hence so does $\chi_e$. 
If an edge $e\in E_{j+1}{\setminus}\{e^{j+1}_\gamma\}$ shares the endpoint $q_j$ with $e^{j+1}_\gamma$, then $e$ is not contained in $\calT(P)$, 
which means $e$ is one of the edges in $\{e_1, \dots, e_{n-d}\}$ and $\chi_e$ takes the form of a unit vector. 
Thus, for any edge $e$ in $E_{j+1}{\setminus}\{e^{j+1}_\gamma\}$, we see that $\chi_e=(0,\dots,0,\overset{i}{\check{1}},0,\dots,0)$ for some $i$. 
By \eqref{relation_character_j}, we have
\[
\chi_{e^{j+1}_\gamma}+\sum_{e\in E_{j+1}{\setminus}\{e^{j+1}_\gamma\}}\chi_{e}=(0, \dots, 0), 
\]
and hence $\chi_{e^{j+1}_\gamma}\in\ZZ^{n-d}_{\le 0}$. 
\end{proof}

Under these settings, we can show the existence of a toric NCCR for the Hibi ring $k[P]$. 

We define the poset ${P_\ga}$ as follows:
\begin{itemize}
\item ${P}_\ga=\{p\in \widehat{P}\mid p \text{ is not an endpoint of an edge in }\gamma\}$ as a set,
\item for $x,y\in{P_\ga}$, let $x\preceq y$ if there is a sequence of cover relations $x\preceq z_1\preceq\cdots\preceq z_m\preceq y$ from $x$ to $y$ in $P$ such that $z_i\in P_\ga$.
\end{itemize}
Note that the Hasse diagram of ${P_\ga}$ is obtained from that of $\widehat{P}$ by deleting the vertices contained in $\ga$ and their adjacent edges.

\begin{Thm}
\label{NCCR_Hibi}
Let $P$, $\ga$ be a pure poset and a maximal chain defined in Setting~\ref{setting_poset}.
\begin{enumerate}
\item ${P_\ga}$ is a disjoint union of some totally ordered sets.
\item The number of the connected components of ${P_\ga}$ is $n-d$.
\item 
Let $s\coloneqq n-d$, and $P_1, \dots, P_s$ be the connected components of ${P_\ga}$, which are totally ordered sets. 
Let $n_i\coloneqq |P_i|$ for $i=1, \dots, s$.
Put 
\[X_\ga=\Proj(k[P_1])\times\cdots\times\Proj(k[P_s])=\PP^{n_1}\times\cdots\times\PP^{n_s},\] 
and let $E_\ga=\bigoplus_{e\in\ga}\O_{X_\ga}(\chi_e)$ be a vector bundle over $X_\ga$. Then $\rank E_\ga=r$ and $\Wedge^rE_\ga=\om_{X_\ga}$.
\item Let $Y_\ga$ be the total space of $E_\ga$. Then $k[P]$ is isomorphic to the global section ring of $Y_\ga$. 
Therefore the Hibi ring $k[P]$ has a toric NCCR. 
\end{enumerate}
\end{Thm}

\begin{proof}
(1)  By our choice of $\ga$, each element in $P_\ga$ covers (resp. is covered by) at most one element in $P_\ga$. Then the assertion is a consequence of the following observation: If $\La$ is a connected poset such that each element covers (resp. is covered by) at most one element, then $\La$ is totally ordered.  

(2) 
By the definition of ${P_\ga}$ and Lemma~\ref{lem_weight_Hibi}, 
for any edge $e$ of $\widehat{P}$ with at least one endpoint contained in ${P_\ga}$, we have $\chi_e=(0,\dots,0,\overset{i}{\check{1}},0,\dots,0)\in\ZZ^{n-d}$ for some $i=1, \dots, n-d$. 
By the relations \eqref{relation_divisor}, the characters corresponding to edges that have an endpoint contained in the same connected component of ${P_\ga}$ are identical. 
Thus, we have a bijection between the connected components of ${P_\ga}$ and generators of $\rmX(G_{P})\cong \ZZ^{n-d}$.

(3)  
It is clear that $\rank E_\ga=|\ga|$, which is equal to the length $r$ of a maximal chain of $\widehat{P}$. 
Let $E(P_i)=\{e : \text{an edge of } \widehat{P} \mid \text{at least one endpoint of $e$ is contained in $P_i$} \}$ for $i=1, \dots, s$. 
By the definition of ${P_\ga}$, the minimal (resp. maximal) element of $P_i$ covers (is covered by) a single element of $\widehat{P}$ appearing in $\ga$. Thus, we have $|E(P_i)|=n_i+1$ for any $i=1, \dots, s$. 
Let $e_1, \dots, e_{n-d}$ be the edges of $\widehat{P}$ not contained in the spanning tree $\calT(P)$. 
By the argument in (2), we may assume that $e_i\in E(P_i)$, in which case 
the characters corresponding to the edges in $E(P_i)$ are identical to $\chi_{e_i}=(0,\dots,0,\overset{i}{\check{1}},0,\dots,0)\in\ZZ^{n-d}$ 
for any $i=1, \dots, s$. 
Since the sum of all characters corresponding to edges of $\widehat{P}$ is trivial by \eqref{relation_character_j}, 
we have $\sum_{e\in\ga}\chi_e=(-n_1-1,\ldots,-n_s-1)$ and hence $\Wedge^rE_\ga=\O_{X_\ga}(-n_1-1,\ldots,-n_s-1)=\om_{X_\ga}$.

(4)  The Hibi ring $k[P]$ is isomorphic to the invariant subring under the action of $G_{P}$ defined by the (weight of) characters observed in Lemma~\ref{lem_weight_Hibi}. 
Thus, $k[P]$ is a special case of Setting~\ref{setting_toric_PxP}, and the assertion follows from Theorem~\ref{NCCR_toric_PxP}. 
\end{proof}

\begin{Rem}
Note that Theorem~\ref{NCCR_Hibi} recovers \cite[Theorem~3.6]{HigNak_conic} and \cite[Theorem~3.3]{Nak_HibiNCCR}, 
since the posets discussed in these papers satisfy the conditions in Setting~\ref{setting_poset}. 
\end{Rem}

\begin{figure}[H]
\begin{center}
\newcommand{\edgewidth}{0.08cm} 
\newcommand{\nodewidth}{0.08cm} 
\newcommand{\noderad}{0.26} 
\newcommand{\lettersize}{\Huge} 
\newcommand{\ciclewidth}{0.06cm} 

{\scalebox{0.35}{
\begin{tikzpicture}
\newcommand{\treewidth}{0.4cm} 
\newcommand{\treecolor}{red!50} 
\foreach \a/\b in {1/1,1/2,1/3,1/4,2/1,2/3,2/4,3/1,3/3,3/4,4/1,4/2,4/3,4/4} {\coordinate (V\a\b) at (\a*3,\b*3); };
\coordinate (Min) at (7.5,0); \coordinate (Max) at (7.5,15); \coordinate (Min-) at (7.5,-1); \coordinate (Max+) at (7.5,16);
\coordinate (Q2) at (2.5*3,2*3); \coordinate (Q3) at (3.5*3,3*3); 

\foreach \s/\t in {12/13,13/14,23/24,33/34,41/42,42/43} {\draw[line width=\edgewidth]  (V\s)--(V\t);};
\foreach \t in {1,2,3,4} {\draw[line width=\edgewidth,lightgray]  (Min)--(V\t1);};
\foreach \s in {1,2,3,4} {\draw[line width=\edgewidth,lightgray]  (V\s4)--(Max);}; 
\foreach \s/\t in {11/12,43/44} {\draw[line width=\edgewidth,lightgray]  (V\s)--(V\t);};
\foreach \s/\t in {1/1,2/1,3/1,2/3,3/3} {\draw[line width=\edgewidth,lightgray] (Q2)--(V\s\t); };
\draw[line width=\edgewidth,lightgray] (Q2)--(Q3); \draw[line width=\edgewidth,lightgray] (Q3)--(V44); 

\foreach \a/\b in {1/2,1/3,1/4,2/1,2/3,2/4,3/1,3/3,3/4,4/1,4/2,4/3} {
\draw [line width=\nodewidth, fill=white] (V\a\b) circle [radius=\noderad]; };
\foreach \a/\b in {1/1,4/4} {\draw [line width=\nodewidth, lightgray, fill=white] (V\a\b) circle [radius=\noderad]; };
\draw [line width=\nodewidth, lightgray, fill=white] (Min) circle [radius=\noderad]; 
\draw [line width=\nodewidth, lightgray, fill=white] (Max) circle [radius=\noderad]; 
\draw [line width=\nodewidth, lightgray, fill=white] (Q2) circle [radius=\noderad]; 
\draw [line width=\nodewidth, lightgray, fill=white] (Q3) circle [radius=\noderad]; 

\coordinate (Q1-) at (1*3-1,1*3); \coordinate (Q2-) at (2.5*3-1,2*3); 
\coordinate (Q3-) at (3.5*3-0.7,3*3); \coordinate (Q4-) at (4*3-1,4*3); 
\foreach \a in {1,2,3,4} {\node[lightgray] at (Q\a-) {\lettersize$q_{\a}$}; };
\node[lightgray] at (Min-) {\lettersize$\widehat{0}$}; \node[lightgray] at (Max+) {\lettersize$\widehat{1}$};

\draw[red,line width=\ciclewidth] (V21) circle[radius=1.1]; \draw[red,line width=\ciclewidth] (V31) circle[radius=1.1]; 
\draw[red,line width=\ciclewidth] (V13) circle (1 and 4); \draw[red,line width=\ciclewidth] (V42) circle (1 and 4);
\draw[red,line width=\ciclewidth] (6,10.5) circle (1 and 2.6); \draw[red,line width=\ciclewidth] (9,10.5) circle (1 and 2.6); 
\foreach \a/\b/\n in {4.8/3/P_1, 7.8/3/P_2, 11/6/P_3, 2/9/P_4, 5/10.5/P_5, 8/10.5/P_6}{
\filldraw[white] (\a,\b) circle [radius=0.55]; \node[red] at (\a,\b) {\lettersize$\n$}; };
\end{tikzpicture}
} } 

\end{center}
\caption{The Hasse diagram of ${P}_\gamma$ for $P$ and $\gamma$ given in Example~\ref{ex_sptree}}
\label{Fig_Pgamma}
\end{figure}
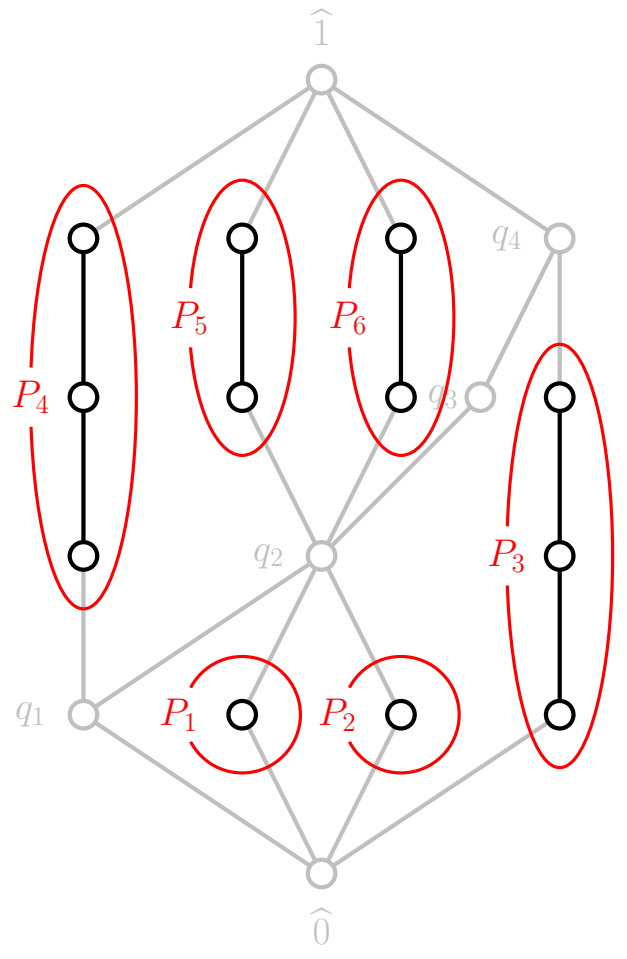

\begin{Ex}
We consider the pure poset $P$ and the maximal chain $\gamma$ given in Example~\ref{ex_sptree}. 
Then, the poset ${P}_\gamma=P_1\sqcup \cdots \sqcup P_6$,  where $P_1, \dots, P_6$ are totally ordered sets 
whose Hasse diagrams take the forms as in Figure~\ref{Fig_Pgamma}. 
In this case, we put 
\[
X_\ga=\Proj(k[P_1])\times\cdots\times\Proj(k[P_6])=\PP^1\times\PP^1\times\PP^3\times\PP^3\times\PP^2\times\PP^2, 
\]
and consider the vector bundle $E_\ga=\bigoplus_{e\in\ga}\O_{X_\ga}(\chi_e)$ where 
$\chi_e$ $(e\in \ga)$ are the characters given in \eqref{ex_char}. 
Then, the Hibi ring $k[P]$ is isomorphic to the section ring of the total space of $E_\ga$, and it has a toric NCCR by Theorem~\ref{NCCR_Hibi}. 
\end{Ex}

\section{Calabi--Yau algebras from quotient stacks}\label{quotstack}
We give a variant of Theorem \ref{CY} for quotient stacks. 

Let $A$ be a finitely generated (commutative) $k$-algebra graded by a finitely generated abelian group $L$. We fix a homogeneous ideal $I\subset A$ and consider the quotient category
\[ \qmod^L\!A:=\mod^L\!A/\mod^L_I\!A \]
of the category of finitely generated $L$-graded modules $\mod^L\!A$ by its full subcategory $\mod_I^L\!A$ consisting of modules supported on $I$. The category $\qmod^L\!A$ should be regarded as the category of coherent sheaves over the quotient stack
\[ \X=[(\Spec A\setminus V(I))/\Spec k[L]], \]
where $V(I)=\Spec A/I$, thus $\coh\X=\qmod^L\!A$.
Similarly, consider the large version
\[ \Qcoh\X:=\QMod^L\!A:=\Mod^L\!A/\Mod^L_I\!A, \]
the quotient of the category of all $L$-graded $A$-modules by the category of $I$-torsion modules.

We denote by $\O_\X\in\coh\X$ the object corresponding to $A\in\qmod^L\!A$. Also, we denote the morphism spaces in $\QMod^L\!A$ by $\Hom_\X(-,-)$, and similarly for the derived functors.
Since $A$ is commutative, we have a bifunctor $-\otimes_A-\colon\Mod^L\!A\times\Mod^L\!A\to\Mod^L\!A$. This induces a bifunctor on the quotient $\QMod^L\!A$, which is denoted by $-\otimes_\X-$. Similarly, the functor induced by $\Hom_A(-,-)\colon(\mod^L\!A)^\op\times\Mod^L\!A\to\Mod^L\!A$ is denoted by $\cHom_\X(-,-)$.

In the rest of this section, we assume the following on our algebra $A$ and the ideal $I$.
\begin{Assum}\label{stack}
Let $A$ be an $L$-graded algebra, and let $I\subset A$ be a homogeneous ideal.
\begin{enumerate}
\alphenum
\item $A=k[x_1,\ldots,x_n]$ is an $L$-graded polynomial ring with $\deg x_i=l_i$ with $L$-graded canonical module $\om_A=A(-l)$, where $l=\sum_{i=1}^nl_i$. 
\item The derived category $\Db(\qmod^L\!A)$ has a Serre functor $-\otimes_A\om_A[d]$ for some $d\geq0$. 
\end{enumerate}
\end{Assum}
Since $A$ is a polynomial ring, we have an equivalence $\per^L\!A\simeq\Db(\mod^L\!A)$, which shows that in the quotient, $\rD(\QMod^L\!A)$ is compactly generated by $\Db(\qmod^L\!A)$.

We next consider a vector bundle over a stack. 
Let $B=A[y_1,\dots,y_r]$ be the $L$-graded polynomial algebra over $A$ with $\deg y_i=b_i$, and consider the categories
\[ \coh\Y:=\qmod^L\!B:=\mod^L\!B/\mod^L_I\!B \quad \text{ and } \quad \Qcoh\Y:=\QMod^L\!B:=\Mod^L\!B/\Mod^L_I\!B. \]
This should be regarded as the category of (quasi-)coherent sheaves over the total space $\Y$ of the vector bundle $\bigoplus_{i=1}^r\O_\X(b_i)$ over the stack $\X$.
We have an adjoint pair $(\pi^\ast,\pi_\ast)$ defined by
\[ \pi^\ast=-\otimes_AB\colon\QMod^L\!A\to \QMod^L\!B, \qquad \pi_\ast=\res\colon\QMod^L\!B\to\QMod^L\!A. \]
For every $M\in\QMod^L\!A$ and $N\in\QMod^L\!B$, we have the projection formula $M\otimes_A\pi_\ast N\xsimeq \pi_\ast(\pi^\ast M\otimes_BN)$.
It is easy to see, using adjunctions, that the functor $\pi^\ast\colon\rD(\Qcoh\X)\to\rD(\Qcoh\Y)$ preserves compact generators.

\medskip

The main result of this section is the following variant of Theorem \ref{CY} which we give for one-sided modules.
\begin{Thm}\label{CY2}
Let $A$ and $I$ be as in Assumption \ref{stack}, let $E=\bigoplus_{j=1}^rA(b_j)$ be such that $\Wedge^rE\cong\om_A$ in $\qmod^L\!A$, and put $B=A[y_1,\ldots,y_r]$ with $\deg y_i=b_i$. Let $T\in\rD(\Qcoh\X)$ be a compact generator and put
\[ \La=\REnd_\X(T), \quad \Ga=\REnd_\Y(\pi^\ast T)=\RHom_\X(T,T\lotimes_\X\O_\Y), \]
and regard $\Ga$ as an Adams graded dg algebra.
\begin{enumerate}
\item There is a functorial isomorphism $\RHom_\Ga(M,\Ga)[d+r]\cong DM(r)$ for $M\in\rD^\Z_\La(\Ga)$, where $\rD^\Z_\La(\Ga)$ is the full subcategory of graded derived category of $\Ga$ consisting of objects whose cohomology is perfect over $\La$.  
\item If $\Ga$ is concentrated in (cohomological) degree $0$ and is a module-finite over a commutative Gorenstein ring $R$ of dimension $d+r$, then $\Ga$ is a locally symmetric $R$-order and is a bimodule $(d+r)$-CY algebra of Gorenstein parameter $r$.
\end{enumerate}
\end{Thm}

We start with the following compatibility of dualities. For a compactly generated triangulated category $\T$, we denote by $\T^c$ the subcategory of compact objects.
\begin{Lem}\label{com}
There is a commutative diagram
\[ \xymatrix@C=25mm{
	\Db(\coh\Y)\ar[r]^-{\RHom_\Y(\pi^\ast T,-)}\ar[d]_-{\RcHom_\Y(-,\O_\Y)}&\per\Ga\ar[d]^-{\RHom_\Ga(-,\Ga)}\\
	\Db(\coh\Y)\ar[r]^-{\RHom_\Y(\pi^\ast T^\vee,-)}&\per\Ga^\op. } \]
\end{Lem}
\begin{proof}
   Since $\pi^\ast T\in\Db(\coh\Y)$ is a generator and $\Ga=\REnd_\Y(\pi^\ast T)$, we have the following isomorphisms: 
   \begin{align*}
   &\RHom_\Ga(\RHom_\Y(\pi^\ast T,-),\Ga)=\RHom_\Ga(\RHom_\Y(\pi^\ast T,-),\RHom_\Y(\pi^\ast T,\pi^\ast T))\\
   &\ysimeq\RHom_\Y(-,\pi^\ast T)\xsimeq\RHom_\Y(\pi^\ast T^\vee,(-)^\vee) =\RHom_\Y(\pi^\ast T^\vee,\RcHom_\Y(-,\O_\Y)).
   \end{align*}
   Therefore the assertion follows.
\end{proof}

We denote by $\om_\X\in\coh\X$ the object corresponding to $\om_A\in\qmod^L\!A$. By Assumption \ref{stack}, the functor $-\otimes_\X\om_\X[d]$ gives a Serre functor on $\Db(\coh\X)$.
We say that $F\in\coh\Y$ is {\it coherent over $\X$} if $\pi_\ast F\in\coh\X$. Similarly, an object in $\Db(\coh\Y)$ is coherent over $\X$ if each cohomology is coherent over $\X$.
\begin{Lem}\label{FK}
Let $F\in\Db(\coh\Y)$ be coherent over $\X$. There is an isomorphism $\RcHom_\X(\pi_\ast F,\O_\X)\cong\pi_\ast\RcHom_\Y(F,\pi^\ast\om_\X^{-1})[r]$ in $\rD(\X)$.
\end{Lem}
\begin{proof}
	Consider the Koszul complex
	\[ \xymatrix{ 0\ar[r]&\O_\Y\otimes_\X\Wedge^rE^\vee\otimes_\X\O_\Y\ar[r]&\cdots\ar[r]&\O_\Y\otimes_\X E^\vee\otimes_\X\O_\Y\ar[r]&\O_\Y\otimes_\X\O_\Y\ar[r]&0 } \]
	which is a resolution of $\O_\Y$. Applying the functor $F\otimes_\Y-$ followed by $(-)^\dagger:=\RcHom_\Y(-,\om_\X^{-1}\otimes_\X\O_\Y)[r]$, we get a complex in the first row below which computes $\RcHom_\Y(F,\om_\X^{-1}\otimes_\X\O_\Y)[r]$.
	\[\small \xymatrix{
		(\pi_\ast F\otimes_\X\O_\Y)^\dagger\ar[r]\ar[d]^-\rsimeq&\cdots\ar[r]&(\pi_\ast F\otimes_\X\Wedge^{r-1}E^\vee\otimes_\X\O_\Y)^\dagger\ar[r]\ar[d]^-\rsimeq&(\pi_\ast F\otimes_\X\Wedge^rE^\vee\otimes_\X\O_\Y)^\dagger\ar[d]^-\rsimeq \\
		\RcHom_\X(\pi_\ast F,\O_\X)\otimes_\X\Wedge^rE^\vee\otimes_\X\O_\Y\ar[r]&\cdots\ar[r]&\RcHom_\X(\pi_\ast F,\O_\X)\otimes_\X E^\vee\otimes_\X\O_\Y\ar[r]&\RcHom_\X(\pi_\ast F,\O_\X)\otimes_\X\O_\Y. } \]
	Now, each term in the first row is isomorphic to 
	\begin{align*}
	(\pi_\ast F\otimes_\X\Wedge^{r-l}E^\vee\otimes_\X\O_\Y)^\dagger &=\RcHom_\X(\pi_\ast F\otimes_\X\Wedge^{r-l}E^\vee,\om_\X^{-1}\otimes_\X\Y)\\ 
	&=\RcHom_\X(\pi_\ast F,\Wedge^lE^\vee)\otimes_\X\O_\Y=\RcHom_\X(\pi_\ast F,\O_\X)\otimes_\X\Wedge^lE^\vee\otimes_\X\O_\Y,
	\end{align*}
	which is therefore obtained by applying $\RcHom_\X(\pi_\ast F,\O_\X)\otimes_\X-$ to the (one-sided) Koszul complex, as in the second row. It is easy to see that each square commutes, and hence we deduce that $\RcHom_\X(\pi_\ast F,\O_\X)\cong\pi_\ast\RcHom_\Y(F,\om_\X^{-1}\otimes_\X\O_\Y)[r]$, as desired.
\end{proof}

We denote by $\mathrm{D^b_\X}(\coh\Y)$ the subcategory of $\Db(\coh\Y)$ consisting of complexes each of whose cohomology is coherent over $\X$. 
Also, we denote by $\mathrm{D^b_\La}(\Ga)$ the subcategory of $\rD(\Ga)$ consisting of dg $\Ga$-modules which are perfect over $\La$. 
\begin{Lem}\label{subs}
The equivalence $\RHom_\Y(\pi^\ast T,-)\colon\rD(\Qcoh\Y)\to\rD(\Ga)$ restricts to an equivalence \\ $\mathrm{D^b_\X}(\coh\Y)\xsimeq\mathrm{D^b_\La}(\Ga)$.
\end{Lem}
\begin{proof}
	We have to show that $F\in\coh\Y$ is coherent over $\X$ if and only if $\RHom_\Y(\pi^\ast T,F)$ is perfect over $\La$.
	Since we have an equivalence $\RHom_\X(T,-)\colon\rD(\Qcoh\X)\xsimeq\rD(\La)$ which restricts to $\Db(\coh\X)\xsimeq\per\La$, an object $G\in\Qcoh\X$ is coherent if and only if $\RHom_\X(T,G)$ is perfect over $\La$.
	It follows that for $F\in\rD(\Qcoh\Y)$, the dg $\La$-module $\RHom_\Y(\pi^\ast T,F)=\RHom_\X(T,\pi_\ast F)$ is perfect if and only if $F$ is coherent over $\X$, which proves the assertion.
\end{proof}

We are now ready to prove the main result of this section.
\begin{proof}[Proof of Theorem \ref{CY2}]
	(1)  Let $F\in\mathrm{D^b_\X}(\coh\Y)$. We apply the functor $\RHom_\X(T^\vee,-)$ to the isomorphism in Lemma \ref{FK}, which gives an isomorphism
	\[ \xymatrix{ \RHom_\X(T^\vee,\RcHom_\X(\pi_\ast F,\O_\X))\ar@{-}[r]^-\cong&\RHom_\X(T^\vee,\pi_\ast\RcHom_\Y(F\otimes_\Y\pi^\ast\om_\X,\O_\Y))[r] }. \]
	Now, the left-hand-side is isomorphic to $\RHom_\X(\pi_\ast F,T)\cong D\RHom_\X(T,\pi_\ast F\otimes_\X\om_\X[d])$ by Serre duality, and thus to $D\RHom_\Y(\pi^\ast T,F\otimes_\Y\pi^\ast\om_\X)[-d]$ by adjunction and the projection formula. 
	On the other hand, the right-hand-side is $\RHom_\Y(\pi^\ast T^\vee,\RcHom_\Y(F\otimes_\Y\pi^\ast\om_\X,\O_\Y))[r]$ by adjunction, and thus is isomorphic to $\RHom_\Ga(\RHom_\Y(\pi^\ast T,F\otimes_\Y\pi^\ast\om_\X),\Ga)[r]$ by Lemma \ref{com}. If $F$ varies in $\mathrm{D^b_\X}(\coh\Y)$, then the $\Ga$-module $\RHom_\Y(\pi^\ast T,F\otimes_\Y\pi^\ast\om_\X)$ varies along $\mathrm{D^b_\La}(\Ga)$ by Lemma \ref{subs}, which shows that there is a functorial isomorphism $DM[-d]\cong\RHom_\Ga(M,\Ga)[r]$ for $M\in\mathrm{D^b_\La}(\Ga)$.
	
	If one incorporates the grading, the Koszul complex in the proof of Lemma \ref{FK} becomes the one as in the proof of Theorem \ref{CY}, and the isomorphism becomes $\RcHom_\X(\pi_\ast F,\O_\X)(r)\cong\pi_\ast\RcHom_\Y(F,\pi^\ast\om_\X^{-1})[r]$. It follows from the above argument that we have an isomorphism $DM(r)[-d]\cong\RHom_\Ga(M,\Ga)[r]$ in $\rD^\Z(\Ga^\op)$.
	
	(2)  This follows from (1) and \cite[Theorem 3.3(2)]{IR_FZ}.
\end{proof}

As in Corollary \ref{CY-tilting}, when $\pi^\ast T$ is tilting, we obtain a CY algebra concentrated in degree $0$, and therefore an NCCR of the section ring. 
\begin{Cor}\label{CY2-tilting}
Under the setting of Theorem \ref{CY2}(2), suppose that $\pi^\ast T\in\Db(\coh \Y)$ is tilting.
\begin{enumerate}
\item The algebra $\Ga=\End_\Y(\pi^\ast T)=\bigoplus_{i\geq0}\Hom_\X(T,T\otimes_\X\Sym^iE^\vee)$ 
is $(d+r)$-CY of Gorenstein parameter $r$.
\item If $T$ contains a line bundle, and the corresponding idempotent $e$ gives a height $0$ progenerator $e\Ga$, then $\Ga$ is an NCCR of $\End_\Y(\O_\Y)$.
\end{enumerate}
\end{Cor}
\begin{proof}
(1)  This is immediate from Theorem \ref{CY2}.

(2)  This follows from (1) and Proposition \ref{eBe}.
\end{proof}

\section{NCCRs arising from toric Deligne--Mumford stacks}
\label{sec_toric_stack}

In this section, we apply Theorem~\ref{CY2} to the case where $\calX$ is a toric Deligne--Mumford stack, 
and produce some toric NCCRs. 

\subsection{Preliminaries on toric Deligne--Mumford stacks}
\label{subsec_stacksetting}

In this subsection, we review toric Deligne--Mumford stacks and their fundamental properties following \cite{BCS_toricstack}, \cite{FMN_toricstack}. 

Let $N$ be a finitely generated abelian group of rank $d$, and $N_{\rm tor}$ be the subgroup of $N$ consisting of torsion elements. 
For each $b\in N$, we denote the image of $b$ under the natural projection $N\rightarrow N/N_{\rm tor}$ by $\overline{b}$. 
Then, we consider the $d$-dimensional vector space $N_\RR\coloneqq N\otimes_\ZZ\RR$. 
Let $\Sigma$ be a complete rational simplicial fan in $N_\RR$, in which any cone in $\Sigma$ is generated by linearly independent vectors. 
Let $\Sigma(1)=\{\rho_1,\cdots,\rho_n\}$ be the set of rays (one dimensional cones) in $\Sigma$. 
We assume that the rays $\rho_1,\cdots,\rho_n$ span $N_\RR$. 
Let $\beta\colon \ZZ^n\rightarrow N$ be a group morphism such that for each $i=1, \dots, n$ 
the element $\overline{\beta(\e_i)}$ lies on the ray $\rho_i$, 
where $\e_1, \dots, \e_n$ form the canonical basis of $\ZZ^n$. We note that the cokernel of $\beta$ is finite by our assumption. 
The triple $\bfSigma\coloneqq(N,\Sigma, \beta)$ is called a \emph{stacky fan}. 

Using a stacky fan $\bfSigma$, we define a toric Deligne--Mumford stack. The construction is similar to Cox's construction of toric varieties (see \cite{Cox}). 
First, since $N$ is a finitely generated abelian group of rank $d$, we have a projective resolution of the $\ZZ$-module $N$ given as 
$0\rightarrow \ZZ^l \xrightarrow{Q}  \ZZ^{d+l}\rightarrow N\rightarrow 0$. 
We consider the map $B\colon \ZZ^n\rightarrow \ZZ^{d+l}$ which lifts $\beta\colon \ZZ^n\rightarrow N$, 
which fits into a commutative diagram with exact rows.
\[ \xymatrix{
	0\ar[r]&\Z^{n-d}\ar[r]\ar[d]&\Z^n\ar[r]^-\be\ar[d]_-B&N\ar@{=}[d]\\
	0\ar[r]&\Z^l\ar[r]^-Q&\Z^{d+l}\ar[r]&N\ar[r]&0 } \]
Then, we define the \emph{dual group} $\DG(\beta)\coloneqq \Coker[B\,|\,Q]^*= (\ZZ^{n+l})^*/\image[B\,|\,Q]^*$, where $(-)^*\coloneqq \Hom_\ZZ(- , \ZZ)$, and let $\beta^\vee\colon (\ZZ^n)^*\rightarrow \DG(\beta)$ be the composition of the inclusion map 
$(\ZZ^n)^*\rightarrow(\ZZ^{n+l})^*$ and the quotient map $(\ZZ^{n+l})^*\rightarrow \DG(\beta)$. 
Thus we have a commutative diagram with exact rows, in which the right square is a push-out.
\[ \xymatrix{
	0\ar[r]&N^*\ar[r]\ar@{=}[d]&\Z^{d+l}\ar[d]_-{B^*}\ar[r]^-{Q^*}&\Z^l\ar[d]\\
	0\ar[r]&N^*\ar[r]&\Z^n\ar[r]^-{\be^\vee}&\DG(\be) } \]
Note that if $X$ is an ordinary toric variety defined by the fan $\Si$ in $N_\RR$, the second sequence is precisely the one $0\to M\to \Z^n\to\Cl(X)\to 0$ computing the divisor class group.

Let $k=\CC$. We consider the ideal $J_\Sigma$ generated by monomials associated with each cone $\sigma$ in $\Sigma$: 
\[ J_\Sigma\coloneqq\Big\langle\prod_{\rho_i\not\subset\sigma}x_i\in k[x_1,\cdots,x_n] \mid \sigma\in\Sigma\,\Big\rangle, \] 
and let $U_\Sigma\coloneqq\AA^n{\setminus}V(J_\Sigma)$, where $V(J_\Sigma)$ is the vanishing locus of $J_\Sigma$. 
Applying $\Hom_\ZZ(-,k^\times)$ to $\beta^\vee\colon (\ZZ^n)^*\rightarrow \DG(\beta)$, we have a morphism $G\coloneqq\Hom_\ZZ(\DG(\beta),k^\times)\rightarrow(k^\times)^n$. 
Since we can consider the natural action of $(k^\times)^n$ on $\AA^n=\Spec k[x_1,\cdots,x_n]$, it induces the action of $G$ on $\AA^n$ and $U_\Sigma$. 
Then we define the quotient stack $\calX_\bfSigma\coloneqq[U_\Sigma/G]$, which is a smooth separated Deligne--Mumford stack
whose coarse moduli space is the toric variety $X_\Sigma$ defined by the simplicial fan $\Sigma$ (see \cite[Section~3]{BCS_toricstack}). 
Thus, we call $\calX_\bfSigma$ the \emph{toric Deligne--Mumford (DM) stack} associated to a stacky fan $\bfSigma=(N,\Sigma, \beta)$. 
If $N$ is torsion free, then $\calX_\bfSigma$ is also called the \emph{toric DM orbifold} associated to $\bfSigma$. 

Next, we discuss line bundles on $\calX_\bfSigma$ (cf.\ \cite[Section~3]{BH_Kingconj}, \cite[Section~3]{JT_toricstack}). 
Since there exists an equivalence $\coh\calX_\bfSigma\simeq \coh^GU_\Sigma$, where $\coh^GU_\Sigma$ is the category of $G$-equivariant coherent sheaves of $U_\Sigma$ (see \cite[(7.21)]{Vis_moduli}), 
we consider $G$-equivariant line bundles on $U_\Sigma$. 
For a subset $I\subset \{1, \dots, n\}$, we let $P_I\coloneqq \langle x_i \mid i\in I\rangle$. 
The ideal $J_\Sigma$ can be described as $J_\Sigma=\bigcap_I P_I$ where the intersection runs over the subsets $I\subset\{1, \dots ,n\}$ such that $\overline{\beta(\e_i)}\in N_\RR$ $(i\in I)$ do not generate any cone in $\Sigma$. 
Since the element $\overline{\beta(\e_i)}$ generates the ray $\rho_i\in\Sigma$ for each $i$, we have $\mathrm{ht}(J_\Sigma)\ge 2$, 
and hence $\mathrm{codim}_{\AA^n}V(J_\Sigma)\ge 2$, which implies $\Pic(U_\Sigma)\cong \Pic(\AA^n)=0$. 
Thus, any $G$-equivariant line bundle on $U_\Sigma$ can be described as the trivial line bundle $U_\Sigma\times k$ 
equipped with a $G$-action determined by a character $\chi\in G^\vee$, given by 
$g\cdot (u, t)=(gu, \chi(g)t)$, where $g\in G$ and $(u,t)\in U_\Sigma\times k$. 
From these arguments, we can identify a line bundle on $\calX_\bfSigma$ with a $G$-equivariant trivial line bundle on $U_\Sigma$,
and hence we obtain an isomorphism $\Pic(\calX_\bfSigma)\cong G^\vee\cong\DG(\beta)$. 
Let $D_i$ be the toric divisor on $\calX_\bfSigma$ corresponding to the ray $\rho_i$ in $\Sigma(1)$ $(i=1, \dots, n)$, see e.g., \cite[Subsection~2.2]{OU_frobenius}.

\medskip

These geometric inputs give rise to our algebraic setup from Section~\ref{quotstack}.
Given a stacky fan $\bSi=(N,\Si,\be)$ as above, let $L=\mathrm{DG}(\be)$ be the grading group, and let
\[ A=k[x_1,\ldots,x_n] \text{ with } \deg x_i:=a_i:=\be^\vee(\e_i)\in L, \text{ and } I=J_\Si. \]
\begin{Prop}\label{AB}
The $L$-graded algebra $A$ and the ideal $I$ satisfy Assumption \ref{stack}.
\end{Prop}
\begin{proof}
Since the graded canonical module of $A=k[x_1,\ldots,x_n]$ is $\om_A=A(-c)$ with $c=\sum_{i=1}^na_i$, this shows (a).

Now consider the category $\qmod^L\!A$ and its derived category. 
The toric DM stack $\X=\X_{\bSi}$ associated with the stacky fan $\bSi$ satisfies $\coh\X=\coh^G\!U$ for $U=\Spec A\setminus V(I)$ and $G=\Spec k[L]$; thus we have an equivalence $\coh\X=\qmod^L\!A$ for $\qmod^L\!A=\mod^L\!A/\mod_{I}^L\!A$ from Section \ref{quotstack}.
Now, the canonical bundle $\om_\X=\O_\X(-\sum_{i=1}^n{D}_i)$ on $\X$ corresponds under the equivalence $\coh\X\simeq\qmod^L\!A$ to $\bigoplus_{i=1}^nA(-a_i)=\om_A$, and since $-\otimes_\X\om_\X[d]$ for $d=\dim\X$ gives a Serre functor on $\Db(\coh\X)$ by \cite[Theorem 2.20]{Nironi}, we deduce that $-\otimes_A\om_A[d]$ is the Serre functor of $\Db(\qmod^L\!A)$, proving (b).
\end{proof}

In what follows, we consider vector bundles given by $E=\bigoplus_{i=1}^rA(b_i)$ whose total space is CY, which is characterized as follows.
\begin{Lem}\label{CC}
Let $E=\bigoplus_{i=1}^rA(b_i)$. Then it satisfies $\Wedge^rE\cong\om_A$ in $\qmod^L\!A$ if $\sum_{i=1}^na_i+\sum_{i=1}^rb_i=0$.
\end{Lem}
\begin{proof}
This is because $\om_A=A(-\sum_{i=1}^na_i)$ and $\Wedge^rE=A(\sum_{i=1}^rb_i)$.
\end{proof}

\subsection{Constructions of NCCRs from weighted projective stacks}
\label{subsec_WPS}

We keep the notation in Subsection~\ref{subsec_stacksetting}. 
In this subsection, we construct some toric NCCRs using tilting bundles of toric DM stacks and Corollary~\ref{CY2-tilting}. 

We first note that a toric NCCR of a toric ring whose class group is $\ZZ$, which was already given in \cite[Theorem~8.9]{VdB_NCCR}, 
can be constructed by our framework. 

\begin{Setup}
\label{setting_toric_Z}
Let $k=\CC$. 
Let $n\geq2$ and $r\geq1$ be integers, and consider the polynomial ring
\[ S=k[x_1,\ldots,x_n, y_1,\ldots,y_r]. \] 
Let the algebraic torus $G\coloneqq k^\times$ act on $S$ by weights 
$q_i\in\ZZ_{>0}$ on $x_i$ $(i=1, \dots, n)$ and $a_j\in\ZZ_{<0}$ on $y_j$ $(j=1, \dots, r)$ 
such that $\sum_{i=1}^nq_i+\sum_{j=1}^ra_j=0$. Then, let $R\coloneqq S^G$ 
be the invariant subring, which is a Gorenstein toric ring of dimension $n+r-1$. 
Here, we assume that $\gcd(q_1, \dots, q_n, a_1, \dots, a_r)=1$ without loss of generality. 
\end{Setup}

To construct a toric NCCR of the toric ring $R$ as in Setting~\ref{setting_toric_Z}, we use a toric DM stack whose Picard group is $\ZZ$, 
which is a weighted projective stack (see \cite[Subsection~7.3]{FMN_toricstack}). 
First, let $m\coloneqq \gcd(q_1, \dots, q_n)$, let $q_i=m q_i^\prime$ for $i=1, \dots, n$, and define the stacky fan $\bfSigma=(N,\Sigma, \beta)$ as follows (cf. \cite[Subsection~2.1]{wang_toricHir}). 
Let
\[ N=\ZZ^{n-1}\oplus\ZZ/m\ZZ \]
and $\Sigma\subset N_\RR$ be the fan of the weighted projective space $\PP(q_1^\prime, \dots, q_n^\prime)$  (see e.g., \cite[Proposition~2.5]{wang_toricHir}). 
Let $w_1, \dots, w_n\in N/N_{\rm tor}$ be the primitive generators of rays in $\Sigma$. 
Then, we let $b_i=\frac{m_i}{m}w_i$ $(i=1, \dots, n)$, where $m_i\coloneqq \gcd(q_1, \dots, q_{i-1}, q_{i+1}, \dots, q_n)$. 
Note that $b_1, \dots, b_n$ are also generators of rays in $\Sigma$. 
We define a group morphism $\beta\colon \ZZ^n\rightarrow N$ as $\beta(\e_i)={\small \begin{pmatrix}b_i\\ c_i\end{pmatrix}}$ for $i=1, \dots, n$, 
where $c_1, \dots, c_n$ are integers with $0\le c_1, \dots, c_n\le m-1$ and \\ $c_1q_1^\prime+\cdots+c_nq_n^\prime\equiv 1 \pmod m$. 
The lifting of $\beta$ is given by the matrix $B=(\beta(\e_1), \dots, \beta(\e_n))={\small \begin{pmatrix}b_1&\cdots&b_n\\c_1&\cdots&c_n\end{pmatrix}}$. 
We then consider the matrix $[B\,|\,Q]$ where $Q={\small \begin{pmatrix}{\bf 0}\\m\end{pmatrix}}$. 
We see that \[\DG(\beta)= (\ZZ^{n+1})^*/\image [B\,|\,Q]^*\cong\ZZ\] and the integer kernel of $[B\,|\,Q]$ is generated by 
a vector with the form ${}^t(q_1, \dots, q_n, u)={}^t(m q_1^\prime, \dots, m q_n^\prime, u)$, where $u$ is an integer with $c_1q_1+\cdots+c_nq_n+um=0$. 
Thus, we see that $\beta^\vee\colon\ZZ^n\rightarrow \ZZ$ is given by $(q_1, \dots, q_n)$, 
and $g\in G=k^\times$ acts on $(u_1, \dots, u_n)\in U_\Sigma=\AA^n{\setminus}\{(0, \dots, 0)\}$ by 
\[ g\cdot (u_1, \dots, u_n)=(g^{q_1}u_1, \dots, g^{q_n}u_n),\]
in which case $\calX_\bfSigma$ is the weighted projective stack $\calP(q_1, \dots, q_n)$, see \cite[Proposition~2.7]{wang_toricHir}. 
We also see that $\Pic(\calX_\bfSigma)\cong \ZZ$ and $\calX_\bfSigma$ has a tilting bundle $\bigoplus_{l=0}^{q}\O_{\calX_\bfSigma}(l)$, 
where $q\coloneqq q_1+\cdots+q_n-1$, see \cite[Proposition~5.1]{BH_Kingconj}. 

By Proposition \ref{AB}, the algebra $A=k[x_1,\ldots,x_n]$ and the ideal $I=J_\Si$ for the stacky fan $\bSi$ satisfy the conditions in Assumption \ref{stack}. Also by Lemma \ref{CC}, the vector bundle $E=\bigoplus_{j=1}^r\O_\X(a_j)$ satisfies $\Wedge^rE\cong\om_A$. This leads to the following result.

\begin{Thm}[{cf. \cite[Theorem~8.9]{VdB_NCCR}, see also \cite[Corollary~3.11]{Han_NCCR}}]
\label{NCCR_toricZ_revisited}
Let $\calX$ be the weighted projective stack $\calP(q_1, \dots, q_n)$ and $T$ be the tilting bundle $\bigoplus_{l=0}^{q}\O_\calX(l)$ of $\calX$ as above. 
For the vector bundle $E=\bigoplus_{j=1}^r\O_\calX(a_j)$, we consider the total space $\Y$ of $E$ and the map $\pi\colon \Y\rightarrow \calX$, 
where $a_j$ is the weight defined in Setting~\ref{setting_toric_Z}. 
\begin{enumerate}
\item $\pi^\ast T$ is a tilting bundle on $\Y$ and $\Ga:=\End_\Y(\pi^\ast T)=\bigoplus_{i\geq0}\Hom_\X(T,T\otimes_\X\Sym^iE^\vee)$ is \\ an $(n-1+r)$-CY algebra of Gorenstein parameter $r$.
\item $\Ga$ is an NCCR of the toric ring $R$ given in Setting~\ref{setting_toric_Z}.
\end{enumerate}
\end{Thm}

\begin{proof}
By the same argument as the proof of Theorem~\ref{NCCR_toric_PxP}, it suffices to prove that 
$\Ext^i_\calX(T, T\otimes\Sym E^\vee)=0$ for any $i>0$ to show that $\pi^\ast T$ is a tilting bundle on $\Y$. 
Since the canonical bundle of $\calX$ is $\O_\calX(-q-1)$, we see that $\O_\calX(l)$ is acyclic if $l>-q-1$ by \cite[Propositions~4.3 and 4.5]{BH_Kingconj}.
Since $0\le l\le q$ and $a_j\in\ZZ_{<0}$, we have $\Ext^i_\calX(T, T\otimes\Sym E^\vee)=0$. 
Thus, by Corollary~\ref{CY2-tilting}, $\End_\Y(\pi^\ast T)$ is $(n-1+r)$-CY of Gorenstein parameter $r$, and is an NCCR of the section ring $\bigoplus_{i\geq0}\Hom_\calX(\O_\calX,\O_\calX\otimes\Sym^iE^\vee)$, which is isomorphic to $R$. 
\end{proof}

The following example shows that, for a tilting bundle $T$ on the base space $\X$, 
the quiver of $\End_\X(T)$ may be disconnected, whereas the corresponding quiver for the total space $\Y$ is connected. 

\begin{Ex}
Let $n=2$, $r=2$. We apply Theorem~\ref{NCCR_toricZ_revisited} to the weights $q_1=q_2=2$ and $a_1=-1, a_2=-3$, 
thus $\X=\P(2,2)$ and $E=\O_\X(-1)\oplus\O_\X(-3)$. We note that $\gcd(q_1,q_2)=2$, while $\gcd(q_1,q_2,a_1,a_2)=1$. 
For the tilting bundle $T=\bigoplus_{l=0}^3\O_\X(l)$, the quiver of $\End_\X(T)$ is as depicted on the left-hand side of the figure below. 
In particular, it is disconnected. 
Meanwhile, the quiver of $\End_\Y(\pi^*T)$, which is a $3$-CY algebra of Gorenstein parameter $2$ and an NCCR of the section ring of $\Y$, 
is as depicted on the right-hand side of the figure, 
where the triple arrow pointing northwest represents the multiplication maps by $x_1y_2$, $x_2y_2$, and $y_1$, respectively.
In particular, it is connected. 
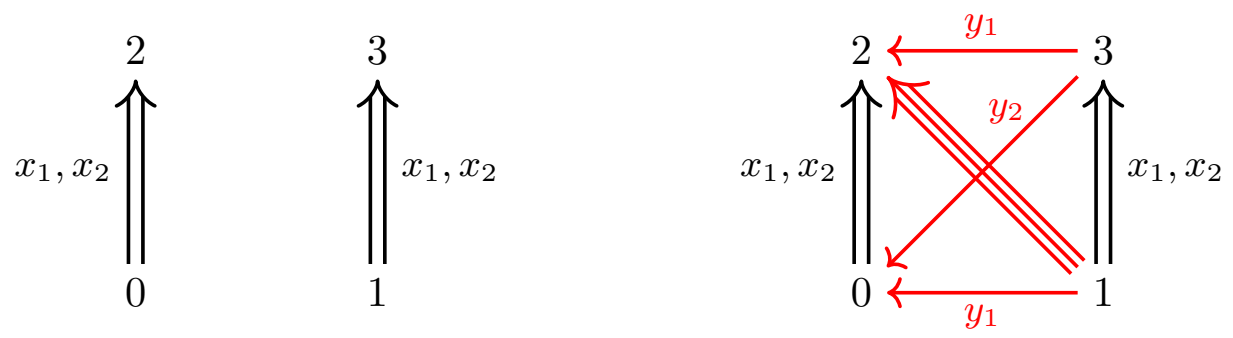
\begin{center}
{\scalebox{0.8}{
\begin{tikzpicture}
\newcommand{\arrowwidth}{0.03cm} 

\node at (0,0) {
\begin{tikzpicture}
\foreach \a/\b/\n in {0/0/0,1/0/1,0/1/2,1/1/3} {\node (V\a\b) at (\a*2,\b*2)  {$\n$};}; 

\foreach \s/\t in {00/01,10/11} {\draw[double distance=3pt, scaling nfold=2, -Implies, line width=\arrowwidth] (V\s)--(V\t);}; 
\node at (-0.6,1) {\small $x_1, x_2$}; \node at (2.6,1) {\small $x_1, x_2$}; 
\end{tikzpicture}}; 

\node at (6,0) {
\begin{tikzpicture}
\foreach \a/\b/\n in {0/0/0,1/0/1,0/1/2,1/1/3} {\node (V\a\b) at (\a*2,\b*2)  {$\n$};}; 

\foreach \s/\t in {00/01,10/11} {\draw[double distance=3pt, scaling nfold=2, -Implies, line width=\arrowwidth] (V\s)--(V\t);};  
\node at (-0.6,1) {\small $x_1, x_2$}; \node at (2.6,1) {\small $x_1, x_2$}; 

\foreach \s/\t in {10/00,11/01} {\draw[->,line width=\arrowwidth, red] (V\s)--(V\t);}; 
\draw[->, line width=\arrowwidth, red] (V11)--(V00); 
\draw[double distance=1.8pt, scaling nfold=3, -Implies, line width=\arrowwidth, red]  (V10)--(V01);  
\node[red] at (1,-0.2) {\small $y_1$}; \node[red] at (1,2.2) {\small $y_1$}; \node[red] at (1.2,1.5) {\small $y_2$}; 
\end{tikzpicture}}; 

\end{tikzpicture}
} }
\end{center}
\end{Ex}

\subsection{Some examples of NCCRs from toric Fano DM stacks of Picard rank two}

We keep the notation in Subsection~\ref{subsec_stacksetting}. 
In this section, we show that we can also construct toric NCCRs using a tilting bundle of a toric Fano DM stack whose Picard group is $\ZZ^2$ and Corollary~\ref{CY2-tilting}. 
Here, we say that a toric DM stack $\calX_\bfSigma$ is \emph{Fano} if $\overline{\beta(\e_1)}, \dots, \overline{\beta(\e_n)}$ are precisely the vertices of a simplicial convex polytope in $N_\RR$. 

The existence of a tilting bundle (or full strong exceptional collection) of a toric Fano DM stack $\calX_\bfSigma$ 
with Picard rank two was shown in \cite[Theorem~5.11]{BH_Kingconj}. 
Although we will use a tilting bundle different from the ones in \cite{BH_Kingconj}, we heavily use ideas introduced in \cite[Section~5]{BH_Kingconj}. 
Thus, we here review some materials used in \cite{BH_Kingconj}. 
Tilting bundles of $\calX_\bfSigma$ given in \cite{BH_Kingconj} can be obtained from (shifts of) a ``window $P$" in $\Pic(\calX_\bfSigma)_\RR\cong\RR^2$, 
that is, line bundles contained in $P$ as elements of $\Pic(\calX_\bfSigma)$ form a full strong exceptional collection of $\calX_\bfSigma$ 
and the direct sum of such bundles is a tilting bundle. 
The $\Ext$-vanishing condition for line bundles contained in $P$ can be checked by observing \emph{forbidden cones} in $\Pic(\calX_\bfSigma)_\RR$. 
To define the forbidden cones, we first consider non-zero rational numbers $\alpha_i$ such that 
\[
\sum_{i=1}^n\alpha_i\overline{\beta(\e_i)}={\bf 0} \text{ and } \sum_{i=1}^n\alpha_i=0, 
\]
which are determined uniquely up to scalar multiplication, see \cite[Proposition~5.4]{BH_Kingconj}. 
Let \[I_+=\{i\in\{1, \dots, n\} \mid \alpha_i>0\}, \quad I_-=\{i\in\{1, \dots, n\} \mid \alpha_i<0\}.\] 
Then, there are precisely three forbidden cones $F_{\varnothing}, F_{I_+}, F_{I_-}$ in $\Pic(\calX_\bfSigma)_\RR$ 
taking the following form (see \cite[Proposition~5.7]{BH_Kingconj}): 
\begin{align*}
	&F_{\varnothing}=-\sum_{i=1}^nD_i-\sum_{i=1}^n\RR_{\ge0}D_i, \\
	&F_{I_+}=-\sum_{i\in I_-}D_i+\sum_{i\in I_+}\RR_{\ge0}D_i-\sum_{i\in I_-}\RR_{\ge0}D_i, \\
	&F_{I_-}=-\sum_{i\in I_+}D_i+\sum_{i\in I_-}\RR_{\ge0}D_i-\sum_{i\in I_+}\RR_{\ge0}D_i, 
\end{align*}
and any line bundle whose image in $\Pic(\calX_\bfSigma)_\RR$ does not lie in any of these forbidden cones is acyclic, see \cite[Proposition~4.5]{BH_Kingconj}. 
Thus, identifying the forbidden cones would be beneficial to detect a strong exceptional collection. 
Also, by \cite[Corollary~4.8]{BH_Ktheory}, it is enough to show that all line bundles on $\calX_\bfSigma$ can be generated by line bundles contained in $P$, 
to see that such line bundles form a generator of $\Db(\coh \calX_\bfSigma)$.

\medskip

We will give examples of NCCRs for some families of toric rings using a tilting bundle of a toric Fano DM stack of Picard rank two.

\subsubsection{Example~1}\label{Subsec_ex_toricDM1}

\begin{Setup}
	\label{setting_toric_Z^2-1}
	Let $k=\CC$. Let $n, m\geq2$ and $r\geq1$ be integers, and consider the polynomial ring 
	\[ S=k[x_1,\ldots, x_n, y_1,\ldots,y_m, z_1,\ldots, z_r]. \] 
	Let $p_1, \dots, p_n, q_1, \dots, q_m$ be positive integers, 
	and suppose that the algebraic torus $G\coloneqq (k^\times)^2$ acts on $S$ by weights $(p_i,0)$ on $x_i$ $(i=1, \dots,n)$, $(0,q_j)$ on $y_j$ $(j=1,\dots,m)$, and $(c_l, d_l)\in\ZZ_{\leq 0}^2$ on $z_l$ $(l=1,\dots,r)$ respectively,  
	where $c_l, d_l$ are non-positive integers such that 
	\[\sum_{i=1}^n(p_i,0)+\sum_{j=1}^m(0,q_j)+\sum_{l=1}^r(c_l,d_l)=(0,0)\] and $(c_l,d_l)\neq(0,0)$. 
	Then, let $R\coloneqq S^G$ be the invariant subring, which is a Gorenstein toric ring of dimension $n+m+r-2$. 
    We assume that the weights $(p_i,0), (0,q_j), (c_l,d_l)$ span the lattice $\Z^2$, otherwise the category $\qmod^{\Z^2}\!S$ splits into several connected components. Since these components are mutually equivalent, this assumption does not lose generality.
\end{Setup}

To construct NCCRs of the toric ring defined in Setting~\ref{setting_toric_Z^2-1}, we consider the following toric Fano DM stack. 
Let $p_1, \dots, p_n, q_1, \dots, q_m$ be positive integers given in Setting~\ref{setting_toric_Z^2-1}, and let $p_i=n_pp_i^\prime$ $(i=1,\dots,n)$, $q_j=n_qq_j^\prime$ $(j=1,\dots,m)$, 
where $n_p=\gcd(p_1, \dots, p_n)$, $n_q=\gcd(q_1, \dots, q_m)$. 
We take integers $s_1, \dots, s_n, t_1, \dots, t_m$ such that $0\le s_1, \dots, s_n\le n_p-1$, $0\le t_1, \dots, t_m\le n_q-1$ and 
$s_1p_1^\prime+\cdots+s_np_n^\prime\equiv 1 \pmod {n_p}$, $t_1q_1^\prime+\cdots+t_mq_m^\prime\equiv 1 \pmod {n_q}$. 
Let $\Sigma_p$ (resp. $\Sigma_q$) be the fan of the weighted projective space $\PP(p_1^\prime, \dots, p_n^\prime)$ 
(resp. $\PP(q_1^\prime, \dots, q_m^\prime)$), see e.g., \cite[Proposition~2.5]{wang_toricHir}. 
Then, let
\[ N=\ZZ^{n+m-2}\oplus \ZZ/n_p\ZZ\oplus\ZZ/n_q\ZZ \]
and let $\beta\colon \ZZ^{n+m}\to N$ be a morphism such that $\beta_i\coloneqq\beta(\e_i)$ $(i=1, \dots, n)$ and $\beta_{n+j}\coloneqq\beta(\e_{n+j})$ $(j=1,\dots,m)$ are given by 
\[\beta_i={}^t({\bf p}_i,{\bf 0},s_i,0) \hspace{7pt} (i=1,\dots,n), \quad
\beta_{n+j}={}^t({\bf 0},{\bf q}_j, 0,t_j) \hspace{7pt} (j=1,\dots,m)
\]
where ${\bf p}_1, \dots, {\bf p}_n\in\ZZ^{n-1}$ (resp. ${\bf q}_1, \dots, {\bf q}_m\in\ZZ^{m-1}$) are generators of rays in the fan $\Sigma_p$ (resp. $\Sigma_q$) defined similarly to $b_i$ in Subsection~\ref{subsec_WPS}. 
Note that $\overline{\beta}_1, \dots, \overline{\beta}_{n+m}$ generate the rays of the fan $\Sigma_p\times\Sigma_q$. 
For the stacky fan $\bfSigma_{p,q}\coloneqq(N, \Sigma_p\times\Sigma_q, \beta)$, 
the lifting of $\beta$ is given by the matrix $B={\small \begin{pmatrix}\beta_1&\cdots&\beta_{n+m}\end{pmatrix}}$. 
We consider the matrix $[B\,|\,Q]$ where $Q={\small \begin{pmatrix}{\bf 0}&{\bf 0}\\{\bf 0}&{\bf 0}\\n_p&0\\0&n_q\end{pmatrix}}$. 
Then, we see that $\DG(\beta)= (\ZZ^{n+m+2})^*/\image [B\,|\,Q]^*\cong\ZZ^2$ and the integer kernel of $[B\,|\,Q]$ is generated by 
${}^t{\small \begin{pmatrix}p_1&\cdots&p_n&0&\cdots&0&u&0\\0&\cdots&0&q_1&\cdots&q_m&0&v\end{pmatrix}}$
where $u, v$ are integers with $s_1p_1+\cdots+s_np_n+un_p=0$, $t_1q_1+\cdots+t_mq_m+vn_q=0$.
Thus, $\beta^\vee\colon \ZZ^{n+m}\to \ZZ^2$ is given by 
${\small \begin{pmatrix}p_1&\cdots&p_n&0&\cdots&0\\0&\cdots&0&q_1&\cdots&q_m\end{pmatrix}}$ 
and this determines the action of $G=\Hom_\ZZ(\DG(\beta),k^\times)\cong(k^\times)^2$ on $U_{\Sigma_p\times\Sigma_q}$. 
The toric DM stack defined by the stacky fan $\bfSigma_{p,q}$, which we denote by $\calX_{p,q}$, is Fano and $\Pic(\calX_{p,q})\cong\ZZ^2$. 
We can take a basis of $\Pic(\calX_{p,q})$ so that prime divisors $D_1, \dots, D_{n+m}$ are represented by 
\[[D_i]=(p_i,0) \hspace{7pt} (i=1,\dots,n), \quad [D_{n+j}]=(0,q_j) \hspace{7pt} (j=1,\dots,m)\]
in $\Pic(\calX_{p,q})$, and we will fix such a basis in the following argument.

\medskip

The above geometric inputs can be translated to an algebraic setup in Section \ref{quotstack} as follows.
Let $A=k[x_1,\ldots,x_n,y_1,\ldots,y_m]$ on which we give a grading by $L=\Z^2$ specified by the weights as in Setting \ref{setting_toric_Z^2-1}. Also, let $I\subset A$ be the ideal $J_{\Si_p\times\Si_q}$ for the fan $\Si_{p}\times\Si_{q}$. Then by Proposition \ref{AB}, these satisfy the conditions in Assumption \ref{stack}, and we have $\qmod^{\Z^2}\!A=\coh^{(k^\times)^2}\!(\Spec A\setminus V(I))=\coh\X_{p,q}$.

We consider a tilting bundle of $\calX_{p,q}$ given in the following form, which will be found to be compatible with the construction of NCCRs. 

\begin{Prop}
	\label{tilting_toric_Z^2-1}
	Let $\calX=\calX_{p,q}$ be the toric Fano DM stack as above. 
	For positive integers $p_1, \dots, p_n, q_1, \dots, q_m$ given in Setting~\ref{setting_toric_Z^2-1}, let $p=p_1+\dots+p_n-1$, $q=q_1+\cdots+q_m-1$ and 
	consider the set 
	\[\S=\{(a,b)\in\Pic(\calX)\cong\ZZ^2 \mid 0\le a\le p, \, 0\le b\le q\}. \]
	Then $T=\bigoplus_{(a,b)\in\S}\O_\calX(a,b)$ is a tilting bundle of $\calX$. 
\end{Prop}

\begin{proof}
	We may show that the set of line bundles $\{\O_\calX(a,b) \mid (a,b)\in\S\}$ is a full strong exceptional collection. 
	We show that $\Ext_\calX^i(L, L^\prime)=0$ for any $i>0$ and $L, L^\prime\in \{\O_\calX(a,b) \mid (a,b)\in\S\}$, 
	which is equivalent to the statement that any line bundle $\O_\calX(a,b)$ with $(a,b)\in\widetilde{\S}$ is acyclic where 
	$\widetilde{\S}=\{(a,b)\in\ZZ^2 \mid -p\le a\le p, \, -q\le b\le q\}$. 
	By \cite[Proposition~4.5]{BH_Kingconj}, it is enough to show that $(a,b)\in\widetilde{\S}$ does not lie 
	in any of the forbidden cones $F_\varnothing$, $F_{I_+}$, and $F_{I_-}$. 
	Since there is a relation 
	\[\frac{1}{p+1}\sum_{i=1}^np_i\overline{\beta}_i
	-\frac{1}{q+1}\sum_{j=1}^mq_j\overline{\beta}_{n+j}={\bf 0},\]
	we have $I_+=\{1,\dots, n\}$, $I_-=\{n+1, \dots, n+m\}$ and the forbidden cones are 
	\begin{align}\label{eq_fcone_ab}
		\nonumber&F_\varnothing=(-p-1, -q-1)-\RR_{\ge 0}(1,0)-\RR_{\ge 0}(0,1)\\ 
		&F_{I_+}=(0, -q-1)+\RR_{\ge 0}(1,0)-\RR_{\ge 0}(0,1)\\
		\nonumber&F_{I_-}=(-p-1, 0)-\RR_{\ge 0}(1,0)+\RR_{\ge 0}(0,1).
	\end{align}
	Thus, any element in $\widetilde{\S}$ does not lie in any of $F_\varnothing$, $F_{I_+}$, and $F_{I_-}$, and hence $\{\O_\calX(a,b) \mid (a,b)\in\S\}$ forms a strong exceptional collection.
	
	Then, we show that all line bundles are contained in the subcategory $\D$ of $\Db(\coh \calX)$ generated by the line bundles in $\{\O_\calX(a,b) \mid (a,b)\in\S\}$. 
	By \cite[the proof of Proposition~4.5]{BH_Ktheory}, we have an exact sequence 
	{\small\begin{equation}\label{resol_linebdl}
			0\to \O_\calX\Big((s,t)-\sum_{i\in I}[D_i]\Big)\to \bigoplus_{\substack{J\subset I \\ |J|=|I|-1}}\O_\calX\big((s,t)-\sum_{i\in J}[D_i]\big)
			\to \cdots \to \bigoplus_{\substack{J\subset I \\ |J|=1}}\O_\calX\big((s,t)-\sum_{i\in J}[D_i]\big)\to \O_\calX(s,t)\to 0
	\end{equation}}
	for any line bundle $\O_\calX(s,t)$, where $I=I_+$ or $I=I_-$. 
	We now consider the sequence \eqref{resol_linebdl} for $I=I_+=\{1, \dots, n\}$. 
	For $(s^\prime,t)\in\ZZ^2$ with $s-p-1\le s^\prime\le \max\{s-p_1, \dots, s-p_n\}$ and $0\le t\le q$, we assume that $\O_\calX(s^\prime,t)$ is contained in $\D$, 
	which means that all terms of the sequence \eqref{resol_linebdl} except the rightmost one are contained in $\D$, 
	and hence the rightmost term $\O_\calX(s,t)$ is also contained in $\D$. 
	Since $\D$ is generated by $\{\O_\calX(a,b) \mid (a,b)\in\S\}$, we see that $\O_\calX(s,t)$ is contained in $\D$ for any $s\ge 0$ and $0\le t\le q$ by the inductive argument. 
	We also have the exact sequence with the form 
	{\small\[
		0\to \O_\calX(s,t) \to \bigoplus_{\substack{J\subset I \\ |J|=1}}\O_\calX\big((s,t)+\sum_{i\in J}[D_i]\big)\to\cdots\to 
		\bigoplus_{\substack{J\subset I \\ |J|=|I|-1}}\O_\calX\big((s,t)+\sum_{i\in J}[D_i]\big)\to \O_\calX\Big((s,t)+\sum_{i\in I}[D_i]\Big)\to 0
		\]}
	from \eqref{resol_linebdl}. 
	For $(s^\prime,t)\in\ZZ^2$ with $\min\{s+p_1, \dots, s+p_n\} \le s^\prime\le s+p+1$ and $0\le t\le q$, we assume that $\O_\calX(s^\prime,t)$ is contained in $\D$,
	which means that all terms of the sequence except the leftmost one are contained in $\D$, 
	and hence the leftmost term $\O_\calX(s,t)$ is also contained in $\D$. 
	Since $\O_\calX(s,t)$ is contained in $\D$ for any $s\ge 0$ and $0\le t\le q$ by the former argument, 
	we see that $\O_\calX(s,t)$ is contained in $\D$ for any $s\in\ZZ$ and $0\le t\le q$. 
	Then, we consider the sequence \eqref{resol_linebdl} for $I_-=\{n+1, \dots, n+m\}$. 
	By a similar inductive argument, we see that $\O_\calX(s,t)$ is contained in $\D$ for any $(s,t)\in\ZZ^2\cong\Pic(\calX)$, 
	which means that $\D=\Db(\coh \calX)$ by \cite[Corollary~4.8]{BH_Ktheory}. 
\end{proof}

Now we consider the total space of a vector bundle $E=\bigoplus_{l=1}^r\O_\calX(c_l,d_l)\in\coh\X$, whose total space $\Y$ is given as $\coh\Y=\qmod^{\Z^2}\!B$ for $B=A[z_1,\ldots,z_r]=S$ with $\deg z_l=(c_l,d_l)$.
\begin{Thm}
	\label{NCCR_toric_Z^2-1}
	Let $\calX=\calX_{p,q}$ be the toric Fano DM stack as above 
	and $T=\bigoplus_{(a,b)\in\S}\O_\calX(a,b)$ be the tilting bundle of $\calX$ given in Proposition~\ref{tilting_toric_Z^2-1}. 
	For the vector bundle $E=\bigoplus_{l=1}^r\O_\calX(c_l,d_l)$, we consider the total space $\Y$ of $E$ and the map $\pi\colon \Y\rightarrow \calX$, 
	where $c_l$ and $d_l$ are the weights defined in Setting~\ref{setting_toric_Z^2-1}. 
\begin{enumerate}
\item  $\pi^\ast T$ is a tilting bundle on $\Y$ and $\Ga:=\End_\Y(\pi^\ast T)=\bigoplus_{i\geq0}\Hom_\X(T,T\otimes_\X\Sym^iE^\vee)$ is an \\ $(n+m-2+r)$-CY algebra of Gorenstein parameter $r$.
\item $\Ga$ is an NCCR of the toric ring $R$ given in Setting~\ref{setting_toric_Z^2-1}.
\end{enumerate}
\end{Thm}

\begin{proof}
	By the same argument as the proof of Theorem~\ref{NCCR_toric_PxP}, it suffices to prove that 
	$\Ext^i_\calX(T, T\otimes\Sym E^\vee)=0$ for any $i>0$ to show that $\pi^\ast T$ is a tilting bundle on $\Y$. 
	We consider $\widetilde{\S}$ used in the proof of Proposition~\ref{tilting_toric_Z^2-1}. 
	Then we see that $\Ext^i_\calX(T, T\otimes\Sym E^\vee)\cong H^i(\calX, T^\vee\otimes T\otimes \Sym E^\vee)$, 
	where $T^\vee\otimes T\otimes \Sym E^\vee$ is the direct sum of line bundles with the form: 
	\[ \O_\calX\Big(a-\sum_{l=1}^rm_lc_l, b-\sum_{l=1}^rm_l d_l\Big)\]
	for some $m_l \in\ZZ_{\ge 0}$ $(l=1, \dots, r)$ and $(a,b)\in\widetilde{\S}$. 
	Since $a\ge -p$, $b\ge -q$, and $c_l, d_l\le 0$ for any $l$, we see that 
	any lattice point $(a-\sum_{l=1}^rm_lc_l, b-\sum_{l=1}^rm_l d_l)\in\ZZ^2$ does not lie in the forbidden cones 
	$F_\varnothing$, $F_{I_+}$, $F_{I_-}$ shown in \eqref{eq_fcone_ab}, and hence $\Ext^i_\calX(T, T\otimes\Sym E^\vee)=0$ by \cite[Proposition~4.5]{BH_Kingconj}. 
    This shows that $\pi^\ast T\in\Db(\coh\Y)$ is a tilting object. The remaining assertions then follow from Corollary \ref{CY2-tilting}. 
\end{proof}

\subsubsection{Example~2}\label{Subsec_ex_toricDM2}

\begin{Setup}
	\label{setting_toric_Z^2-2}
	Let $k=\CC$. Let $r$ be a positive integer, and consider the polynomial ring
	\[ S=k[x_1,x_2, x_3, x_4, y_1,\ldots,y_r]. \] 
	Let $a, b, c, d$ be positive integers, and let 
	$a=n_1 a^\prime$, $b=n_2 b^\prime$, $c=n_1 c^\prime$, $d=n_2 d^\prime$ where $n_1=\gcd(a,c)$, $n_2=\gcd(b,d)$. 
	We assume that $c^\prime\le b^\prime$, $\gcd(d^\prime,n_2)=1$, and $d^\prime \,|\,ac^\prime$. 
	Suppose that the algebraic torus $G\coloneqq (k^\times)^2$ acts on $S$ by weights $(a,0), (\frac{ac^\prime}{d^\prime},b), (c,0), (0,d) \in\ZZ_{\ge 0}^2$ 
	on $x_1, x_2, x_3, x_4$ respectively, and $(e_j, f_j)\in\ZZ_{<0}\times\ZZ_{\leq 0}$ on $y_j$ $(j=1,\dots,r)$, 
	where $e_j, f_j$ are integers satisfying \[\Big(a+c+\frac{ac^\prime}{d^\prime},b+d\Big)+\sum_{j=1}^r(e_j,f_j)=(0,0) \quad\text{and}\quad 
	0\le \displaystyle\frac{f_j}{e_j}\le \frac{bd^\prime}{ac^\prime}\] for any $j$. 
	Then, let $R\coloneqq S^G$ be the invariant subring, which is a Gorenstein toric ring of dimension $r+2$. 
    We assume that the weights $(a,0), (\frac{ac'}{d^\prime},b), (c,0),(0,d),(e_j,f_j)$ span the lattice $\Z^2$, otherwise the category $\qmod^{\Z^2}\!S$ splits into several connected components. As before, these components are mutually equivalent, and thus this assumption does not lose generality.
\end{Setup}

To construct NCCRs of the toric ring defined in Setting~\ref{setting_toric_Z^2-2}, we consider the following toric Fano DM stack. 
For the positive integers $a^\prime, b^\prime, c^\prime, d^\prime$ given in Setting~\ref{setting_toric_Z^2-2}, we can take integers $s_1, s_2, t$ such that 
$0\le s_1, s_2\le n_1-1$, $0\le t\le n_2-1$ and $s_1a^\prime+s_2c^\prime\equiv 1 \pmod {n_1}$, $td^\prime\equiv 1 \pmod {n_2}$.
Let
\[ N=\ZZ^2\oplus \ZZ/n_1\ZZ\oplus\ZZ/n_2\ZZ \]
and let $\beta\colon \ZZ^4\to N$ be a morphism such that $\beta_i\coloneqq\beta(\e_i)$ $(i=1, 2, 3,4)$ are given by 
\[\beta_1={}^t(-c^\prime,c^\prime, s_1, 0), \quad \beta_2={}^t(0,-d^\prime, 0, 0), \quad 
\beta_3={}^t(a^\prime,0, s_2, 0),\quad \beta_4={}^t(0,b^\prime, 0, t).\]
Let $\Sigma\subset N_\RR$ be the fan that has four rays generated by $\overline{\beta_1}, \overline{\beta_2}, \overline{\beta_3}, \overline{\beta_4}\in N/N_{\rm tor}$.
Note that $\Sigma$ is the fan of the Hirzebruch surface of degree one. 
For the stacky fan $\bfSigma=(N, \Sigma, \beta)$, the lifting of $\beta$ is given by the matrix ${\small \begin{pmatrix}\beta_1&\beta_2&\beta_3&\beta_4\end{pmatrix}}$. 
We consider the matrix $[B\,|\,Q]$ where $Q={\small \begin{pmatrix}0&0\\0&0\\n_1&0\\0&n_2\end{pmatrix}}$. 
Then, we see that $\DG(\beta)= (\ZZ^{6})^*/\image [B\,|\,Q]^*\cong\ZZ^2$ and the integer kernel of $[B\,|\,Q]$ is generated by 
${}^t{\small \begin{pmatrix}a&\frac{ac^\prime}{d^\prime}&c&0&u&0\\0&b&0&d&0&v\end{pmatrix}}$ 
where $u, v$ are integers with $s_1a+s_2c+un_1=0$, $td+vn_2=0$.
Thus, $\beta^\vee\colon \ZZ^4\to \ZZ^2$ is given by ${\small \begin{pmatrix}a&\frac{ac^\prime}{d^\prime}&c&0\\0&b&0&d\end{pmatrix}}$ 
and this determines the action of $G=\Hom_\ZZ(\DG(\beta),k^\times)\cong(k^\times)^2$ on $U_\Sigma$. 
The toric DM stack defined by the stacky fan $\bfSigma=(N, \Sigma, \beta)$, which we denote by $\calX_{a,b,c,d}$, is Fano and $\Pic(\calX_{a,b,c,d})\cong\ZZ^2$. 
We can take a basis of $\Pic(\calX_{a,b,c,d})$ so that prime divisors $D_1, \dots, D_4$ are represented by 
\[[D_1]=(a,0), \quad [D_2]=\Big(\frac{ac^\prime}{d^\prime},b\Big), \quad [D_3]=(c,0), \quad [D_4]=(0,d)\] 
in $\Pic(\calX_{a,b,c,d})$, and we will fix such a basis in the following argument.

\medskip

Again, the above inputs yield an algebraic setup in Section \ref{quotstack} as follows.
Let $A=k[x_1,x_2,x_3,x_4]$ on which we give a grading by $L=\Z^2$ specified by the weights as in Setting \ref{setting_toric_Z^2-2}. Also, let $I\subset A$ be the ideal $J_{\Si}$ for the fan $\Si$. Then by Proposition \ref{AB}, these satisfy the conditions in Assumption \ref{stack}, and we have $\qmod^{\Z^2}\!A=\coh^{(k^\times)^2}\!(\Spec A\setminus V(I))=\coh\X_{a,b,c,d}$.

Then we consider a tilting bundle of $\calX_{a,b,c,d}$ given in the following form. 

\begin{Prop}
	\label{tilting_toric_Z^2-2}
	Let $a,b,c,d$ be positive integers given in Setting~\ref{setting_toric_Z^2-2} and 
	$\calX=\calX_{a,b,c,d}$ be the toric Fano DM stack as above. 
	For a non-negative integer $q$, we let $r_q\coloneqq \min\{n\in\ZZ_{\ge 0} \mid q<(n+1)\frac{bd^\prime}{ac^\prime}\}$. 
	We consider the set 
	\[\S=\bigcup_{q=0}^{b+d-1}\{(p,q)\in\Pic(\calX_{a,b,c,d})\cong\ZZ^2 \mid r_q\le p\le r_q+a+c-1\}. \]
	Then $T=\bigoplus_{(p,q)\in\S}\O_\calX(p,q)$ is a tilting bundle of $\calX$. 
\end{Prop}

\begin{proof}
	First, let $Q$ be the convex hull of $\S$ in $\RR^2$ and consider the Minkowski sum $\widetilde{Q}=Q+(-Q)$. 
	Then, we see that $\widetilde{\S}=\widetilde{Q}\cap\ZZ^2=\widetilde{\S}_1\cup\widetilde{\S}_2$, where 
	\begin{align*}
		&\widetilde{\S}_1=\bigcup_{q=0}^{b+d-1}\{(p,q)\in\ZZ^2 \mid r_q-(a+c-1)\le p\le m_q+a+c-1\}\\
		&\widetilde{\S}_2=\bigcup_{q=0}^{b+d-1}\{(p,-q)\in\ZZ^2 \mid -m_q-(a+c-1)\le p\le -r_q+a+c-1\}, 
	\end{align*}
	where we let $\displaystyle m_q=\max_{1\le i\le b+d-q-1}\{r_{q+i}-r_i,r_q\}$ for $q=0, 1, \dots, b+d-2$ and $m_{b+d-1}=r_{b+d-1}$. 
	
	We now show that the set of line bundles $\{\O_\calX(p,q) \mid (p,q)\in\S\}$ is a full strong exceptional collection. 
	First, we show that $\Ext_\calX^i(L, L^\prime)=0$ for any $i>0$ and any line bundles $L, L^\prime$ appearing in $T$, 
	which means that any line bundle $\O_\calX(p,q)$ with $(p,q)\in\widetilde{\S}$ is acyclic. 
	By \cite[Proposition~4.5]{BH_Kingconj}, it is enough to show that $(p,q)\in\widetilde{\S}$ does not lie 
	in any of the forbidden cones $F_\varnothing$, $F_{I_+}$, and $F_{I_-}$. 
	Since there is a relation 
	\[(b^\prime+d^\prime)a^\prime\,\overline{\beta_1}-(a^\prime(b^\prime-c^\prime)+b^\prime c^\prime)\,\overline{\beta_2}+(b^\prime+d^\prime)c^\prime\,\overline{\beta_3}-\Big(a^\prime+c^\prime+\frac{a^\prime c^\prime}{d^\prime}\Big)d^\prime\,\overline{\beta_4}={\small \begin{pmatrix}0\\0\end{pmatrix}}\]
	and $a^\prime, d^\prime>0$, $b^\prime\ge c^\prime>0$, we have $I_+=\{1,3\}$, $I_-=\{2,4\}$ and the forbidden cones are 
	\begin{align}\label{eq_fcone_abc}
		\nonumber&F_\varnothing=\Big(-a-c-\frac{ac^\prime}{d^\prime}, -b-d\Big)-\RR_{\ge 0}(1,0)-\RR_{\ge 0}(0,1)\\
		&F_{I_+}=\Big(-\frac{ac^\prime}{d^\prime}, -b-d\Big)+\RR_{\ge 0}(1,0)-\RR_{\ge 0}\Big(\frac{ac^\prime}{d^\prime},b\Big)\\
		\nonumber&F_{I_-}=(-a-c, 0)-\RR_{\ge 0}(1,0)+\RR_{\ge 0}\Big(\frac{ac^\prime}{d^\prime},b\Big), 
	\end{align}
	see Figure~\ref{fig_forbidden_cone}. 
	Then, we will show that any $(p,q)\in\widetilde{\S}$ satisfies $q> -b-d$ and $q< \frac{bd^\prime}{ac^\prime}(p+a+c)$. 
	The first inequality is trivial. To show the latter one, we note that $r_q$ satisfies $q<\frac{bd^\prime}{ac^\prime}(r_q+1)$ and $q\ge \frac{bd^\prime}{ac^\prime}\,r_q$ by definition. Then, we have the following. 
	\begin{itemize}
		\setlength{\parskip}{0pt} 
		\setlength{\itemsep}{3pt}
		\item If $p\ge r_q-(a+c-1)$, then $\frac{bd^\prime}{ac^\prime}(p+a+c)\ge \frac{bd^\prime}{ac^\prime}(r_q+1)>q$, and hence any $(p,q)\in\widetilde{\S}_1$ satisfies the latter inequality. 
		\item If $p\ge -m_q-(a+c-1)$, then $\frac{bd^\prime}{ac^\prime}(p+a+c)\ge \frac{bd^\prime}{ac^\prime}(1-m_q)$. 
		\begin{itemize}
			\setlength{\parskip}{0pt} 
			\setlength{\itemsep}{3pt}
			\item[$\ast$] For the case where $m_q=r_q$, we have 
			\[\frac{bd^\prime}{ac^\prime}(1-m_q)=\frac{bd^\prime}{ac^\prime}(1-r_q)\ge \frac{bd^\prime}{ac^\prime}-q>-q.\] 
			\item[$\ast$] For the case where $m_q=r_{q+i}-r_i$, we have 
			\[\frac{bd^\prime}{ac^\prime}(1-m_q)=\frac{bd^\prime}{ac^\prime}(1-r_{q+i}+r_i)=\frac{bd^\prime}{ac^\prime}(r_i+1)-\frac{bd^\prime}{ac^\prime}r_{q+i}>i-(q+i)=-q.\]
		\end{itemize}
		Thus, any $(p,-q)\in\widetilde{\S}_2$ satisfies the latter inequality. 
	\end{itemize}
	In conclusion, any element in $\widetilde{\S}$ does not lie in any of $F_\varnothing$, $F_{I_+}$, and $F_{I_-}$,
	and hence $\{\O_\calX(p,q) \mid (p,q)\in\S\}$ forms a strong exceptional collection. 
	
	Also, by the inductive argument in $\Pic(\calX)$ similar to the proof of Proposition~\ref{tilting_toric_Z^2-1}, 
	we see that all line bundles are contained in the subcategory $\D$ of $\Db(\coh \calX)$ generated by the line bundles $\O_\calX(p,q)$ with $(p,q)\in\S$. 
	By \cite[Corollary~4.8]{BH_Ktheory}, $T$ is a generator of $\Db(\coh \calX)$.
\end{proof}

\begin{figure}[H]
	\begin{center}
		\scalebox{0.45}{
			\begin{tikzpicture}
				\newcommand{\edgewidth}{0.07cm} 
				\newcommand{\noderad}{0.13} 
				\newcommand{\lettersize}{\huge} 
				
				\draw (0,0) circle [ultra thick, radius=0.25]; 
				\filldraw[fill=black] (0,0) circle [radius=\noderad]; 
				
				\foreach \a/\b/\n in {-7/0/a, -5/-5/b, -12/-5/c, -16/0/a-, 3/6/a+, -12/-9.2/b-, 4/-5/b+, -16/-5/c-, -12/-9.2/c+} {\coordinate (F\n) at (\a,\b); }
				\filldraw [lightgray!60] (Fa-)--(Fa)--(Fa+)--(-16,6)--(Fa-); 
				\filldraw [lightgray!60] (Fb+)--(Fb)--(Fb-)--(4,-9.2)--(Fb+); 
				\filldraw [lightgray!60] (Fc-)--(Fc)--(Fc+)--(-16,-9.2)--(Fc-); 
				\foreach \n in {a,b,c} {\node at (F\n) {\huge$\times$}; }
				\node at (-12,4) {\lettersize$F_{I_-}$}; \node at (-2,-8) {\lettersize$F_{I_+}$}; \node at (-14,-7) {\lettersize$F_\varnothing$}; 
				\node at (-4.8,3.3) {\lettersize$q=\displaystyle \frac{bd^\prime}{ac^\prime}(p+a+c)$}; \node at (0,2.-7.5) {\lettersize$q=-b-d$}; 
				\node at (-11,0.5) {\lettersize$q=0$}; \node at (-6.3,-7) {\lettersize$q=\displaystyle \frac{bd^\prime}{ac^\prime}\,p-d$}; 
				
				\path (Fa) ++(-90:0.6cm) node {\lettersize$(-a-c,0)$}; \path (Fb) ++(90:0.7cm) node {\lettersize$(-\frac{ac^\prime}{d^\prime},-b-d)$};
				\path (Fc) ++(90:0.7cm) node {\lettersize$(-a-c-\frac{ac^\prime}{d^\prime},-b-d)$};
				\foreach \s/\t in {a/a-, a/a+, b/b-, b/b+, c/c-, c/c+}{\draw [line width=\edgewidth] (F\s)--(F\t);}
				
		\end{tikzpicture}}
	\end{center}
	\caption{Forbidden cones $F_\varnothing$, $F_{I_+}$, and $F_{I_-}$ for $\calX_{a,b,c,d}$, where the double circle stands for the origin}
	\label{fig_forbidden_cone}
\end{figure}

Now we consider the total space of a vector bundle $E=\bigoplus_{j=1}^r\O_\calX(e_j,f_j)\in\coh\X$, whose total space $\Y$ is given as $\coh\Y=\qmod^{\Z^2}\!B$ for $B=A[y_1,\ldots,y_r]=S$ with $\deg y_j=(e_j,f_j)$.
\begin{Thm}
	\label{NCCR_toric_Z^2-2}
	Let $\calX=\calX_{a,b, c,d}$ be the toric Fano DM stack as above 
	and $T=\bigoplus_{(p,q)\in\S}\O_\calX(p,q)$ be the tilting bundle of $\calX$ given in Proposition~\ref{tilting_toric_Z^2-2}. 
	For the vector bundle $E=\bigoplus_{j=1}^r\O_\calX(e_j,f_j)$, we consider the total space $\Y$ of $E$ and the map $\pi\colon \Y\rightarrow \calX$, 
	where $e_j$ and $f_j$ are the weights defined in Setting~\ref{setting_toric_Z^2-2}.
\begin{enumerate}
\item $\pi^\ast T$ is a tilting bundle on $\Y$ and $\Ga=\End_\Y(\pi^\ast T)=\bigoplus_{i\geq0}\Hom_\X(T,T\otimes_\X\Sym^iE^\vee)$ is an $(2+r)$-CY algebra of Gorenstein parameter $r$.
\item $\Ga$ is an NCCR of the toric ring $R$ given in Setting~\ref{setting_toric_Z^2-2}.
\end{enumerate}
\end{Thm}

\begin{proof}
	The proof is similar to the one of Theorem~\ref{NCCR_toric_Z^2-1}. 
	In particular, it suffices to show that any lattice point of the form $(p-\sum_{j=1}^rm_je_j, q-\sum_{j=1}^rm_j f_j)\in\ZZ^2$, 
	where $m_j \in\ZZ_{\ge 0}$ and $(p,q)\in\widetilde{\S}$, 
	does not lie in the forbidden cones $F_\varnothing$, $F_{I_+}$, $F_{I_-}$ shown in \eqref{eq_fcone_abc}, see also Figure~\ref{fig_forbidden_cone}. 
	This claim follows from the fact that any $(p,q)\in\widetilde{\S}$ satisfies $q> -b-d$ and $q< \frac{bd^\prime}{ac^\prime}(p+a+c)$
	(see the proof of Proposition~\ref{tilting_toric_Z^2-2}) and the assumption $0\le\frac{f_j}{e_j}\le \frac{bd^\prime}{ac^\prime}$. 
    This shows that $\pi^\ast T\in\Db(\coh\Y)$ is a tilting object. 
    The remaining assertions then follow from Corollary \ref{CY2-tilting}. 
\end{proof}

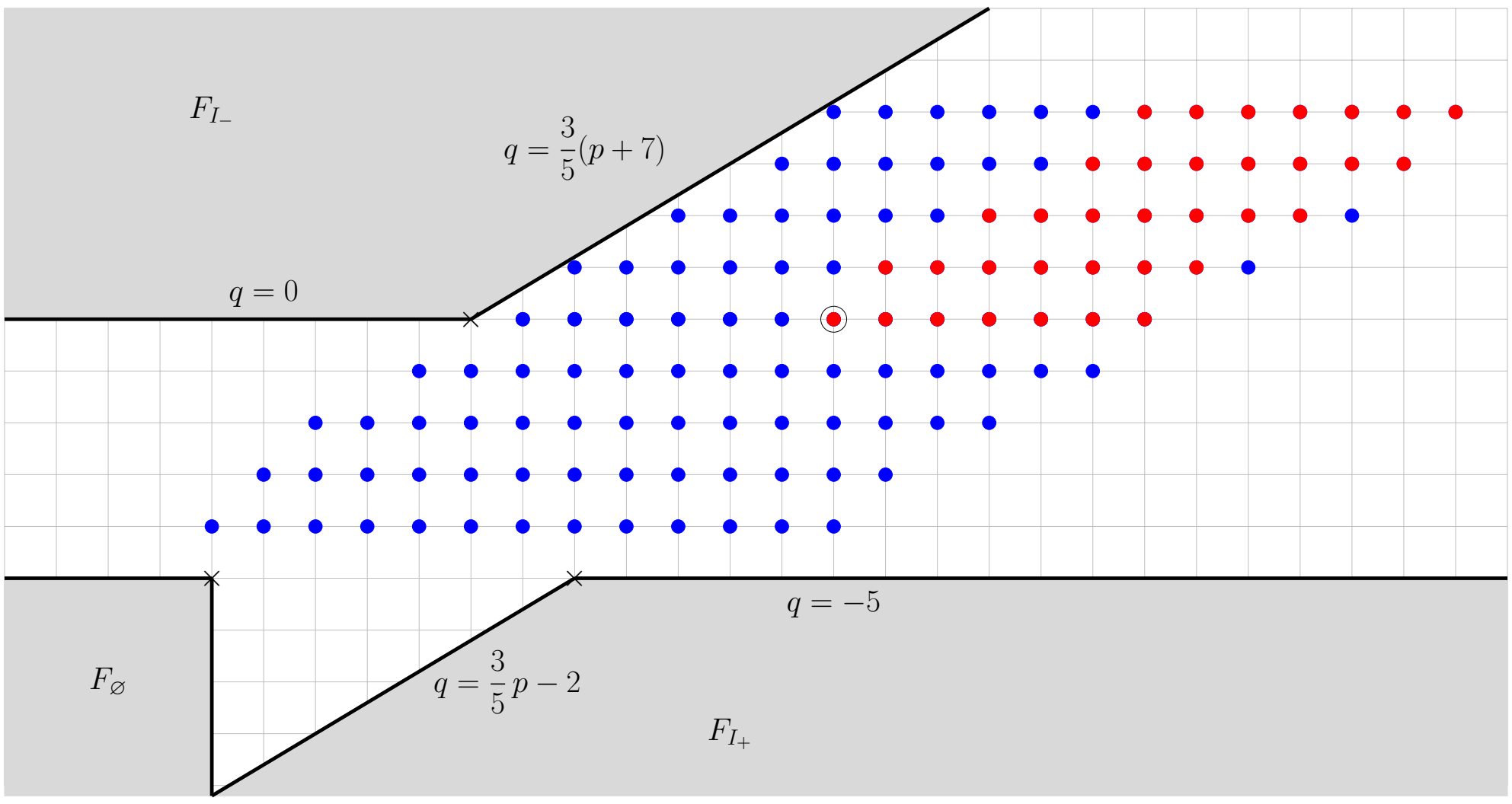
\begin{figure}[H]
	\begin{center}
		\scalebox{0.45}{
			\begin{tikzpicture}
				\newcommand{\edgewidth}{0.07cm} 
				\newcommand{\noderad}{0.13} 
				\newcommand{\lettersize}{\huge} 
				
				\draw [step=1, gray!50] (-16,-9) grid (13,6); 
				\draw (0,0) circle [ultra thick, radius=0.25]; 
				
				\foreach \a/\b/\n in {-7/0/a, -5/-5/b, -12/-5/c, -16/0/a-, 3/6/a+, -12/-9.2/b-, 13/-5/b+, -16/-5/c-, -12/-9.2/c+} {\coordinate (F\n) at (\a,\b); }
				\filldraw [lightgray!60] (Fa-)--(Fa)--(Fa+)--(-16,6)--(Fa-); 
				\filldraw [lightgray!60] (Fb+)--(Fb)--(Fb-)--(13,-9.2)--(Fb+); 
				\filldraw [lightgray!60] (Fc-)--(Fc)--(Fc+)--(-16,-9.2)--(Fc-); 
				\foreach \n in {a,b,c} {\node at (F\n) {\huge$\times$}; }
				\node at (-12,4) {\lettersize$F_{I_-}$}; \node at (-2,-8) {\lettersize$F_{I_+}$}; \node at (-14,-7) {\lettersize$F_\varnothing$}; 
				\node at (-4.8,3.3) {\lettersize$q=\displaystyle \frac{3}{5}(p+7)$}; \node at (0,2.-7.5) {\lettersize$q=-5$}; 
				\node at (-11,0.5) {\lettersize$q=0$}; \node at (-6.3,-7) {\lettersize$q=\displaystyle \frac{3}{5}\,p-2$}; 
				
				\foreach \s/\t in {a/a-, a/a+, b/b-, b/b+, c/c-, c/c+}{\draw [line width=\edgewidth] (F\s)--(F\t);}
				\foreach \b/\r in {0/0,1/1,2/3,3/5,4/6}{ 
					\foreach \a in {-6,-5,-4,-3,-2,-1,0,1,2,3,4,5,6}{\filldraw[blue] (\a+\r,\b) circle [radius=\noderad]; }
					\foreach \a in {-6,-5,-4,-3,-2,-1,0,1,2,3,4,5,6}{\filldraw[blue] (-\a-\r,-\b) circle [radius=\noderad]; }}
				\filldraw[blue] (8,1) circle [radius=\noderad]; \filldraw[blue] (10,2) circle [radius=\noderad]; 
				\filldraw[blue] (-8,-1) circle [radius=\noderad]; \filldraw[blue] (-10,-2) circle [radius=\noderad]; 
				\foreach \b/\r in {0/0,1/1,2/3,3/5,4/6}{ 
					\foreach \a in {0,1,2,3,4,5,6}{\filldraw[red] (\a+\r,\b) circle [radius=\noderad]; }}
		\end{tikzpicture}}
	\end{center}
	\caption{Forbidden cones $F_\varnothing$, $F_{I_+}$, and $F_{I_-}$ for $\calX_{5,3,2,2}$, where the double circle stands for the origin.  
		The red dots are elements in $\S$ which produce the tilting bundle $T$, and the red and blue dots are elements in $\widetilde{\S}$.}
	\label{fig_forbidden_cone532}
\end{figure}

\begin{Ex}
	We here consider the toric Fano DM orbifold $\calX_{a,b,c,d}$ for $a=5, b=3, c=2, d=2$. 
	By Proposition~\ref{tilting_toric_Z^2-2}, $\calX=\calX_{5,3,2,2}$ admits the tilting bundle $T=\bigoplus_{(p,q)\in\S}\O_\calX(p,q)$, where 
	\[\S=\bigcup_{q=0}^{4}\{(p,q)\in\ZZ^2 \mid r_q\le p\le r_q+6\},
	\quad r_0=0, \, r_1=1, \, r_2=3, \, r_3=5, \, r_4=6. \]
	Figure~\ref{fig_forbidden_cone532} shows lattice points in $\S$ and the forbidden cones. 
	For $(e_j, f_j)\in \ZZ_{<0}\times\ZZ_{\leq 0}$ such that $\sum_{j=1}^r(e_j,f_j)=(-12,-5)$ and $0\le \displaystyle\frac{f_j}{e_j}\le \frac{3}{5}$ for any $j$, 
	we consider the vector bundle $E=\bigoplus_{j=1}^r\O_\calX(e_j,f_j)$ and the total space $\Y$ of $E$. 
	Then, by Theorem~\ref{NCCR_toric_Z^2-2}, $\pi^\ast T$ is a tilting bundle on $\Y$ and $\End_\Y(\pi^\ast T)$ is an NCCR of the toric ring 
	defined as the invariant subring of $S=\CC[x_1,x_2, x_3, x_4, y_1,\ldots,y_r]$ under the action of $(\CC^\times)^2$ 
	by the weights $(5,0), (5,3), (2,0), (0,2)\in\ZZ_{\ge 0}^2$ for $x_1, x_2, x_3, x_4$ respectively, and $(e_j, f_j)\in\ZZ_{\leq 0}^2$ for $y_j$ $(j=1,\dots,r)$.

	There are some choices of non-positive integers $e_j, f_j$ satisfying the above conditions. For example, if $r=2$, then these integers can be chosen as 
	{\small\begin{align*}
			\{(e_1,f_1), (e_2,f_2)\}=&\{(-1,0), (-11,-5)\}, \{(-2,0), (-10,-5)\}, \{(-2,-1), (-10,-4)\}, \{(-3,0), (-9,-5)\},\\ 
			&\{(-3,-1), (-9,-4)\}, \{(-4,-1), (-8,-4)\}, \{(-4,-2), (-8,-3)\}, \\
			&\{(-5,-1), (-7,-4)\}, \{(-5,-2), (-7,-3)\}, \{(-6,-2), (-6,-3)\}. 
	\end{align*}}
	We also write down some choices of these integers for the case $r=3$: 
	\begin{align*}
		\{(e_j,f_j)\}_{j=1,2,3}=&\{(-1,0), (-2,0), (-9,-5)\}, \{(-3,-1), (-4,-2), (-5,-2)\}, \\
		&\{(-2,-1), (-5,-2), (-5,-2)\}, \{(-4,-1), (-4,-2), (-4,-2)\}, \,\, \text{etc.} 
	\end{align*}
\end{Ex}

\section{NCCRs arising from Geigle--Lenzing projective spaces}
\label{sec_GLspace}

Let us first recall some fundamentals on Geigle--Lenzing theory \cite{GL} and its higher dimensional generalization \cite{HIMO}.

Let $k$ be a field, let $d\geq0$ and $n\geq1$ be integers, let $p=(p_1,\ldots,p_n)$ be a sequence of positive integers and consider the ring
\[ R=k[T_0,\ldots,T_d,X_1,\ldots,X_n]/(X_1^{p_1}-l_1,\ldots,X_n^{p_n}-l_n), \]
where $l_1,\ldots,l_n$ are linear forms in $T_0,\ldots,T_d$ taken as linearly independent as possible, that is, every subset of $l_1,\ldots,l_n$ with at most $d+1$ elements is linearly independent.
The ring $R$ is a complete intersection of dimension $d+1$, which is graded by the abelian group
\[ \LL:=\Z\vec{c}\oplus\textstyle{\bigoplus}_{i=1}^n\Z\vec{x}_i/(p_1\vec{x}_1=\cdots=p_n\vec{x}_n=\vec{c}) \]
by setting $\deg T_i=\vec{c}$ and $\deg X_i=\vec{x}_i$.
This group $\LL$ is naturally an ordered group, by requiring that $\vec{x}\in\LL$ satisfies $\vec{x}\geq0$ if and only if $\vec{x}$ lies in 
the submonoid $\LL_+$ of $\LL$ generated by $(\vec{x}_i)_{i}$ and $\vec{c}$.

The fact that the ring $R$ has an $\LL$-grading amounts to saying that the scheme $\Spec R$ has an action of a reductive algebraic group $G:=\Spec k[\LL]$.
The {\it Geigle--Lenzing (GL) projective space} is the quotient stack
\[ \XX=[(\Spec R\setminus \{R_+\})/G] \]
for the $G$-invariant maximal ideal $R_+=(T_0,\ldots,T_d,X_1,\ldots,X_n)$.
Therefore, the category of coherent sheaves $\coh \XX$ on $\XX$ is nothing but the category of $G$-equivariant sheaves on $\Spec R\setminus\{R_+\}$. Note that $G$-equivariant coherent sheaves on $\Spec R$ (resp. supported on $\Spec R_0$) are just the $\LL$-graded finitely generated modules over $R$ (resp. finite length $\LL$-graded $R$-modules). 
We denote by $\mod^\LL\!R$ the former and by $\fl^\LL\!R$ the latter, and we thus have
\[ \coh \XX=\qmod^\LL\!R:=\mod^\LL\!R/\fl^\LL\!R. \]
In other words, one can view the category $\coh \XX$ as an ($\LL$-graded) (non-)commutative projective scheme \cite{Ser,AZ}.
Let $\vec{\om}=(n-d-1)\vec{c}-\sum_{i=1}^n\vec{x}_i\in\LL$ and $\om_\XX=\O_\XX(\vec{\om})$.

The following result shows that the category $\Db(\coh \XX)$ always has a tilting bundle, whose endomorphism ring $A$ is called the {\it $d$-canonical algebra}.
\begin{Thm}[{\cite[Theorem 1.13]{HIMO}}]\label{GLtilting}
	Let $\XX$ be a GL projective space defined as above. Then $T\coloneqq\bigoplus_{0\leq\vec{x}\leq d\vec{c}}\O_\XX(\vec{x})$ is a tilting bundle whose endomorphism ring $A:=\End_\XX(T)$ has finite global dimension. Therefore we have a triangle equivalence $\Db(\coh \XX)\simeq\Db(\mod A)$.
\end{Thm}

Suppose in what follows $n=d+1\geq2$, so that we have
\begin{equation}\label{R_n=d+1}
R=k[X_1,\ldots,X_n] \end{equation}
graded by the abelian group $\LL=\bigoplus_{i=1}^n\Z\vec{x}_i/(p_1\vec{x}_1=\cdots=p_n\vec{x}_n\eqqcolon \vec{c})$. 
In this case, $\vec{\om}=-\sum_{i=1}^n\vec{x}_i\in\LL$, and the $d$-canonical algebra $A$ defined in Theorem~\ref{GLtilting} is a $d$-representation infinite algebra (see \cite[Theorem~1.14]{HIMO}).

Let us observe that in the case $n=d+1$, a GL projective space $\XX$ can be described as a toric DM stack as follows. 
Let $N=\Z^d=\Z^{n-1}$ and let $\Si\subset N_\RR\cong\RR^{n-1}$ be the complete fan whose rays are $\{\rho_1,\ldots,\rho_n\}$, where $\rho_i$ is spanned by the unit vector $\e_i$ for $1\leq i\leq n-1$, and $\rho_n$ by $-\e_1-\cdots-\e_{n-1}$. Then define the map $\be\colon\Z^n\to N=\Z^{n-1}$ by the matrix
\[ \be=\begin{pmatrix} p_1&&&&-p_n \\ & p_2&&&-p_n \\ &&\ddots&&\vdots \\ &&& p_{n-1}&-p_n\end{pmatrix}. \]
It is easily seen that the  toric DM stack $\calX_\bSi$ associated to the stacky fan $\bSi=(N,\Si,\be)$
is the quotient stack of $\AA^n\setminus\{0\}$ by the reductive group $G=\{ (t_1,\ldots,t_n)\in (k^\times)^n \mid t_1^{p_1}=\cdots=t_n^{p_n}\}$ whose character group is $\DG(\beta)=\LL$. 
It follows that $\calX_\bSi$ is nothing but the GL projective space associated with $R$ \eqref{R_n=d+1}, and $\coh\X_{\bSi}=\qmod^\LL\!R$. 

By Proposition \ref{AB}, this setup satisfies Assumption \ref{stack} for the algebra $A=R\cong k[X_1,\ldots,X_n]$ and the maximal ideal $I=(X_1,\ldots,X_n)$.

Now, consider a vector bundle $E=\bigoplus_{j=1}^r\O_\XX(-\vec{y}_j)$ of rank $r$ on $\XX$ such that $\vec{y}_1, \dots, \vec{y}_r\in\LL_+$. Suppose also that $\vec{x}_1+\cdots+\vec{x}_n=\vec{y}_1+\cdots+\vec{y}_r$ so that we have $\Wedge^rE=\om_\XX$. The following result is a generalization of \cite[Corollary 1.20]{HIMO} which studies the case $E=\om_\XX$ and produces an NCCR of the Veronese subring $R^{(\vec{\om})}$. 
\begin{Thm}\label{NCCRfromGL}
Let $\XX$ be a GL projective space for $n=d+1$ and $p=(p_1, \dots, p_n)$, let $E$ be a vector bundle of rank $r$ as above 
with $\pi\colon \YY\to \XX$ the total space. 
Let $T=\bigoplus_{0\leq\vec{x}\leq d\vec{c}}\O_\XX(\vec{x})$ be the tilting bundle on $\XX$ given in Theorem~\ref{GLtilting}.
\begin{enumerate}
\item $\pi^\ast T$ is a tilting bundle on $\YY$, and 
$\Ga=\End_\YY(\pi^\ast T)$ is $(d+r)$-CY of Gorenstein parameter $r$.
\item $\Ga$ is an NCCR of the section ring $\bigoplus_{i\geq0}\Hom_\XX(\O_\XX,\O_\XX\otimes\Sym^iE^\vee)$.
\end{enumerate}
\end{Thm}

\begin{proof}
We first show that $\pi^\ast T$ is a tilting bundle on $\YY$ (the proof is similar to \cite[Theorem~3.13]{IW_WPL}). 
Since $\pi^\ast T\in\Db(\coh\Y)$ is a generator, we only have to show that $\Ext^p_\YY(\pi^\ast T, \pi^\ast T)=0$ for any $p>0$. By the adjunction and the projection formula, we have 
\begin{align*}
\Ext^p_\YY(\pi^\ast T, \pi^\ast T)&\cong\Ext^p_\XX(T, \pi_\ast\pi^\ast T)\cong\Ext^p_\XX(T, T\otimes \Sym E^\vee) \\
&=\Ext^p_\XX(T, T\otimes \bigoplus_{i\ge 0}\Sym^i E^\vee)=
\Ext^p_\XX\Big(T, T\otimes \bigoplus_{i\ge 0}\bigoplus_{\substack{(m_1, \dots, m_r)\in\ZZ_{\ge 0}^r\\ 
\sum_{j=1}^r m_j=i}}\O_\XX\Big(\sum_{j=1}^r m_j\vec{y}_j\Big)\Big)\\
&\cong\bigoplus_{i\ge 0}\bigoplus_{\substack{(m_1, \dots, m_r)\in\ZZ_{\ge 0}^r\\ \sum_{j=1}^r m_j=i}}\Ext^p_\XX\Big(T, T\otimes \O_\XX\Big(\sum_{j=1}^r m_j\vec{y}_j\Big)\Big). 
\end{align*}
By \cite[Proposition~5.5]{HIMO}, for $p>0$, we have 
\[
\Ext^p_\XX\Big(T, T\otimes \O_\XX\Big(\sum_{j=1}^r m_j\vec{y}_j\Big)\Big)
=\begin{cases}
D\Big(\displaystyle\bigoplus_{0\le s\le d\vec{c}}\,\bigoplus_{0\le t\le d\vec{c}} R_{\vec{\om}-s+t-\sum m_j\vec{y}_j}\Big)  & (p=d)\\
0& (p>0, p\neq d), 
\end{cases}
\]
which is concentrated in degree $\geq -\vec{\omega}-d\vec{c}$ (where the possible lower bound is attained when $i=0$, $m_1=\cdots=m_r=0$, $s=0$, $t=d\vec{c}$). 
Since $n=d+1$ and $\vec{c}=p_1\vec{x}_1=\cdots=p_n\vec{x}_n$, we see that 
\[
\vec{\om}+d\vec{c}=d\vec{c}-(\vec{x}_1+\cdots+\vec{x}_n)=(p_1-1)\vec{x}_1+\cdots+(p_{n-1}-1)\vec{x}_{n-1}-\vec{x}_n\not\in\LL_+, 
\]
and hence $0\not\geq -\vec{\omega}-d\vec{c}$. 
In conclusion, we have $\Ext^p_\XX\Big(T, T\otimes \O_\XX\Big(\sum_{j=1}^r m_j\vec{y}_j\Big)\Big)=0$ for any $p>0$. 
Thus, $\pi^\ast T$ is a tilting bundle on $\YY$, and the remaining assertions follow from Corollary \ref{CY2-tilting}.
\end{proof}

\begin{Cor}\label{NCCRfromGL2}
Let the notation be the same as in Theorem~\ref{NCCRfromGL}. 
Let $R=k[X_1,\ldots,X_n]$ be the $\LL$-graded ring given as \eqref{R_n=d+1}. 
We consider the $\LL$-graded ring $S\coloneqq R[Y_1, \dots, Y_r]$ defined by $\deg Y_j=-\vec{y}_j$ $(j=1, \dots, r)$. 
Then, the degree zero part $S_0$ is isomorphic to the section ring $\bigoplus_{i\geq0}\Hom_\XX(\O_\XX,\O_\XX\otimes\Sym^iE^\vee)$, 
and hence it admits an NCCR. 
\end{Cor}

\begin{proof}
We see that $\Sym^iE^\vee=\bigoplus\O_\XX(m_1\vec{y}_1+\cdots+m_r\vec{y}_r)$, where the direct sum runs over all $(m_1, \dots, m_r)\in\ZZ_{\ge 0}^r$ 
with $\sum_{j=1}^r m_j=i$. 
Thus, $\Hom_\XX(\O_\XX,\O_\XX\otimes\Sym^iE^\vee)$ consists of morphisms from $\O_\XX$ in $\coh\XX$ whose degree coincides with 
$m_1\vec{y}_1+\cdots+m_r\vec{y}_r$ in $\LL$. 
Since $\vec{x}_1, \dots, \vec{x}_n, \vec{y}_1, \dots, \vec{y}_r\in\LL_+$, generators of the section ring correspond to those of $S_0$ and vice versa. 
\end{proof}

\begin{Rem}
In Corollary~\ref{NCCRfromGL2}, we can consider the $\LL$-grading of $S$ as (the weight of) the action of $G=\Spec k[\LL]$ on $S$, 
in which case $S_0$ is the invariant subring $S^G$, and it is a toric ring. 
On the other hand, let $R^\prime\coloneqq k[Y_1^\prime, \dots, Y_r^\prime]$ be the $\LL$-graded ring defined by 
$\deg Y_j^\prime=\vec{y}_j$ $(j=1, \dots, r)$. 
Then, $S_0=S^G$ is isomorphic to the Segre product $R\seg R^\prime=\bigoplus_{l\in\LL}R_l\otimes R^\prime_l$. 
\end{Rem}

\begin{figure}[H]
\begin{center}
\begin{tikzpicture}[sarrow/.style={black, -latex},tarrow/.style={black, latex-}]
\newcommand{\arrowwidth}{0.03cm} 

\node at (0,0) {\scalebox{0.75}{
\begin{tikzpicture}[sarrow/.style={black, -latex},tarrow/.style={black, latex-}]

\foreach \a/\b/\n in {0/0/0, 2/1/\vec{x}_1, 4/1/2\vec{x}_1, 6/1/\cdots, 9/1/(p_1-1)\vec{x}_1, 
2/-1/\vec{x}_2, 4/-1/2\vec{x}_2, 6/-1/\cdots, 9/-1/(p_2-1)\vec{x}_2, 12/0/\vec{c}}{\node (V\a\b) at (\a,\b*1.5) {$\n$};}; 

\foreach \a/\b in {00/21,21/41,41/61,61/91,91/120} 
{\draw[sarrow, line width=\arrowwidth] (V\a)--(V\b)node[midway,xshift=0cm,yshift=0.25cm] {\tiny $X_1$}; };  
\foreach \a/\b in {00/2-1,2-1/4-1,4-1/6-1,6-1/9-1,9-1/120} 
{\draw[sarrow, line width=\arrowwidth] (V\a)--(V\b)node[midway,xshift=0cm,yshift=0.25cm] {\tiny $X_2$}; };  

\end{tikzpicture}
} }; 

\node at (0,-3.5) {\scalebox{0.75}{
\begin{tikzpicture}[sarrow/.style={black, -latex},tarrow/.style={black, latex-}]

\foreach \a/\b/\n in {0/0/0, 2/1/\vec{x}_1, 4/1/2\vec{x}_1, 6/1/\cdots, 9/1/(p_1-1)\vec{x}_1, 
2/-1/\vec{x}_2, 4/-1/2\vec{x}_2, 6/-1/\cdots, 9/-1/(p_2-1)\vec{x}_2, 12/0/\vec{c}}{\node (V\a\b) at (\a,\b*1.5) {$\n$};}; 

\begin{scope}[transform canvas={yshift=0.1cm}]
\foreach \a/\b in {00/21,21/41,41/61,61/91,91/120} 
{\draw[sarrow, line width=\arrowwidth] (V\a)--(V\b)node[midway,xshift=0cm,yshift=0.25cm] {\tiny $X_1$}; };  
\foreach \a/\b in {00/2-1,2-1/4-1,4-1/6-1,6-1/9-1,9-1/120} 
{\draw[sarrow, line width=\arrowwidth] (V\a)--(V\b)node[midway,xshift=0cm,yshift=0.25cm] {\tiny $X_2$}; };  
\end{scope}

\begin{scope}[transform canvas={yshift=-0.1cm}]
\foreach \a/\b in {00/21,21/41,41/61,61/91,91/120} {\draw[tarrow, line width=\arrowwidth,red] (V\a)--(V\b)
node[midway,xshift=0cm,yshift=-0.3cm] {\tiny $X_2Y$};}; 
\foreach \a/\b in {00/2-1,2-1/4-1,4-1/6-1,6-1/9-1,9-1/120} {\draw[tarrow, line width=\arrowwidth,red] (V\a)--(V\b)
node[midway,xshift=0cm,yshift=-0.3cm] {\tiny $X_1Y$};};
\end{scope}

\end{tikzpicture}
}}; 

\node at (0,-7.5) {\scalebox{0.75}{
\begin{tikzpicture}[sarrow/.style={black, -latex},tarrow/.style={black, latex-}]

\foreach \a/\b/\n in {0/0/0, 2/1/\vec{x}_1, 4/1/2\vec{x}_1, 6/1/\cdots, 9/1/(p_1-1)\vec{x}_1, 
2/-1/\vec{x}_2, 4/-1/2\vec{x}_2, 6/-1/\cdots, 9/-1/(p_2-1)\vec{x}_2, 12/0/\vec{c}}{\node (V\a\b) at (\a,\b*1.5) {$\n$};}; 

\begin{scope}[transform canvas={yshift=0.1cm}]
\foreach \a/\b in {00/21,21/41,41/61,61/91,91/120} 
{\draw[sarrow, line width=\arrowwidth] (V\a)--(V\b)node[midway,xshift=0cm,yshift=0.25cm] {\tiny $X_1$}; };  
\foreach \a/\b in {00/2-1,2-1/4-1,4-1/6-1,6-1/9-1,9-1/120} 
{\draw[sarrow, line width=\arrowwidth] (V\a)--(V\b)node[midway,xshift=0cm,yshift=0.25cm] {\tiny $X_2$}; };  
\end{scope}

\begin{scope}[transform canvas={yshift=-0.1cm}]
\foreach \a/\b in {00/21,21/41,41/61,61/91} {\draw[tarrow, line width=\arrowwidth,red] (V\a)--(V\b)
node[midway,xshift=0cm,yshift=-0.25cm] {\tiny $Y_1$};}; 
\draw[tarrow, line width=\arrowwidth,red] (V91)--(V120) node[midway,xshift=0cm,yshift=-0.15cm] {\tiny $Y_1$};
\foreach \a/\b in {00/2-1,2-1/4-1,4-1/6-1,6-1/9-1,9-1/120} {\draw[tarrow, line width=\arrowwidth,red] (V\a)--(V\b)
node[midway,xshift=0cm,yshift=-0.25cm] {\tiny $Y_2$};};
\end{scope}

\foreach \a in {21,41,91} {\draw [line width=\arrowwidth,sarrow, red] (V\a) to [in=60,out=120,loop] ();};  
\foreach \a in {21,41,91} {
\path (V\a) ++(50:1.1cm) coordinate (V\a L); \node[red] at (V\a L) {\tiny $X_2Y_2$}; }
\foreach \a in {2-1,4-1,9-1} {\draw [line width=\arrowwidth,sarrow, red] (V\a) to [in=300,out=240,loop] ();}; 
\foreach \a in {2-1,4-1,9-1} {
\path (V\a) ++(-50:1.1cm) coordinate (V\a L); \node[red] at (V\a L) {\tiny $X_1Y_1$}; }
\end{tikzpicture}
}}; 

\end{tikzpicture}
\caption{The quiver descriptions of $\End_\XX(T)$, $\End_{\YY_1}(\pi_1^\ast T)$ and $\End_{\YY_2}(\pi_2^\ast T)$}
\label{fig_canonicalalg}
\end{center}
\end{figure}
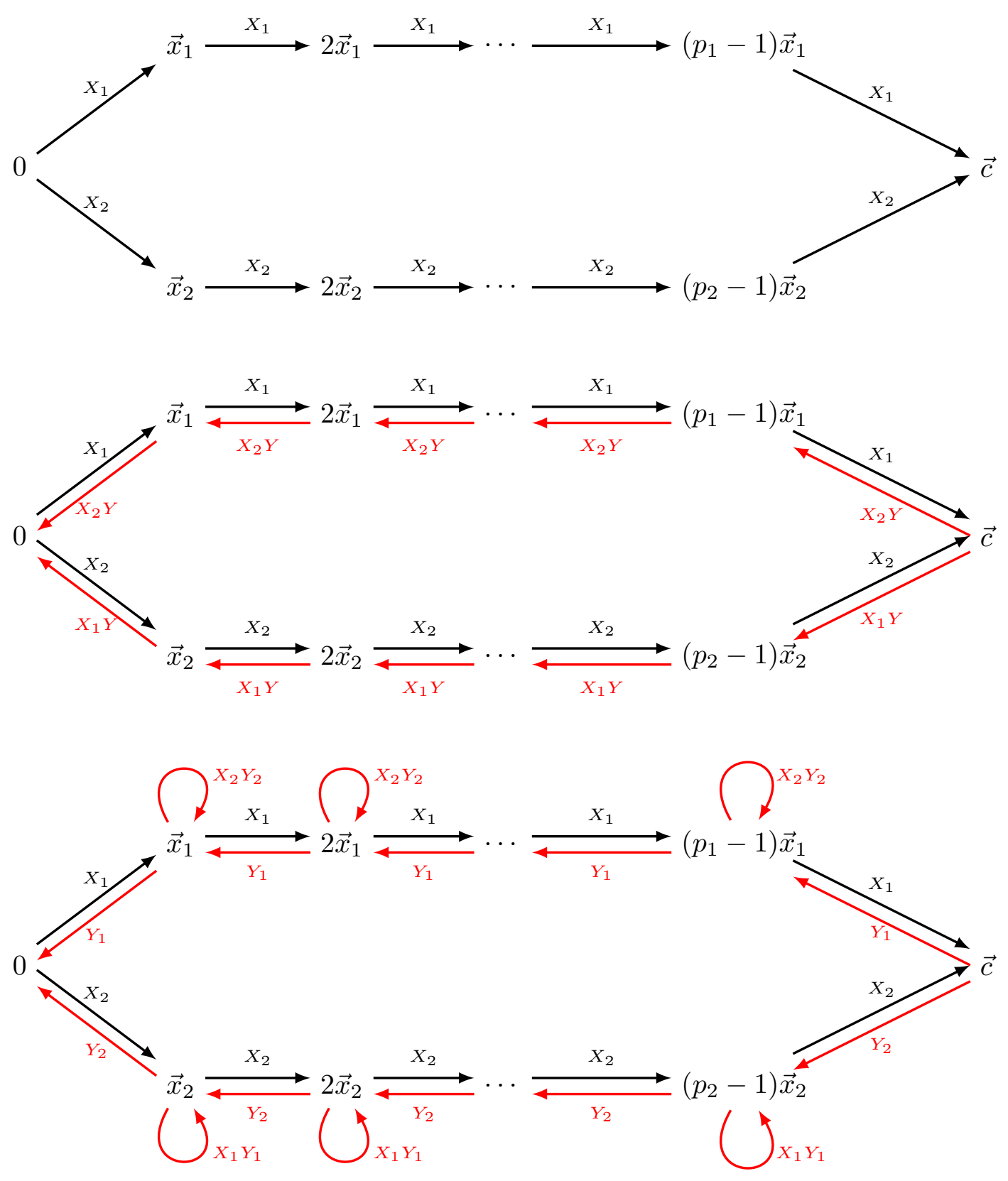

\begin{Ex}\label{Ex_WPL}
For $n=d+1=2$ and $p=(p_1,p_2)$, we have the ring $R=k[X_1,X_2]$ graded by $\LL=\Z\vec{x}_1\oplus\Z\vec{x}_2/(p_1\vec{x}_1=p_2\vec{x}_2\eqqcolon\vec{c})$ 
by setting $\deg X_i=\vec{x}_i$ $(i=1,2)$. Let $\vec{\om}=-\vec{x}_1-\vec{x}_2$. 
In this case, the associated GL projective space $\XX$ is the \emph{weighted projective line} in the sense of \cite{GL}. 
By Theorem~\ref{GLtilting} (see also \cite[Proposition~4.1]{GL}), $\XX$ has a tilting bundle with the form $T=\bigoplus_{0\leq\vec{x}\leq \vec{c}}\O_\XX(\vec{x})$. 
\begin{itemize}
\setlength{\parskip}{0pt} 
\setlength{\itemsep}{3pt}
\item[(1)] 
For the vector bundle $E_1=\om_\XX=\O_\XX(\vec{\om})$ and its total space $\pi_1\colon \YY_1\coloneqq\Spec (\Sym E_1^\vee)\to \XX$, 
we see that $\End_{\YY_1}(\pi_1^\ast T)$ is an NCCR of $(S_1)_0$ by \cite[Corollary 1.20]{HIMO} or Theorem~\ref{NCCRfromGL}, 
where $S_1=R[Y]$ is an $\LL$-graded ring defined by $\deg Y=-\vec{x}_1-\vec{x}_2$.
Note that $(S_1)_0$ is the Veronese subring $R^{(\vec{\om})}$ in the sense of \cite{HIMO}, and 
\[(S_1)_0=k[X_1^{p_1+p_2}Y^{p_2}, X_2^{p_1+p_2}Y^{p_1}, X_1X_2Y]\cong k[A,B,C]/(AB-C^{p_1+p_2})\] is 
the Kleinian $A_{p_1+p_2-1}$ singularity, see also \cite[Proposition~8.4]{GL2}. 
\item[(2)] 
For the vector bundle $E_2=\O_\XX(-\vec{x}_1)\oplus \O_\XX(-\vec{x}_2)$ and its total space $\pi_2\colon \YY_2\coloneqq\Spec (\Sym E_2^\vee)\to \XX$, 
we see that $\End_{\YY_2}(\pi_2^\ast T)$ is an NCCR of $(S_2)_0$ by Theorem~\ref{NCCRfromGL}, 
where $S_2=R[Y_1, Y_2]$ is an $\LL$-graded ring defined by $\deg Y_i=-\vec{x}_i$ $(i=1,2)$. 
Note that 
\[
(S_2)_0=k[X_1^{p_1}Y_2^{p_2}, X_2^{p_2}Y_1^{p_1}, X_1Y_1, X_2Y_2]\cong k[A,B,C,D]/(AB-C^{p_1}D^{p_2})
\]
is a toric compound $A_{p_1+p_2-1}$ singularity. 
An NCCR of a (toric) compound $A_{p_1+p_2-1}$ singularity was also discussed in \cite{IW_cAn,Nag,Nak_cAn,Wem_MMP}. 
\end{itemize}
Figure~\ref{fig_canonicalalg} shows the quiver descriptions of the canonical algebra $\End_\XX(T)$ and NCCRs $\End_{\YY_1}(\pi_1^\ast T)$, $\End_{\YY_2}(\pi_2^\ast T)$, respectively. 
Here, we denote only the weights $\vec{x}$ as the vertices of the quiver instead of writing $\O_\XX(\vec{x})$, $\O_{\YY_1}(\vec{x})$, $\O_{\YY_2}(\vec{x})$. 
\end{Ex}

\begin{figure}[H]
\begin{center}\scalebox{0.8}{
\begin{tikzpicture}[sarrow/.style={black, -latex},tarrow/.style={black, latex-}]
\newcommand{\arrowwidth}{0.03cm} 
\newcommand{\distance}{2.5} 

\node (V000) at (0:0) {$0$}; 
\foreach \a/\b/\c in {0/120/240}{
\foreach \v/\d/\n in {001/1/\vec{x}_3,002/2/2\vec{x}_3,003/3/\vec{c}}{\node (V\v) at (\a:\d*\distance) {$\n$};}; 
\foreach \v/\d/\n in {100/1/\vec{x}_1,200/2/\vec{c}}{\node (V\v) at (\b:\d*\distance) {$\n$};}; 
\foreach \v/\d/\n in {010/1/\vec{x}_2,020/2/\vec{c}}{\node (V\v) at (\c:\d*\distance) {$\n$};}; 

\foreach \v/\d/\n in {101/1/\vec{x}_1+\vec{x}_3,102/2/\vec{x}_1+2\vec{x}_3,103/3/\vec{x}_1+\vec{c}}{\path (V100) ++(\a:\d*\distance) node (V\v) {$\n$};}; 
\foreach \v/\d/\n in {201/1/\vec{x}_3+\vec{c},202/2/2\vec{x}_3+\vec{c},203/3/2\vec{c}}{\path (V200) ++(\a:\d*\distance) node (V\v) {$\n$};}; 

\foreach \v/\d/\n in {011/1/\vec{x}_2+\vec{x}_3,012/2/\vec{x}_2+2\vec{x}_3,013/3/\vec{x}_2+\vec{c}}{\path (V010) ++(\a:\d*\distance) node (V\v) {$\n$};}; 
\foreach \v/\d/\n in {021/1/\vec{x}_3+\vec{c},022/2/2\vec{x}_3+\vec{c},023/3/2\vec{c}}{\path (V020) ++(\a:\d*\distance) node (V\v) {$\n$};}; 

\foreach \v/\d/\n in {110/1/\vec{x}_1+\vec{x}_2,210/2/\vec{x}_2+\vec{c}}{\path (V010) ++(\b:\d*\distance) node (V\v) {$\n$};}; 
\foreach \v/\d/\n in {120/1/\vec{x}_1+\vec{c},220/2/2\vec{c}}{\path (V020) ++(\b:\d*\distance) node (V\v) {$\n$};}; 
}

\begin{scope}[transform canvas={yshift=0.1cm}]
\foreach \a/\b in {000/001,001/002,002/003,100/101,101/102,102/103,200/201,201/202,202/203,010/011,011/012,012/013,020/021,021/022,022/023} 
{\draw[sarrow, line width=\arrowwidth] (V\a)--(V\b)node[midway,xshift=0cm,yshift=0.2cm] {\tiny $X_3$}; };  
\end{scope}
\begin{scope}[transform canvas={xshift=0.1cm}]
\foreach \a/\b in {000/100,100/200,001/101,101/201,002/102,102/202,003/103,103/203,010/110,110/210,020/120,120/220} 
{\draw[sarrow, line width=\arrowwidth] (V\a)--(V\b)node[midway,xshift=0.3cm,yshift=0cm] {\tiny $X_1$}; };  
\foreach \a/\b in {000/010,010/020,001/011,011/021,002/012,012/022,003/013,013/023,100/110,110/120,200/210,210/220} 
{\draw[sarrow, line width=\arrowwidth] (V\a)--(V\b)node[midway,xshift=0.3cm,yshift=0cm] {\tiny $X_2$}; }; 
\end{scope}

\begin{scope}[transform canvas={yshift=-0.1cm}]
\foreach \a/\b in {000/001,001/002,002/003,100/101,101/102,102/103,200/201,201/202,202/203,010/011,011/012,012/013,020/021,021/022,022/023} 
{\draw[tarrow, line width=\arrowwidth,red] (V\a)--(V\b)node[midway,xshift=0cm,yshift=-0.2cm,red] {\tiny $Y_3$}; };  
\end{scope}

\begin{scope}[transform canvas={xshift=-0.1cm}]
\foreach \a/\b in {000/100,100/200,001/101,101/201,002/102,102/202,003/103,103/203,010/110,110/210,020/120,120/220} 
{\draw[tarrow, line width=\arrowwidth,red] (V\a)--(V\b)node[midway,xshift=-0.3cm,yshift=0cm,red] {\tiny $Y_1$}; };  
\foreach \a/\b in {000/010,010/020,001/011,011/021,002/012,012/022,003/013,013/023,100/110,110/120,200/210,210/220} 
{\draw[tarrow, line width=\arrowwidth,red] (V\a)--(V\b)node[midway,xshift=-0.3cm,yshift=0cm,red] {\tiny $Y_2$}; }; 
\end{scope}

\draw [line width=\arrowwidth,sarrow, red] (V110) to [in=195,out=165,loop] (); 
\path (V110) ++(160:1.3cm) node[red] {\tiny $X_3Y_3$}; 

\foreach \v in {101,102} {\draw [line width=\arrowwidth,sarrow, red] (V\v) to [in=45,out=90,loop] ();} 
\foreach \a in {101,102} {\path (V\a) ++(60:1.2cm) node[red] {\tiny $X_2Y_2$}; }
\foreach \v in {011,012} {\draw [line width=\arrowwidth,sarrow, red] (V\v) to [in=-45,out=-90,loop] ();} 
\foreach \a in {011,012} {\path (V\a) ++(-60:1.2cm) node[red] {\tiny $X_1Y_1$}; }

\end{tikzpicture}}
\caption{The quiver description of $\End_{\YY}(\pi^\ast T)$ for $p=(p_1,p_2,p_3)=(2,2,3)$}
\label{fig_quiver_223}
\end{center}
\end{figure}
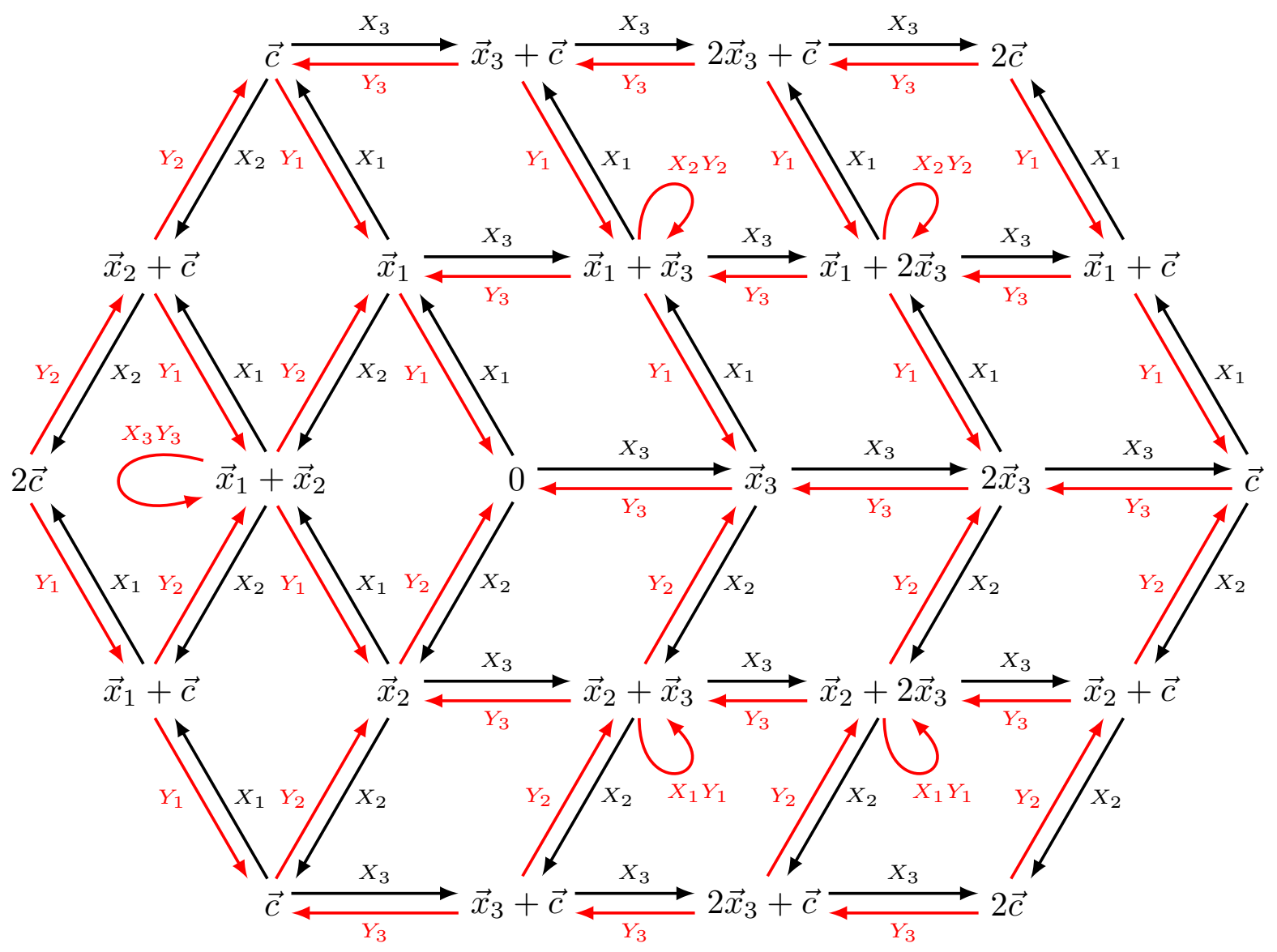

\begin{Ex}
For $n=d+1=3$ and $p=(p_1,p_2,p_3)$, we have the ring $R=k[X_1,X_2,X_3]$ graded by 
$\LL=\Z\vec{x}_1\oplus\Z\vec{x}_2\oplus\Z\vec{x}_3/(p_1\vec{x}_1=p_2\vec{x}_2=p_3\vec{x}_3\eqqcolon\vec{c})$ 
by setting $\deg X_i=\vec{x}_i$ $(i=1,2,3)$. 
In this case, the GL projective space $\XX$  has a tilting bundle with the form $T=\bigoplus_{0\leq\vec{x}\leq 2\vec{c}}\O_\XX(\vec{x})$ by Theorem~\ref{GLtilting}.
Let us consider the vector bundle $E=\O_\XX(-\vec{x}_1)\oplus \O_\XX(-\vec{x}_2)\oplus \O_\XX(-\vec{x}_3)$ and 
its total space $\pi\colon \YY\coloneqq\Spec (\Sym E^\vee)\to \XX$. 
We see that $\End_{\YY}(\pi^\ast T)$ is an NCCR of $S_0$ by Theorem~\ref{NCCRfromGL}, 
where $S\coloneqq R[Y_1, Y_2, Y_3]$ is an $\LL$-graded ring defined by $\deg Y_i=-\vec{x}_i$ $(i=1,2, 3)$. 
Here, $S_0$ is a five-dimensional Gorenstein toric ring with the form: 
\[
S_0=k[X_1Y_1, X_2Y_2, X_3Y_3, X_1^{p_1}Y_2^{p_2}, X_2^{p_2}Y_1^{p_1}, 
X_1^{p_1}Y_3^{p_3}, X_3^{p_3}Y_1^{p_1}, X_2^{p_2}Y_3^{p_3}, X_3^{p_3}Y_2^{p_2} ]. 
\]
For example, Figure~\ref{fig_quiver_223} is the quiver description of $\End_{\YY}(\pi^\ast T)$ for $p=(p_1,p_2,p_3)=(2,2,3)$. 
Here, we denote only the weights $\vec{x}$ as the vertices of the quiver instead of writing $\O_{\YY}(\vec{x})$, 
and the vertices labeled by the same weight are identified. 
Note that we have the quiver description of $\End_\XX(T)$ by removing red arrows from this quiver. 
\end{Ex}

\subsection*{Acknowledgements} 
The authors would like to thank Wahei Hara for pointing out the importance of Corollary~\ref{CY-tilting}(3).
The authors would also like to thank Osamu Iyama and Ryu Tomonaga for stimulating discussions. 

The authors are supported by the World Premier International Research Center Initiative (WPI initiative), MEXT, Japan. 
The first author is supported by JSPS KAKENHI Grant Numbers JP22KJ0737 and JP25K17233.
The second author is supported by JSPS Grant-in-Aid for Scientific Research (C) JP24K06698.

\subsection*{AI disclosure } 
The mathematical ideas were developed by the authors without the use of AI. 
ChatGPT 5.6 Sol and ChatGPT 6 Astra were used to assist with language editing, literature searches including some basic results on stacks, and to check mathematical correctness. 
It pointed out an inaccuracy in the statement of Lemma \ref{codim}(2) (which was previously claimed without any conditions on the idempotent) by giving a (trivial) counter-example of a disconnected ring. It also suggested a proof that this condition is satisfied in the setting of our application, which we adopted as the proof of Corollary~\ref{CY-tilting}(2).


\end{document}